\documentclass[11pt,a4paper]{article}
\usepackage[T1]{fontenc}
\usepackage[utf8]{inputenc}
\usepackage{lmodern}
\usepackage[a4paper,left=26mm,right=26mm,top=25mm,bottom=26mm]{geometry}
\usepackage{amsmath,amssymb,amsthm,mathtools,mathrsfs,bm}
\usepackage{graphicx,booktabs,longtable,array,calc,etoolbox}
\usepackage{fvextra}
\usepackage{microtype}
\usepackage{xcolor}
\definecolor{ReferenceBlue}{RGB}{32,67,117}
\usepackage[unicode,colorlinks=true,linkcolor=ReferenceBlue,citecolor=ReferenceBlue,urlcolor=ReferenceBlue,linktoc=all,bookmarksdepth=2]{hyperref}
\usepackage{bookmark}
\usepackage{xurl}
\usepackage{enumitem}
\setlist{itemsep=2pt,topsep=4pt}
\DeclareUnicodeCharacter{2212}{-}
\DeclareUnicodeCharacter{2264}{\ensuremath{\le}}
\DeclareUnicodeCharacter{2265}{\ensuremath{\ge}}
\DeclareUnicodeCharacter{221E}{\ensuremath{\infty}}
\DeclareUnicodeCharacter{03C0}{\ensuremath{\pi}}
\DeclareUnicodeCharacter{03C6}{\ensuremath{\varphi}}
\DeclareUnicodeCharacter{03B2}{\ensuremath{\beta}}
\DeclareUnicodeCharacter{03B7}{\ensuremath{\eta}}
\DeclareUnicodeCharacter{03B4}{\ensuremath{\delta}}
\DeclareUnicodeCharacter{2192}{\ensuremath{\to}}
\DeclareUnicodeCharacter{00D7}{\ensuremath{\times}}
\DeclareUnicodeCharacter{2208}{\ensuremath{\in}}
\DeclareUnicodeCharacter{2248}{\ensuremath{\approx}}
\DeclareUnicodeCharacter{2218}{\ensuremath{\circ}}
\DeclareUnicodeCharacter{2282}{\ensuremath{\subset}}
\DeclareUnicodeCharacter{2113}{\ensuremath{\ell}}
\DeclareUnicodeCharacter{2011}{-}
\DeclareUnicodeCharacter{2009}{\,}
\usepackage{tikz}
\usetikzlibrary{calc,arrows.meta}
\usepackage[nottoc]{tocbibind}
\begin{document}
\title{The Chen–Yang volume conjecture for long integral fillings of fundamental shadow links}
\author{Ce Shen\\[2pt]\parbox{0.90\textwidth}{\centering\small Beijing Institute of Mathematical Sciences and Applications (BIMSA), Beijing, China}\\[2pt]{\small\href{mailto:shence@bimsa.cn}{\texttt{shence@bimsa.cn}}}}
\hypersetup{pdftitle={The Chen–Yang volume conjecture for long integral fillings of fundamental shadow links},pdfauthor={Ce Shen}}
\date{10 September 2026}
\maketitle
\begin{abstract}
We prove that the Turaev--Viro invariants recover the hyperbolic volume of every sufficiently long integral Dehn filling of a fixed fundamental shadow link complement, establishing the Chen--Yang volume conjecture for these manifolds. More precisely, we obtain a complete asymptotic expansion of the corresponding $SO(3)$ Reshetikhin--Turaev invariant as the odd level tends to infinity. Its exponential term records the volume and Chern--Simons invariant, and we identify the leading coefficient explicitly in terms of adjoint twisted Reidemeister torsion. The main difficulty is cancellation in the signed surgery sum. We control it by proving an exponentially accurate reflection estimate for quantum $6j$-symbols, which identifies the surviving Fourier terms. We also determine the leading effect of fixed even colors on the core curves of the filling solid tori.
\end{abstract}

\begingroup
\makeatletter
\renewcommand{\l@section}[2]{\vskip0.3em\@dottedtocline{1}{0em}{0em}{\bfseries #1}{\bfseries #2}}
\makeatother
{\small\setlength{\parskip}{0pt}\tableofcontents}
\endgroup
\clearpage

\section{1. Introduction}\label{sec-1}

Quantum invariants are defined algebraically, yet their growth at large levels can encode geometry. The Chen--Yang volume conjecture predicts that the leading exponential growth of the Turaev--Viro invariant recovers hyperbolic volume. After Dehn filling, proving this relation requires control of cancellation in the oscillating surgery sum.

Let \(M\) be a closed oriented hyperbolic \(3\)-manifold. For odd \(r\), let \(TV_r(M)\) be its Turaev--Viro invariant at \(q=e^{2\pi i/r}\), in the \(SO(3)\) convention specified in \hyperref[sec-1-3-4]{Section 1.3.4} \cite{TV}, \cite[Section 2.2]{DKY}. In this convention, the conjecture is

\[\lim_{\substack{r\to\infty\\r\text{ odd}}}\frac{2\pi}{r}\log TV_r(M)=\operatorname{Vol}(M).\]

The manifold \(M\) is fixed as \(r\) tends to infinity. Chen and Yang formulated this conjecture in \cite{CY}, following the volume conjecture for knot invariants of Kashaev \cite{K} and Murakami--Murakami \cite{MM}.

We prove the conjecture for sufficiently long integral fillings of any fixed fundamental shadow link complement. Introduced by Costantino--Thurston \cite{CT}, these complements decompose into standard blocks with explicit hyperbolic geometry and quantum weights \cite{C}, \cite[Section 2.4]{WY}. Dehn filling attaches a solid torus to each boundary torus. An integral coefficient \(p_i\) identifies its meridian with \(p_i\mu_i+\lambda_i\), where \((\mu_i,\lambda_i)\) is the chosen meridian--longitude pair.

Our proof uses the \(SO(3)\) Reshetikhin--Turaev invariant \(Z_r(M)\) \cite{W,RT}, \cite[Section 2.1]{WY}, with \(TV_r(M)=|Z_r(M)|^2\). Its surgery formula is a signed sum over colors, with a quantum \(6j\)-symbol for each block and an oscillating filling factor. An exponentially accurate reflection estimate controls cancellation and identifies the surviving Fourier terms.

\subsection{1.1 Main results}\label{sec-1-1}

The block descriptions make fundamental shadow links a useful setting for comparing geometry with quantum asymptotics. Under long fillings, the geometric critical point stays near the complete structure, where both descriptions are explicit. The theorem identifies the surviving exponential term with hyperbolic volume and proves that its leading coefficient is nonzero.

For a closed hyperbolic manifold \(M\), let \(\rho_{\rm geom}\) be its holonomy representation and put \(\mathcal T_M=\operatorname{Tor}(M;\operatorname{Ad}\rho_{\rm geom})\). This is the nonzero adjoint twisted Reidemeister torsion, defined up to sign in the convention of \cite{WY-T}.

\leavevmode\phantomsection\label{stmt-main-theorem}\textbf{Main theorem.} Let \(N\) be a fixed marked fundamental shadow link complement with \(k\) boundary tori and integers \(\iota_1,\ldots,\iota_k\) determined by the mutations in \hyperref[sec-1-3-2]{Section 1.3.2}. There is a constant \(P_0(N)>0\) such that every fixed integral vector \(\mathbf p=(p_1,\ldots,p_k)\) satisfying

\[\min_i|p_i+\iota_i/2|>P_0(N)\]

defines a closed hyperbolic filling \(M=N(\mathbf p)\) for which

\[|Z_r(M)|=\frac{\exp\!\left(r\operatorname{Vol}(M)/(4\pi)\right)}{\sqrt{|\mathcal T_M|}}(1+O_{N,\mathbf p}(r^{-1}))\]

as \(r\) tends to infinity through all odd integers. In particular,

\[\lim_{\substack{r\to\infty\\r\text{ odd}}}\frac{2\pi}{r}\log TV_r(M)=\operatorname{Vol}(M).\]

More precisely, we establish a complete asymptotic expansion and identify its leading coefficient explicitly. Formula \eqref{eq-4-1} records the Chern--Simons term and the exact surgery phases. After these factors are extracted, the remaining expansion is in powers of \(r^{-1}\). \hyperref[sec-3-2]{Section 3.2} gives a formula for every coefficient and works out the first correction \(c_{\mathbf p,1}\).

There is also a version with colored core curves. For fixed even colors \(m_i\) on the core curves \(\gamma_i\) of the filling solid tori, write \(Z_r(M,\gamma;\mathbf m)\) for the relative Reshetikhin--Turaev invariant of the resulting colored pair, where \(\gamma\) is their union and \(\mathbf m=(m_1,\ldots,m_k)\) is the color vector, using the framings in \hyperref[sec-4-3]{Section 4.3}. The symbol \(\operatorname{Sym}^{m_i}\) denotes the symmetric-power representation of dimension \(m_i+1\), and \(\widetilde\rho_{\rm geom}\) is an \(\mathrm{SL}_2(\mathbb C)\) lift of the holonomy representation. Then

\[\frac{Z_r(M,\gamma;\mathbf m)}{Z_r(M)}
=\prod_i\operatorname{tr}\!\left(\operatorname{Sym}^{m_i}\widetilde\rho_{\rm geom}(\gamma_i)\right)+O_{N,\mathbf p,\mathbf m}(r^{-1}).\]

Evenness of \(m_i\) makes this character independent of the choice of lift \(\widetilde\rho_{\rm geom}\) to \(\mathrm{SL}_2(\mathbb C)\). The core framing is specified in \hyperref[sec-4-3]{Section 4.3}.

The filling coefficients may have mixed signs and unrelated magnitudes. The threshold depends on the marked link complement; the filling vector is then held fixed as the level tends to infinity. The theorem makes no assertion about short fillings or a joint limit in the level and filling coefficients. It also uses the integral surgery formula. Rational fillings require a separate analysis of their surgery factors.

\textbf{The reflection estimate behind the theorem.} Reversing any selection of the signed dihedral-angle coordinates leaves the complete asymptotic expansion of a block integral unchanged. We prove that the difference is exponentially small, with a uniform positive gap in its exponential rate. This precision is what allows reflection to control the signed surgery sum after its larger terms have canceled. \hyperref[stmt-theorem-2-2]{Theorem 2.2} gives the local statement, and \hyperref[sec-B]{Appendix B} extends it to compact subsets of the reflection domain defined in \hyperref[sec-A]{Appendix A} and compares the integrals with the finite Racah sums.

\textbf{A separate application.} \hyperref[sec-E]{Appendix E} studies a one-edge state sum in which the six labels of each block agree and the outer color is restricted to a fixed central interval. For every fixed number of blocks \(t\ge2\), we prove that the two terms from the sine quantum dimension have the same complete expansion and cancel to exponential accuracy. For fixed \(t\ge8192\), we also determine the smaller surviving exponential rate and its leading coefficient. This is a result about the specified central sum. A manifold application would additionally require an exact invariant formula and estimates for the omitted outer colors.

\subsection{1.2 Relation to earlier work}\label{sec-1-2}

Belletti--Detcherry--Kalfagianni--Yang \cite{BDKY} proved the Turaev--Viro volume conjecture for fundamental shadow link complements. The behavior under Dehn filling is more delicate because the surgery formula introduces oscillation. Ohtsuki \cite{O} proved the asymptotic expansion and volume conjectures for integral surgeries on the figure-eight knot. Wong--Yang \cite{WY-F8} extended the result to all hyperbolic Dehn surgeries on that knot, using quantum dilogarithms, Poisson summation, and saddle-point analysis.

For fillings of fundamental shadow link complements, Pandey--Wong \cite{PW} obtained asymptotics for sufficiently long rational fillings with cone angles less than \(\pi\). They also discussed the smooth filling value \(2\pi\), at which the singularity along a core curve disappears. The leading-coefficient theorem of Wong--Yang \cite[Theorem 1.2]{WY-A} has a small-cone-angle hypothesis as well. The \hyperref[stmt-main-theorem]{main theorem} here treats the smooth value for sufficiently long integral fillings of each fixed marked link complement. We prove the needed estimates directly for the signed surgery sum at that value.

The reflection property of tetrahedral asymptotics was proposed by Chen--Murakami \cite[Conjecture 3 and Remark 3]{CM}. They conjectured that the complete expansion of their continuous quantum \(6j\)-symbol is invariant under independent reversals of the dihedral-angle signs, and explained its relevance to cancellation. We prove an exponentially accurate reflection estimate for the midpoint interpolation defined in \hyperref[sec-A]{Appendix A}. A comparison with their interpolation on its entire domain is not needed for our argument.

The exact reflection identity for the real-parameter double sine integral is due to Liu--Ming--Sun--Wu--Yang \cite[Proposition 3.4]{LMSWY}; Meng--Yang \cite[Proposition 2.3]{MY} extended it to complex parameter. We use the identity for positive real parameter. The proof compares the saddle expansions of this exactly symmetric integral and our quantum-factorial integral. They have the same coefficients, up to alternating signs. A factorial bound for the remainders then gives an exponentially small reflection error. This last step is necessary for the subsequent signed sum.

The geometric ingredients have earlier origins. Costantino \cite{C} established hyperbolic tetrahedral asymptotics. Chen--Murakami \cite[Theorem 2]{CM} identified the Gram determinant in the leading coefficient, and Belletti--Yang \cite{BY} studied more general admissible tetrahedral weights. We use the absolute growth bounds of \cite{BDKY}, the strict maximum estimates in \cite{PW}, and the signed state-sum and classical potential conventions of Wong--Yang \cite{WY}. Their Gram-matrix and filling formulas for adjoint twisted Reidemeister torsion \cite{WY-T} identify our leading coefficient after the normalization calculation in \hyperref[sec-D-2]{Appendix D.2}.

Cancellation in other filling problems was studied by Ge--Meng--Wang--Yang \cite[Section 4.3]{GMWY}. Factorial asymptotics of the quantum dilogarithm also appear in the resurgence and Borel-summation results of Garoufalidis--Kashaev \cite{GK}. Here the required estimates are proved for the specific finite sums in the surgery formula. The error comparison in \hyperref[sec-2-3]{Section 2.3} states how much exponential accuracy is needed, and \hyperref[stmt-theorem-3-1]{Theorem 3.1} isolates the Fourier argument from its topological application.

\subsection{1.3 Geometric and quantum preliminaries}\label{sec-1-3}

We first recall hyperbolic volume and the geometric gluing operations used in the block construction. We then describe the link complement and its peripheral curves, followed by holonomy, torsion, and the quantum weights. \hyperref[sec-2]{Section 2} introduces the continuous coordinates and states the analytic estimates used to pass from finite weights to integrals.

\subsubsection{1.3.1 Hyperbolic manifolds, polyhedra, and volume}\label{sec-1-3-1}

An oriented hyperbolic \(3\)-manifold carries a Riemannian metric of constant sectional curvature \(-1\): every point has a neighborhood isometric to an open subset of hyperbolic space \(\mathbb H^3\). We consider complete metrics of finite volume. Their volume is the integral of the hyperbolic volume element \(d\operatorname{vol}\),

\[\operatorname{Vol}(M)=\int_M d\operatorname{vol}.\]

Fixing the curvature at \(-1\) fixes the scale. Mostow--Prasad rigidity says that a complete finite-volume hyperbolic metric on a given \(3\)-manifold without boundary is unique up to isometry. Its volume is therefore a topological invariant \cite[Theorem 5.7.2]{TH}.

A closed manifold is compact and has no boundary. A noncompact complete finite-volume hyperbolic \(3\)-manifold has finitely many \emph{cusps}: each end is homeomorphic to \(T^2\times[0,\infty)\), where \(T^2\) is a torus, and has Euclidean torus cross sections that shrink toward infinity. Such an end is infinitely long but has finite volume. Cutting off the ends gives a compact manifold with torus boundary. When \(N\) denotes this compact link complement, \(\operatorname{Vol}(N)\) means the volume of its complete hyperbolic interior, including the cusp ends \cite[Section 5.8]{TH}.

\textbf{Polyhedral pieces and gluing.} A hyperbolic polyhedron is bounded by faces lying in hyperbolic planes. An \emph{ideal vertex} lies on the sphere at infinity of \(\mathbb H^3\) and is omitted from the polyhedron itself. Thus an ideal triangle has three vertices at infinity, although its sides are geodesics. Ideal polyhedra such as the tetrahedron and octahedron have finite volume despite their noncompact ends.

To assemble hyperbolic pieces, identify selected faces by isometries. The identifications must fit consistently around the edges and ends: a smooth interior edge has total dihedral angle \(2\pi\), and the ends must satisfy the completeness conditions. A boundary surface is \emph{totally geodesic} if it is locally a hyperbolic plane. Two pieces can be joined along isometric totally geodesic boundary surfaces; in particular, \emph{doubling} takes two copies and identifies their corresponding boundary points \cite[Chapter 3]{TH}.

Volume is additive in such a decomposition. If the polyhedral pieces \(\mathcal P_1,\ldots,\mathcal P_\ell\) have disjoint interiors and are joined along faces, then

\[\operatorname{Vol}(M)=\sum_{j=1}^{\ell}\operatorname{Vol}(\mathcal P_j).\]

The common faces have zero \(3\)-dimensional volume, so they contribute no extra term. The same principle applies to cutting and gluing along totally geodesic surfaces, and a double has twice the volume of one copy.

\textbf{The regular ideal octahedron.} This polyhedron has six ideal vertices, eight ideal triangular faces, and dihedral angle \(\pi/2\) at each of its twelve edges. We denote it by \(\mathcal O_8\) and write

\[v_8=\operatorname{Vol}(\mathcal O_8)=3.66386\ldots.\]

This is the basic volume unit for the shadow-link complements considered below \cite[Chapters 3 and 7]{TH}. The next subsection constructs the \(D\)-block topologically and then realizes its two halves as regular ideal octahedra. Volume additivity will give \(2v_8\) for one block and \(2cv_8\) for a complement made from \(c\) blocks.

\subsubsection{1.3.2 Hyperbolic geometry of fundamental shadow links}\label{sec-1-3-2}

We describe the construction from its individual pieces. First we form a \(D\)-block from two truncated tetrahedra. We then join two such blocks in an explicit example, before explaining the general assembly, the peripheral markings, and the closed hyperbolic manifolds obtained by Dehn filling.

\textbf{1. The truncated tetrahedron.} Start with a tetrahedron and cut off a small neighborhood of each of its four vertices. The resulting polyhedron, denoted by \(\Delta\), is still a topological \(3\)-ball. It has four triangular faces, one at each removed vertex, and four hexagonal faces, one from each original face. A hexagon has three sides inherited from the original tetrahedron, alternating with three sides created by the cuts.

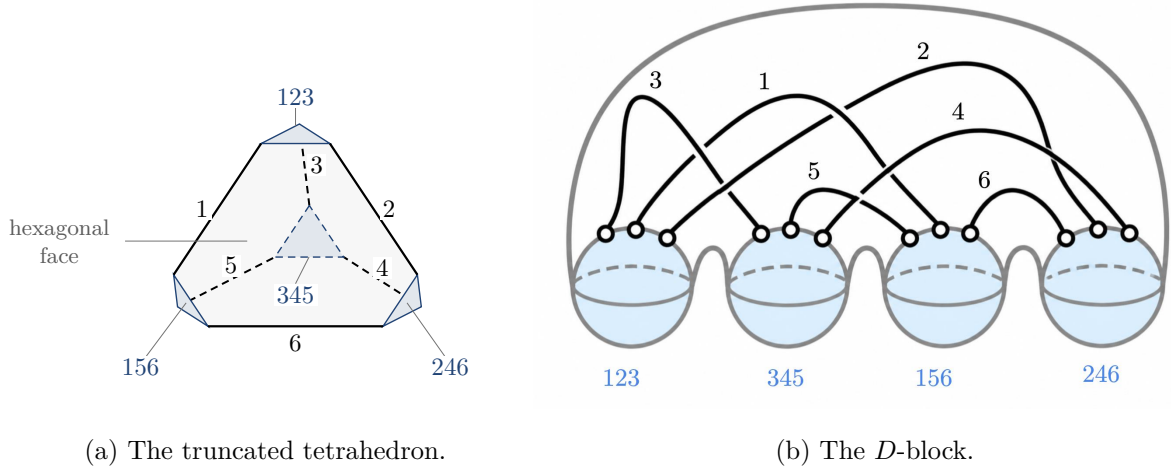
\begin{figure}[!htbp]
\centering
\begin{minipage}[b]{0.40\linewidth}
\centering
\resizebox{\linewidth}{!}{%
\begin{tikzpicture}[x=1cm,y=1cm,font=\small,
  edge/.style={line width=0.9pt},
  hidden edge/.style={line width=0.8pt,densely dashed},
  trunc/.style={draw=ReferenceBlue,fill=ReferenceBlue!13,line width=0.55pt},
  lab/.style={fill=white,inner sep=1.2pt},
  guide/.style={black!55,line width=0.4pt}]
\path[use as bounding box] (-3.8,-2.35) rectangle (2.8,2.35);
\coordinate (A) at (0,2.2);
\coordinate (B) at (-2.2,-1.1);
\coordinate (C) at (2.2,-1.1);
\coordinate (D) at (0.25,0.15);
\foreach \v/\w in {A/B,A/C,A/D,B/A,B/C,B/D,C/A,C/B,C/D,D/A,D/B,D/C}
  \coordinate (\v\w) at ($(\v)!0.22!(\w)$);
\fill[black!3] (AB)--(BA)--(BC)--(CB)--(CA)--(AC)--cycle;
\draw[hidden edge] (AD)--(DA);
\draw[hidden edge] (BD)--(DB);
\draw[hidden edge] (CD)--(DC);
\draw[trunc,densely dashed] (DA)--(DB)--(DC)--cycle;
\draw[trunc] (AB)--(AC)--(AD)--cycle;
\draw[trunc] (BA)--(BC)--(BD)--cycle;
\draw[trunc] (CA)--(CB)--(CD)--cycle;
\draw[edge] (AB)--(BA) node[midway,lab,left=2pt] {$1$};
\draw[edge] (AC)--(CA) node[midway,lab,right=2pt] {$2$};
\draw[edge] (BC)--(CB) node[midway,lab,below=2pt] {$6$};
\node[lab,right=1pt] at ($(AD)!0.50!(DA)$) {$3$};
\node[lab,above=1pt] at ($(BD)!0.55!(DB)$) {$5$};
\node[lab,above=1pt] at ($(CD)!0.50!(DC)$) {$4$};
\draw[guide] (0,1.72)--(0,2.02);
\node[text=ReferenceBlue,lab] at (0,2.13) {$123$};
\draw[guide] (-1.54,-0.77)--(-2.1,-1.52);
\node[text=ReferenceBlue,lab] at (-2.18,-1.68) {$156$};
\draw[guide] (1.57,-0.77)--(2.10,-1.52);
\node[text=ReferenceBlue,lab] at (2.18,-1.68) {$246$};
\draw[guide] (0.19,-0.125)--(0,-0.50);
\node[text=ReferenceBlue,lab] at (0,-0.66) {$345$};
\draw[guide] (-2.40,0.10)--(-0.68,0.10);
\node[text=black!65,font=\footnotesize,align=center,anchor=east] at (-2.48,0.10) {hexagonal\\face};

\end{tikzpicture}
}
\par\smallskip
{\small (a) The truncated tetrahedron.}
\end{minipage}\hfill
\begin{minipage}[b]{0.58\linewidth}
\centering
\includegraphics[width=\linewidth]{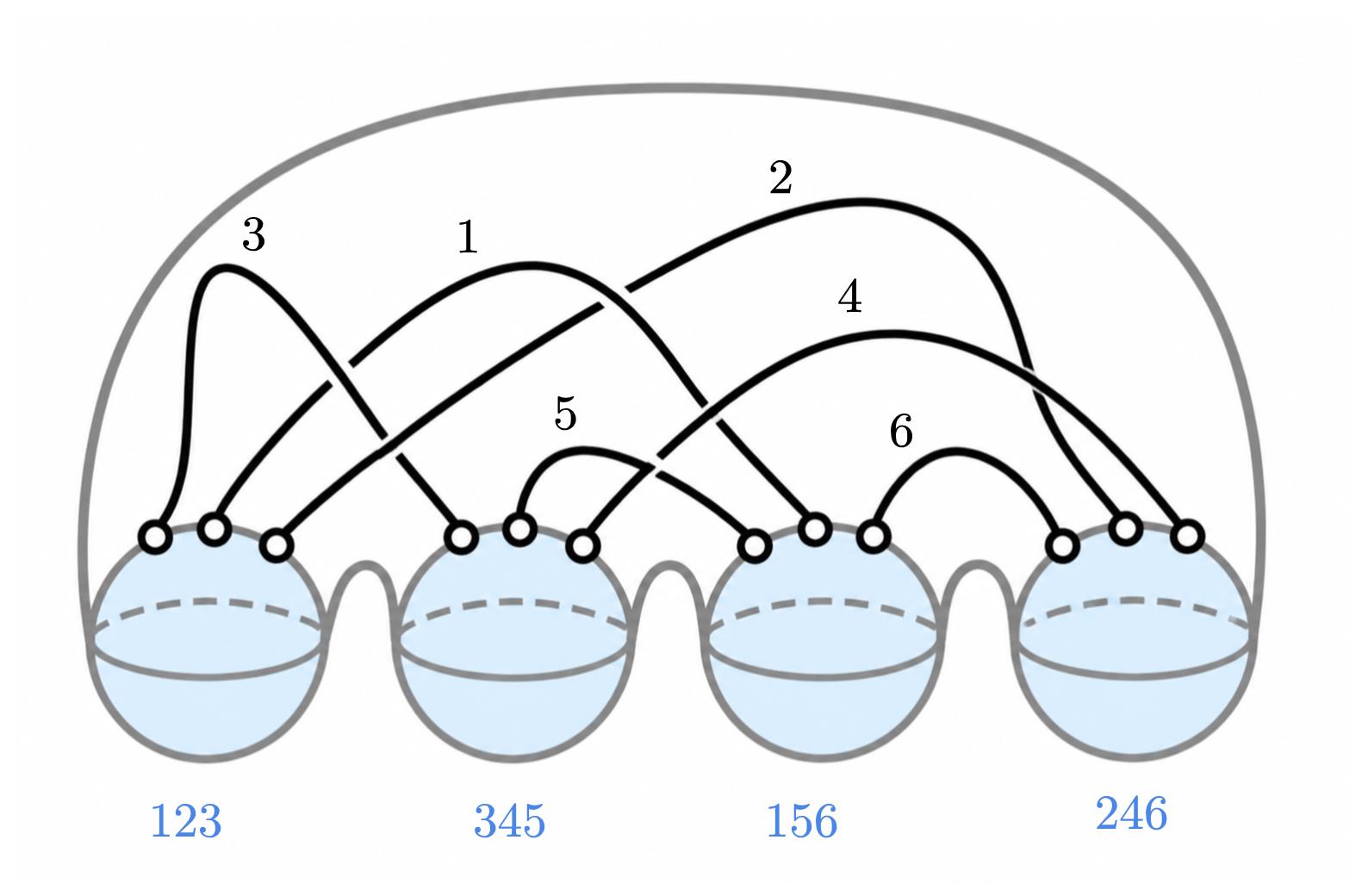}
\par\smallskip
{\small (b) The $D$-block.}
\end{minipage}
\caption{From the truncated tetrahedron to the $D$-block. (a) The six original edges are labeled $1,\ldots,6$, and each shaded truncation triangle is labeled by its three incident edges; the other four faces are hexagons. Dashed lines indicate hidden edges. (b) The four boundary spheres and six edge arcs after the hexagonal-face gluing, with the same labels as in (a). Removing open regular neighborhoods of the black arcs turns the blue spheres into pairs of pants and creates six annuli, giving the $D$-block $B$.}
\label{fig:block-construction}
\end{figure}

Figure \ref{fig:block-construction}(a) shows these two kinds of faces. The shaded triangles are the \emph{triangles of truncation}. The six numbered segments are the remaining portions of the original edges; each lies between two hexagonal faces and has its endpoints on two truncation triangles. We call these six segments the edges of \(\Delta\). The twelve sides created by truncation lie between a triangle and a hexagon and carry no edge labels in our notation.

The label of a truncation triangle records the three numbered edges ending on it. Thus the four triangles are labeled \(123,156,246,345\). For example, edge \(1\) runs from triangle \(123\) to triangle \(156\), since these are precisely the two labels containing \(1\). Every edge label occurs in exactly two triangle labels. These incidences will determine how the link components pass through a block. The figure specifies the topology and labeling; a hyperbolic metric will be introduced after the gluing construction.

\textbf{2. From two tetrahedra to one \(D\)-block.} Take two copies of \(\Delta\), with opposite orientations and identical edge labels. Identify each hexagonal face in one copy with the corresponding hexagonal face in the other, matching all its labeled edges. At this stage the triangular faces remain unglued in their interiors. Each triangle is joined to its counterpart along all three sides, so the pair becomes a sphere. The resulting manifold, denoted by \(\widehat B\), therefore has four spherical boundary components. Topologically, it is a \(3\)-sphere with the interiors of four disjoint balls removed.

The two copies of a numbered edge are identified in this gluing. Their image is a single arc in \(\widehat B\), with its endpoints on two of the boundary spheres. Figure \ref{fig:block-construction}(b) shows these six arcs and the four boundary spheres, ordered from left to right as \(123,345,156,246\). For example, the arc labeled \(1\) connects the spheres labeled \(123\) and \(156\). Each sphere meets exactly three arcs, at the three vertices of its doubled triangle.

Next remove an open regular neighborhood of these six arcs, including small disks around their endpoints on the boundary spheres. Call the remaining compact manifold \(B\), the \emph{\(D\)-block}. Its boundary is divided into four pairs of pants and six annuli. Here a pair of pants is a sphere with three open disks removed. The pair of pants labeled \(123\) has boundary circles labeled \(1,2,3\); the other three are labeled in the same way. The annulus surrounding edge \(1\) joins circle \(1\) on the pair of pants \(123\) to circle \(1\) on the pair of pants \(156\). Each of the other five annuli similarly follows one edge. Thus each annulus has two boundary circles, and every boundary circle of a pair of pants meets exactly one annulus.

The four pairs of pants are the surfaces along which blocks will be glued. The six annuli are the pieces from which the boundary tori of the link complement will be assembled. In the complete hyperbolic metric, the pairs of pants extend to thrice-punctured spheres, and the annuli are cross sections of annular cusps. We use the compact description for the topological gluings. This is the doubled-tetrahedron construction of \cite[Section 2.4]{WY}; throughout the paper, a standard block means this \(D\)-block \(B\).

For the complete metric, each tetrahedral half is realized as a regular ideal octahedron. Combinatorially, this is seen by collapsing each of the six numbered edges to a vertex: the four hexagons become triangles, while the four truncation triangles remain triangles. The resulting polyhedron has six ideal vertices and eight triangular faces. The two octahedra are joined along the four faces inherited from the hexagons. The four remaining pairs of ideal triangles form the totally geodesic thrice-punctured spheres. By volume additivity from \hyperref[sec-1-3-1]{Section 1.3.1}, a complete \(D\)-block therefore has volume \(2v_8\) \cite[Section 2.4]{PW}.

\textbf{3. Example: gluing two \(D\)-blocks.} Take a block \(B^+\) and a mirror copy \(B^-\). In both blocks retain the four boundary labels \(123,156,246,345\). For each label, glue the pair of pants in \(B^+\) to the pair of pants with the same label in \(B^-\) by the mirror identification, matching equally labeled boundary circles. All four pairs of pants in each block are used. Denote the resulting compact link complement by \(N_{\rm dbl}\).

To identify its boundary, follow the annulus labeled \(1\) in \(B^+\). At its \(123\) end it meets the annulus labeled \(1\) in \(B^-\), and the same two annuli meet at their \(156\) ends. Their union is a torus. Repeating this argument for each edge label gives six boundary tori. Equivalently, before the arc neighborhoods are removed, the two arcs labeled \(1\) join end to end to form a closed curve \(L_1\). Each of the other edge labels gives one closed curve \(L_i\). No two different edge labels are identified in this example.

The ambient manifold can also be read from the tetrahedra. Glue one truncated tetrahedron to its mirror along their four truncation triangles. This gives an orientable handlebody \(H\). It retracts onto a graph with two vertices, corresponding to the tetrahedra, and four edges, corresponding to the triangle gluings. Hence its genus is \(4-2+1=3\). The six pairs of tetrahedral edges form the curves \(L_1,\ldots,L_6\) on \(\partial H\). Doubling \(H\) along its boundary then gives

\[L=L_1\cup\cdots\cup L_6\subset\#^3(S^2\times S^1).\]

Here \(\#\) denotes connected sum. Regrouping the four tetrahedral copies in this double by their hexagonal-face gluings gives exactly the two blocks \(B^+,B^-\) constructed above. Thus the example uses two \(D\)-blocks, each containing two tetrahedral halves. In the notation of the theorem, the number of blocks is \(c=2\), the number of link components is \(k=6\), and \(\operatorname{Vol}(N_{\rm dbl})=4v_8\).

There are six independent component colors and twelve edge occurrences. Both blocks therefore contribute the same quantum \(6j\)-symbol, whose definition is given in \hyperref[sec-1-3-5]{Section 1.3.5}, so their product is its square. We return to this example after the finite surgery formula in \hyperref[sec-2-2]{Section 2.2}.

\textbf{4. Assemblies with more blocks.} Now take \(c\) copies of \(B\). Pair their \(4c\) pairs of pants and, for each pair, choose a homeomorphism specifying which of the three boundary circles is matched to which. Use maps that reverse the induced boundary orientations, so that the glued manifold is oriented, and require the assembly to be connected. Two boundary surfaces of the same block may also be paired. Once all \(2c\) gluings have been made, the pairs of pants become internal surfaces. The annuli remain on the boundary and join in cycles, each cycle forming a torus.

One can find the components by following the edge arcs before their neighborhoods are removed. Start on an arc in a block, follow it to a boundary sphere, use the prescribed matching to enter the next block, and continue until the path closes. Each closed path is one link component. A path may visit the same block more than once and may use several of its six edges. Consequently, the number \(k\) of link components depends on the boundary matchings, and a single component color may occur at several edge positions of one block.

For completeness, the block count also determines the ambient manifold. Record the assembly by a graph with one vertex for each block and one edge for each pair of glued boundary surfaces. This connected graph has \(c\) vertices and \(2c\) edges. Before removing the link arcs, each block is a \(3\)-sphere with four balls removed. Gluing first along a spanning tree joins the \(c\) pieces into one punctured \(3\)-sphere. There remain \(2c-(c-1)=c+1\) pairs of boundary spheres; identifying each pair adds one \(S^2\times S^1\) summand. Thus the link lies in \(\#^{c+1}(S^2\times S^1)\). Links obtained by this construction are called \emph{fundamental shadow links} \cite[Section 2.4]{WY}.

We write \(N\) for the compact complement obtained by removing an open tubular neighborhood of the resulting link, and \(k\) for its number of boundary tori. At the complete structure, realize the boundary matchings by isometries of the thrice-punctured spheres. The octahedral pieces then give the complete hyperbolic metric on the interior of \(N\). Since there are two regular ideal octahedra per block, \(\operatorname{Vol}(N)=2cv_8\).

\textbf{5. Peripheral marking, mutations, and filling.} On the torus around \(L_i\), choose a meridian \(\mu_i\) bounding a disk in the removed tubular neighborhood, and a longitude \(\lambda_i\) running once along the component, with oriented intersection \(\mu_i\cdot\lambda_i=1\). This pair is the peripheral marking. Integral filling with coefficient \(p_i\) attaches a solid torus whose meridian is identified with \(p_i\mu_i+\lambda_i\); its core curve is denoted by \(\gamma_i\).

For a fixed block presentation, \(\iota_i\) is the sum of the signed contributions along \(L_i\) from mutations, namely the regluings along thrice-punctured spheres, in the convention of \cite[Section 2.4]{WY}. It records their effect on the preferred longitude and enters the surgery formula through the effective coefficient \(P_i=p_i+\iota_i/2\). Replacing \(\lambda_i\) by \(\lambda_i+m_i\mu_i\), with \(m_i\in\mathbb Z\), changes \(\iota_i\) to \(\iota_i+2m_i\) and the coefficient of the same filling slope to \(p_i-m_i\), leaving \(P_i\) unchanged. For the untwisted two-block example with the longitude supplied by \(\partial H\), \(\iota_i=0\) and \(P_i=p_i\).

Throughout, sufficiently long means that every \(|P_i|\) exceeds a threshold for the fixed marked link complement.

\textbf{6. Closed hyperbolic fillings and volume.} Filling all \(k\) boundary tori in this way gives the closed manifold \(M=N(\mathbf p)\), where \(\mathbf p=(p_1,\ldots,p_k)\). Conversely, removing open tubular neighborhoods of the core curves \(\gamma_1,\ldots,\gamma_k\) recovers \(N\) topologically.

These complements provide a general source of closed manifolds: a theorem of Costantino--Thurston states that every closed oriented \(3\)-manifold can be obtained by integral Dehn filling of a suitable fundamental shadow link complement \cite{CT}, \cite[Theorem 2.6]{PW}. The complement may depend on the manifold, and this existence statement gives no lower bound on the filling coefficients. Our \hyperref[stmt-main-theorem]{main theorem} fixes one marked complement \(N\) and treats its sufficiently long integral fillings, with the precise bound stated in \hyperref[sec-1-1]{Section 1.1}.

For sufficiently long slopes, Thurston's hyperbolic Dehn filling theorem constructs the complete hyperbolic metric on \(M\) by deforming the complete cusp metric on the interior of \(N\). The filled cusps close around short geodesics representing the core curves \cite[Theorem 5.8.2]{TH}. This changes the metric on the complement itself, so the volume of the filled manifold cannot be computed by simply adding solid-torus volumes to the original block volumes. If \(N\) is built from \(c\) blocks, then for these hyperbolic fillings,

\[0<\operatorname{Vol}(N(\mathbf p))<\operatorname{Vol}(N)=2cv_8.\]

Moreover, for fixed \(N\) and fixed peripheral markings,

\[\lim_{\min_i|p_i|\to\infty}\operatorname{Vol}(N(\mathbf p))=\operatorname{Vol}(N)=2cv_8.\]

Thus the volumes of long closed fillings approach the volume of the complete link complement from below \cite[Section 1]{NZ}.

In the volume conjecture studied here, we first choose one filling \(M=N(\mathbf p)\) and then let the quantum level \(r\) tend to infinity with \(\mathbf p\) fixed. The target is \(\operatorname{Vol}(M)\), the volume of this closed manifold. The explicit value \(2cv_8\) supplies the starting point for the deformation calculation. \hyperref[sec-1-3-3]{Section 1.3.3} introduces its holonomy coordinates, and \hyperref[sec-4-1]{Section 4.1} identifies the critical value of the phase at the filling solution with the volume and Chern--Simons invariant of \(M\).

\subsubsection{1.3.3 Holonomy and torsion}\label{sec-1-3-3}

Near the complete structure, logarithmic peripheral holonomies \(u_i,v_i\) are chosen continuously from zero; their eigenvalue ratios are \(e^{u_i},e^{v_i}\). In the branch used in \hyperref[sec-4]{Section 4}, smooth filling is expressed by \(p_i u_i+v_i=2\pi i\). The holonomy representation \(\rho_{\rm geom}\) takes values in \(\mathrm{PSL}_2(\mathbb C)\). Its adjoint action on the traceless \(2\times2\) complex matrices defines the twisted chain complex whose Reidemeister torsion is \(\mathcal T_M\). For the closed hyperbolic fillings considered here this complex is acyclic; its torsion is a nonzero complex number defined up to sign, with the convention of \cite[Section 2.1]{WY-T}. Porti \cite{PORTI} develops the corresponding torsion theory for hyperbolic manifolds and their deformations. Thus the absolute coefficient in the \hyperref[stmt-main-theorem]{main theorem} is unambiguous. The notation \(\operatorname{Sym}^m\) means the \((m+1)\)-dimensional symmetric-power representation; for even \(m\), its character is unchanged on replacing an \(\mathrm{SL}_2(\mathbb C)\) lift by its negative.

\subsubsection{1.3.4 Quantum invariants and colors}\label{sec-1-3-4}

The Reshetikhin--Turaev (RT) invariant is evaluated by a surgery sum over link colors. The Turaev--Viro (TV) invariant is a triangulation state sum built from quantum dimensions and tetrahedral weights. We use the \(SO(3)_{r-2}\) fusion category, with simple colors \(0,2,\ldots,r-3\). The Kauffman-bracket surgery construction is developed in \cite{BHMV}; we use the precise \(SO(3)\) conventions of \cite[Section 2.1]{WY} at the following root of unity:

\[q_r=e^{2\pi i/r},\qquad D_r=\frac{\sqrt r}{2\sin(2\pi/r)}.\]

The \(SO(3)\) RT invariant \(Z_r\) has \(Z_r(S^3)=D_r^{-1}\). For every connected closed oriented manifold \(M\), the \(SO(3)\) TV--RT comparison \cite{ROB}, \cite[Section 3.1]{DKY} gives

\[TV_r(M)=|Z_r(M)|^2.\]

Here \(TV_r\) is the even-color state sum of \cite[Definition 2.8]{DKY}. This convention agrees with the \(SO(3)\) surgery theory used in \cite[Section 2.1]{WY}, \cite[Section 3.2]{WY-A}. For comparison with the original formulation in \cite{CY}, the full-color invariant equals \(2^{b_2(M;\mathbb F_2)-1}TV_r(M)\), where \(b_2(M;\mathbb F_2)=\dim_{\mathbb F_2}H_2(M;\mathbb F_2)\) is the second Betti number with coefficients in the two-element field \cite[Theorem 2.9(1)--(2)]{DKY}. This factor is independent of \(r\), so the two conventions give the same volume limit for each fixed \(M\). All TV formulas in this paper use the \(SO(3)\) convention. Squaring the RT magnitude explains the factors \(4\pi\) and \(2\pi\) in the respective volume-growth normalizations.

\textbf{Colors and quantum factorials.} The level \(r\ge5\) is odd. The set of colors is \(I_r=\{0,2,\ldots,r-3\}\), and

\[[n]=\frac{\sin(2\pi n/r)}{\sin(2\pi/r)},\qquad [n]!=\prod_{j=1}^n[j],\qquad d_a=[a+1],\qquad\theta_a=e^{\pi i a(a+2)/r}.\]

We set \([0]!=1\). The color \(a\) labels an edge or link component, \(d_a\) is its quantum dimension, and \(\theta_a\) is the weight of a unit framing twist. In particular the quantum dimensions are signed. A triple \((a,b,c)\) of colors in \(I_r\) is \(r\)-admissible when its triangle inequalities hold and \(a+b+c\le2(r-2)\). We abbreviate \(r\)-admissible to admissible when the level is fixed. Even colors make all factorial arguments below integers.

\subsubsection{\texorpdfstring{1.3.5 Quantum \(6j\)-symbols and surgery weights}{1.3.5 Quantum 6j-symbols and surgery weights}}\label{sec-1-3-5}

The local weight has six edge labels. The face triples in the formulas below are the triples of edges incident to the triangles of truncation in Figure \ref{fig:block-construction}(a). We use these triples and the quadrilateral index sets

\[\mathcal F=\{123,156,246,345\},\qquad
\mathcal Q=\{1245,1346,2356\},\]

where, for example, \(123\) denotes the set \(\{1,2,3\}\). Each quadrilateral set records the four edges in one cycle of the tetrahedral edge graph. A six-tuple is admissible when all four face triples are admissible. For an admissible triple set

\[\Delta(a,b,c)=\sqrt{\frac{[(a+b-c)/2]![(b+c-a)/2]![(c+a-b)/2]!}{[(a+b+c)/2+1]!}},\]

using \(\sqrt{x}=i\sqrt{|x|}\) for negative real \(x\). For an admissible six-tuple, write \(\mathbf a=(a_1,\ldots,a_6)\) and let \(\mathbf a_f\) denote its restriction to the face \(f\). Set \(T_f=\frac12\sum_{i\in f}a_i\) and \(Q_g=\frac12\sum_{i\in g}a_i\), with \(f\in\mathcal F\) and \(g\in\mathcal Q\). The quantum \(6j\)-symbol in the symmetric normalization is

\[\tau_r(\mathbf a)=i^{-\sum_i a_i}\prod_f\Delta(\mathbf a_f)\sum_z\frac{(-1)^z[z+1]!}{\prod_f[z-T_f]!\prod_g[Q_g-z]!}.\]

The nonzero Racah terms have \(\max_fT_f\le z\le\min(r-2,\min_gQ_g)\). This is the convention of \cite[Definition 2.3]{WY}. The integer \(z\) is an internal summation index for a single tetrahedral weight. The outer surgery sum in \eqref{eq-2-1} then sums over the \(k\) component colors. Each block inherits six labels from those components; a component may occur more than once in a block. Thus there is an internal Racah sum for each block, followed by an external color sum with oscillating framing weights.

\textbf{Surgery phase.} The anomaly is \(\kappa_r=\exp[-\pi i((r+1)/4+3/r)]\). The surgery diagram contributes the multiplier \(\kappa_r^{-\sigma}\), where \(\sigma\) is the linking-matrix signature difference defined in \hyperref[sec-2-2]{Section 2.2}. We retain this multiplier and the signs of all quantum dimensions throughout the calculation; they determine the phase in \eqref{eq-4-1}.

\subsection{1.4 Outline of the proof}\label{sec-1-4}

The main proof has three stages. We first reduce the finite surgery sum to Fourier integrals. We then use exponential reflection to cancel the even Fourier modes. Finally, stationary phase evaluates the surviving odd modes, and their critical data give the geometry of the filled manifold.

\begin{center}
\begin{tabular}{@{}>{\raggedright\arraybackslash}p{0.22\linewidth}>{\raggedright\arraybackslash}p{0.44\linewidth}>{\raggedright\arraybackslash}p{0.26\linewidth}@{}}
\toprule
Stage & Mathematical content & Location \\
\midrule
Algebraic reduction & Surgery factors, centered colors, and Poisson summation & \hyperref[sec-2]{Section 2} \\[3pt]
Analytic input & Exponential reflection and the finite-sum error & \hyperref[sec-2-3]{Section 2.3}; proofs in \hyperref[sec-A]{Appendices A}--\hyperref[sec-B]{B} \\[3pt]
Asymptotics and geometry & The first odd modes, expansion coefficients, and the filling equations & \hyperref[sec-3]{Sections 3}--\hyperref[sec-4]{4} \\
\bottomrule
\end{tabular}
\end{center}

\hyperref[sec-2]{Section 2} begins with a coordinate dictionary and a one-component calculation. The exact surgery formula then shows how the block weights are joined and how the filling coefficients enter. \hyperref[stmt-theorem-2-2]{Theorem 2.2} is the analytic input: it allows us to replace these weights by an even function with an exponentially small error. The parity calculation is exact for that even model.

For long fillings, the leading Fourier orbit has every frequency equal to \(1\) or \(-1\). \hyperref[sec-3]{Section 3} gives its expansion and an explicit first correction. \hyperref[sec-4]{Section 4} identifies the critical equation with smooth Dehn filling, its value with volume and Chern--Simons invariant, and its leading amplitude with adjoint twisted Reidemeister torsion. Fixed even colors on the core curves give the stated character ratios.

The proofs of the analytic inputs are collected separately. \hyperref[sec-A]{Appendix A} constructs the midpoint block integral and proves reflection by comparison with an exactly symmetric double sine integral. \hyperref[sec-B]{Appendix B} treats the hyperideal domain and the endpoints of the finite Racah sum. \hyperref[sec-C]{Appendix C} proves the Fourier estimates used in \hyperref[sec-3]{Section 3}, and \hyperref[sec-D]{Appendix D} verifies the algebraic identities and Gram normalization. The independent one-edge application is in \hyperref[sec-E]{Appendix E}, with its effective constants in \hyperref[sec-F]{Appendix F}.

\section{2. From the surgery sum to Fourier integrals}\label{sec-2}

Each block supplies a quantum \(6j\)-symbol, gluing identifies its edge colors with component colors, and filling introduces a quadratic phase for each boundary torus. We first explain the resulting variables in one dimension. We then give the exact finite formula and state the reflection estimate needed to carry its algebraic reduction through with a controlled error.

\subsection{2.1 Coordinates and a one-component warm-up}\label{sec-2-1}

\textbf{Centered coordinates.} For each component color \(a_i\), set

\[h=\frac1r,\qquad y_i=a_i-\frac{r-2}{2},\qquad
x_i=\frac{y_i}{r}=\frac{a_i+1}{r}-\frac12.\]

The even-color lattice has spacing \(2/r\) in these coordinates. Reflection \(x_i\mapsto-x_i\) corresponds to \(a_i\mapsto r-2-a_i\), which is odd when \(r\) is odd and \(a_i\) is even. Thus reflection will be applied to an interpolating function. The shifted color lattice itself is handled by Poisson summation. A function is independently even if it is unchanged when any subset of its coordinate signs is reversed.

For the surgery sum, \(x=(x_1,\ldots,x_k)\) records component colors. For a block, \(x_b\) is the six-tuple inherited from these components, with repetitions when a component occurs more than once. The internal variable of a block integral records its Racah index and is distinct from the peripheral holonomies \(u_i\).

The following dictionary records the geometric meaning of the quantities that will appear. The identifications involving the critical point are proved in \hyperref[sec-4]{Section 4}; \(S_{\mathbf P}\) denotes the reduced phase and \(\mathcal A_{\mathbf p}\) its leading amplitude, defined in \hyperref[sec-3]{Section 3}.

\begin{center}
\begin{tabular}{@{}>{\raggedright\arraybackslash}p{0.24\linewidth}>{\raggedright\arraybackslash}p{0.31\linewidth}>{\raggedright\arraybackslash}p{0.37\linewidth}@{}}
\toprule
Quantum quantity & Continuous quantity & Geometric interpretation \\
\midrule
Even component color $a_i$ & $x_i=(a_i+1)/r-1/2$ & Centered deformation coordinate; $x_i=0$ is the complete limit \\[4pt]
Color at the geometric saddle & $u_i=4\pi i x_i$ & Logarithmic meridian holonomy \\[4pt]
Exponential part of a Fourier integral & $S_{\mathbf P}(x_{\mathbf p})$ & Volume and Chern--Simons invariant, with the branch and phase conventions in \hyperref[sec-4-1]{Section 4.1} \\[4pt]
Block weights and filling factors & $\mathcal A_{\mathbf p}$, including the Gaussian determinant & The reciprocal square root of adjoint twisted Reidemeister torsion, up to a phase \\
\bottomrule
\end{tabular}
\end{center}

\textbf{One filling variable.} To see the reduction before introducing matrices, specialize to \(k=1\) and put \(P=p+\iota/2\). This is a calculation with one outer color; it does not assume that the complement has only one block. Write its even color as \(a=2m\). Then

\[x=\frac{2m+1}{r}-\frac12
=\frac2r\left(m-\frac{r-2}{4}\right),\qquad -\frac12<x<\frac12.\]

Thus Poisson summation uses a lattice of spacing \(2/r\) with an offset determined by \(r\). The sine dimension and the quadratic filling term become

\[\sin\frac{2\pi(a+1)}r=-\sin(2\pi x),\qquad
\frac{a(a+2)}r=rx^2+rx+\frac r4-\frac1r.\]

After the mutation sign is included, the linear term is constant on the even-color lattice. \hyperref[sec-2-4]{Section 2.4} records this constant together with the other surgery phases. The variable part of the filling term is \(e^{i\pi rPx^2}\).

For the continuous even block model justified in \hyperref[sec-2-3]{Section 2.3}, write \(E_r^{\rm sym}(x)=e^{r\Phi(x)}a_r(x)\). Here \(\Phi\) is the sum of the block critical phases and \(a_r\) contains their amplitudes. Apart from the constant surgery phases, the \(n\)th Fourier integral is

\[I_{r,n}=\int_{\mathbb R}\chi_1(x)a_r(x)\sin(2\pi x)
e^{rS_{P,n}(x)}\,dx,\qquad
S_{P,n}(x)=\Phi(x)+i\pi Px^2-i\pi nx,\]

where \(n\in\mathbb Z\). Here a cutoff is a smooth function that equals one in a smaller neighborhood of zero and vanishes outside a larger fixed neighborhood; \(\chi_1\) is chosen even. The function \(S_{P,n}\) in the exponent is the phase, and \(a_r(x)\sin(2\pi x)\) is the amplitude. The saddle-point calculation expands the phase about a point where its derivative vanishes. Its critical equation is

\[\Phi'(x)+2\pi iPx-i\pi n=0.\]

The three terms have different origins: the block geometry, the filling slope, and the Fourier frequency. Reflection pairs \(n\) with \(-n\) and removes the even frequencies. For large \(|P|\), the leading pair is \(n=\pm1\); for \(n=1\), its critical point is close to \(1/(2P)\). \hyperref[sec-4]{Section 4} identifies its equation with \(pu+v=2\pi i\). The same reduction in \(k\) variables replaces \(Px^2\) by \(\sum_iP_ix_i^2\) and the derivative by the gradient.

\subsection{2.2 The exact finite surgery formula}\label{sec-2-2}

The block decomposition in \hyperref[sec-1-3-2]{Section 1.3.2} determines how the local weights enter the invariant. Each standard block contributes one quantum \(6j\)-symbol. When the blocks are glued, edge occurrences belonging to the same link component carry the same color. The resulting component corresponds to a boundary torus of the link complement. Dehn filling pairs that boundary with a solid torus, contributing its quantum dimension and framing factor before we sum over the component color. Thus the surgery formula has one internal Racah index for each block and one outer color for each boundary torus; repeated appearances of a component are included in the labels of its block factors.

Let \(N\) be a fixed fundamental shadow link complement, presented with \(c\) blocks and \(k\) link components. Let \(\iota_i\) be the integers associated with its mutations and \(p_i\) fixed integral surgery coefficients in the preferred markings. Put

\[P_i=p_i+\iota_i/2,\qquad P_{\min}=\min_i|P_i|.\]

The signed fundamental shadow link formula of \cite[Definition 2.3 and Proposition 2.6]{WY}, followed by the surgery pairing, gives

\[
Z_r(N(\mathbf p))=
\varepsilon_{r,\mathbf p}D_r^c
\sum_{\mathbf a}
\prod_i\left[
\frac2{\sqrt r}\sin\frac{2\pi(a_i+1)}r
(-1)^{\iota_i a_i/2}
 e^{\pi iP_i a_i(a_i+2)/r}
\right]
\prod_b\tau_r(\mathbf a_b),
\tag{2.1}\label{eq-2-1}
\]

where \(\varepsilon_{r,\mathbf p}=\kappa_r^{-\sigma}\) is the fixed signature anomaly of modulus one, with

\[\sigma=\operatorname{sign}(L'\cup L_{\mathbf p})-\operatorname{sign}(L').\]

Here \(L'\) is the auxiliary surgery link for the ambient manifold, \(L_{\mathbf p}\) carries the specified integral filling framings, and \(\operatorname{sign}\) denotes the signature of the linking matrix. The sum is over admissible even component colors. Each six-tuple \(\mathbf a_b\) consists of the labels inherited by block \(b\), including repeated occurrences of a component. Thus \eqref{eq-2-1} contains the internal Racah sums through its factors \(\tau_r(\mathbf a_b)\).

\textbf{The two-block example.} For \(N_{\rm dbl}\) of \hyperref[sec-1-3-2]{Section 1.3.2}, the incidence table gives \(\mathbf a_+=\mathbf a_-=(a_1,\ldots,a_6)\), and all mutation integers vanish. Formula \eqref{eq-2-1} becomes

\[
Z_r(N_{\rm dbl}(\mathbf p))=
\varepsilon_{r,\mathbf p}D_r^2
\sum_{\substack{(a_1,\ldots,a_6)\in I_r^6\\\text{admissible}}}
\tau_r(a_1,\ldots,a_6)^2
\prod_{i=1}^6\left[
\frac2{\sqrt r}\sin\frac{2\pi(a_i+1)}r\,
e^{\pi i p_i a_i(a_i+2)/r}\right].
\]

Here admissibility requires the four triples \(123,156,246,345\) to be \(r\)-admissible. Each block weight still has its own internal Racah sum. The outer sum colors the six components, and the six bracketed factors perform the fillings along \(p_i\mu_i+\lambda_i\). The signature anomaly is computed from the chosen surgery presentation exactly as above. No value for it is assumed in this example.

The same incidence data give \(\Phi(x)=2\phi(x_1,\ldots,x_6)\) in the block expansion of \hyperref[sec-2-3]{Section 2.3}. Each strand multiplicity \(\nu_i\), counting the occurrences of component \(i\) among the block edges, equals two. Since \(\mathbf P=\mathbf p\) here, the phase of the leading Fourier term is

\[
S_{\mathbf P}(x)=2\phi(x_1,\ldots,x_6)
+i\pi\sum_{i=1}^6p_ix_i^2-i\pi\sum_{i=1}^6x_i.
\]

This exhibits the three contributions separately: the two geometric blocks, the six surgery coefficients, and the leading Fourier frequency.

\textbf{Central signs and block normalization.}

We keep track of two pieces separately: the overall phase of the product of block symbols, and its positive magnitude. The calculation below produces the sign \((-1)^{c(r-3)/2}\) and absorbs one factor \(D_r\) into each block. Both pieces reappear in the Fourier formula \eqref{eq-2-4}.

Set \(R=r-2\), \(y_i=a_i-R/2\), and for each block

\[w_b=z_b-R/8-\tfrac14\sum_{j=1}^6a_{b,j}.\]

Here \(z_b\) is the integer Racah index of block \(b\), \(a_{b,j}\) is its \(j\)th edge color, and the scaled internal coordinate is \(u_{z_b}=w_b/r\). When discussing one block we omit its subscript \(b\).

In a small fixed central domain, all nonzero Racah summands in a block have sign \((-1)^{1-(r-1)/2}\): the factorial \([z+1]!\) has \(z+1-(r-1)/2\) negative factors. Matching the face triples in the shadow presentation pairs their square roots. For each paired face, the sign of \(\Delta^2\) is \((-1)^{T_f+1-(r-1)/2}\); its three numerator factorials are positive in this region. There are \(2c\) paired faces, and

\[\sum_{\text{paired }f}T_f=\frac12\sum_{b=1}^c\sum_{j=1}^6a_{b,j}.\]

Indeed each block edge belongs to two faces, and pairing halves the sum over all block faces. Thus their total sign is \((-1)^{\frac12\sum_{b,j}a_{b,j}}\), which cancels \(\prod_b i^{-\sum_j a_{b,j}}\). Strand multiplicities are included in these sums. The remaining total block-symbol sign is

\[(-1)^{c(r-3)/2}.\]

The kernel \(K_h^+\) and its numerator and denominator forms \(e,d\) are defined explicitly in \hyperref[sec-A-2]{Appendix A.2}. With \(H_r(n)=\prod_{j=1}^n|2\sin(2\pi j/r)|\), the identity

\[H_r(r-1-n)=r/H_r(n)\]

rewrites the normalized positive block magnitude \(D_r|\tau_r|\) as \(r^{-1/2}\) times a Racah sum. More explicitly, \(|\tau_r|=(2\sin(2\pi/r)/r)\sum_z K_{1/r}^+(x,u_z)\), so the \(D_r^c\) in \eqref{eq-2-1} is absorbed block by block. The sixteen numerator factors have arguments \(R/4+re\), with power \(1/2\); the eight denominator factors have arguments \(R/8+rd\). Thus their shifts are precisely \(-1/2\) and \(-1/4\) in the \(r\)-scaled variables. Replacing \(H_r\) by \eqref{eq-A-4} at these central integer arguments is exact, and the internal zero Poisson mode of the normalized block is exactly \eqref{eq-A-13}.

The finite formula is the algebraic starting point. Before applying Poisson summation to its outer colors, we need to control two replacements: restriction to central colors and replacement of each internal Racah sum by a continuous block integral.

\subsection{2.3 The analytic input and its error}\label{sec-2-3}

We first explain the absolute estimate that allows us to restrict to central colors and internal indices. In the hyperideal admissible block domain, the classical function governing the absolute size of a Racah summand is bounded by \(v_8\), with equality only at the complete central point in the variables of \cite[Lemma 5.3]{PW}. On the boundary and the complement of that domain, the supremum is strictly smaller; \cite[proof of Proposition 5.24]{PW} states this complementary bound using \cite[Sections 3--4]{BDKY}. In our scale a block contributes \(v_8/(2\pi)\) at the complete point. Taking a finite product of blocks and a fixed central neighborhood therefore gives a strict loss away from their common complete configuration. The error \(O(\log r)\) in the logarithm of each finite factorial is uniform over the admissible indices. There are only polynomially many indices because the shadow presentation is fixed. The framing and mutation factors have absolute value one on real colors, while the quantum dimensions contribute at most a polynomial prefactor. The same exponential estimate therefore applies before summing the oscillating terms.

\leavevmode\phantomsection\label{stmt-proposition-2-1}\textbf{Proposition 2.1 (Localization).} For nested fixed central cutoffs, there are constants \(C_N,K,\kappa_N>0\) such that the contribution of colors or Racah indices outside the central region satisfies

\[|\text{noncentral contribution}|\le C_Nr^Ke^{r(A_N-\kappa_N)},
\quad A_N=\operatorname{Vol}(N)/(4\pi).\tag{2.2}\label{eq-2-2}\]

The strict maximum described above gives the gap \(\kappa_N\). The uniform \(O(\log r)\) factorial error and the polynomial number of summands contribute the factor \(C_Nr^K\).

\textbf{The continuous block model.} For one block, let \(\mathcal B_h(x)\) be the midpoint integral \eqref{eq-A-13}, with \(h=1/r\) and six external coordinates \(x\). It models the normalized positive weight \(D_r|\tau_r|\). More explicitly, \eqref{eq-B-6-1} writes that weight as \(\sqrt h\sum_z K_h^+(x,u_z)\). Since consecutive internal coordinates differ by \(h\), replacing the sum by an integral gives the normalization \(h^{-1/2}\int K_h^+(x,u)\,du\) used in \(\mathcal B_h\). A cutoff retains the neighborhood of its internal saddle.

The internal integration variable varies the Racah index while keeping the six external coordinates fixed. Let \(\phi(x)\) be the critical phase after this internal variable has been eliminated, as in \eqref{eq-A-16}. \hyperref[sec-A]{Appendix A} constructs the integral and proves the uniform expansion

\[\mathcal B_h(x)\sim e^{\phi(x)/h}\sum_{j\ge0}b_j(x)h^j.\]

The functions \(\phi,b_j\) are holomorphic near zero; \(b_0\) is positive on the central real domain and nowhere zero on a smaller complex neighborhood. The quadratic part of \(\phi\) is \(v_8/(2\pi)-\frac\pi2\sum_i x_i^2\). In the next theorem, \(\sigma x=(\sigma_1x_1,\ldots,\sigma_6x_6)\) denotes coordinate reflection. The theorem is the analytic input for the main proof.

\noindent\minipage{\linewidth}

\leavevmode\phantomsection\label{stmt-theorem-2-2}\textbf{Theorem 2.2 (Local independent reflection).}\par\nopagebreak[4]

There are a complex polydisc \(U\) about zero (a product of six small disks) and constants \(C,K,\eta>0\) such that, for all sufficiently small positive \(h\) and all sign vectors \(\sigma\in\{\pm1\}^6\),

\[
|\mathcal B_h(x)-\mathcal B_h(\sigma x)|
\le C h^{-K}\exp[(\Re\phi(x)-\eta)/h],\qquad x\in U.
\tag{2.3}\label{eq-2-3}
\]

Moreover, \(\phi(\sigma x)=\phi(x)\), and every coefficient \(b_j\) is independently even. This is the estimate called \hyperref[eq-2-3]{(ER)} in this paper.

\endminipage

\textbf{Proof reference and use.} \hyperref[sec-A]{Appendix A} proves the theorem by comparing \(\mathcal B_h\) with an exactly symmetric double sine integral. \hyperref[sec-B]{Appendix B} extends the estimate to compact subsets of the reflection domain and compares the integral with the finite Racah sum. No construction of the double sine is needed for the reduction below. The decisive feature of \eqref{eq-2-3} is the positive loss \(\eta\) in the exponential rate; agreement of coefficients to every fixed order alone would not control cancellation.

\textbf{Replacement of the block weights.}

For central real colors, apply the finite-lattice comparison \eqref{eq-B-5-2} and its exact normalization \eqref{eq-B-6-1}, block by block on a fixed compact central set. The finite-sum comparison in \hyperref[sec-B]{Appendix B} includes arbitrary offsets and the finite endpoints, and gives an error \(Ch^{-M}e^{(\phi-\eta)/h}\) with \(h=1/r\). The saddle-centered cutoff there and the fixed cutoff in \eqref{eq-A-13} agree on a common inner saddle neighborhood; on the support of their difference, the real part of the phase is uniformly below its critical value. Thus they can be interchanged at the same exponential precision. \hyperref[stmt-proposition-a-4]{Proposition A.4} bounds the error after taking block products and summing the remaining real colors. The surgery phase has modulus one there, so it introduces only the already counted polynomial weight. The finite endpoints are therefore included in the replacement error. What remains is a finite outer color sum whose block factors are the integrals just described.

The internal zero modes give

\[E_r(x)=\prod_b\mathcal B_{1/r}(x_b),\qquad x_i=y_i/r.\]

By \hyperref[eq-2-3]{(ER)}, replacing this by its independent componentwise symmetrization

\[E_r^{\rm sym}(x)=2^{-k}\sum_{\sigma\in\{\pm1\}^k}E_r(\sigma x)\]

has an absolute error \(C_Nr^Ke^{r(A_N-\eta_N)}\) after the remaining finite real-color sum. Moreover,

\[E_r^{\rm sym}(x)=e^{r\Phi(x)}a_r(x),\qquad
a_r(x)=\sum_{j=0}^{J-1}a_j(x)r^{-j}+O_{N,J}(r^{-J})\]

uniformly on a fixed complex neighborhood, with

\[\Phi=\sum_b\phi(x_b),\quad B_0=\prod_bb_0(x_b),\quad
\Phi=A_N-\tfrac\pi2\sum_i\nu_i x_i^2+O(|x|^4).\]

Here \(a_0=B_0\), and \(\nu_i\) counts the occurrences of component \(i\) among the six edge labels of all blocks. Each occurrence contributes one square \(x_i^2\) to the block expansion, which explains the coefficient \(\nu_i\); in particular \(\sum_i\nu_i=6c\). The expansion holds for every fixed \(J\ge1\). The functions \(a_r,a_j\) are independently even and holomorphic on the neighborhood; \(B_0\) is nonzero near zero. Holomorphy of the fixed-cutoff integrals and the uniform estimates of \hyperref[sec-A-5]{Section A.5} also give the derivative bounds on smaller neighborhoods used below.

\noindent\minipage{\linewidth}

\textbf{Error checkpoint.} All preceding errors have size \(O(r^Ke^{r(A_N-\eta_N)})\) for some fixed \(\eta_N>0\), after decreasing it to include localization. Let \(s_{\mathbf P}\) be the real part of the leading critical value obtained in \hyperref[sec-3]{Section 3}. To recover the original finite invariant from the even model, we require

\[A_N-s_{\mathbf P}<\eta_N.\]

\hyperref[sec-3]{Section 3} proves \(A_N-s_{\mathbf P}=O(\sum_iP_i^{-2})\). Thus sufficiently long fillings satisfy this inequality. This is the point at which the length hypothesis ensures that the reflection and finite-sum errors remain smaller than the term that survives cancellation.

\endminipage

\subsection{2.4 Poisson summation and parity}\label{sec-2-4}

On the even-color lattice the linear term in the surgery phase is constant, because \(p_i+\iota_i\) is integral. In each coordinate, Poisson summation on a lattice of spacing \(2/r\) contributes the density factor \(r/2\). Together with the surgery weight \(2/\sqrt r\), this gives \(\sqrt r\), and hence \(r^{k/2}\) in \(k\) coordinates. The offset of the lattice contributes the frequency-dependent phase shown below. Extracting the constant surgery phases gives

\[
Z_r(N(\mathbf p))=
\Theta_{r,\mathbf p}r^{k/2}
\sum_{n\in\mathbb Z^k}e^{-\pi i(r-2)\sum_i n_i/2}I_{r,n}
+O_{N,\mathbf p}(r^Ke^{r(A_N-\eta_N)}),
\tag{2.4}\label{eq-2-4}
\]

where

\[
I_{r,n}=\int\chi_1(x)E_r^{\rm sym}(x)
\prod_i\sin(2\pi x_i)
 e^{r(i\pi\sum_iP_ix_i^2-i\pi n\cdot x)}\,dx,
\tag{2.5}\label{eq-2-5}
\]

the integral is over \(\mathbb R^k\), \(\chi_1\) is a fixed independently even cutoff, and

\[
\Theta_{r,\mathbf p}=
\varepsilon_{r,\mathbf p}(-1)^{c(r-3)/2+k}
\exp\left[\pi i\Big(\sum_iP_i\Big)(-r/4+1-1/r)\right].
\]

The integrand in \eqref{eq-2-5}, apart from its Fourier factor, is odd in each coordinate: \(E_r^{\rm sym}\) and the cutoff are even, while \(\sin(2\pi x_i)\) is odd. Replacing \(x_i\) by \(-x_i\) therefore negates the integral and reverses the frequency \(n_i\). The ratio of the two lattice phases is \((-1)^{n_i}\), because \(r\) is odd. Their sum is consequently multiplied by \(1-(-1)^{n_i}\). It vanishes when \(n_i\) is even and doubles when \(n_i\) is odd. Applying this in every coordinate removes all orbits with an even coordinate. Each remaining orbit contains \(2^k\) terms that add with the same phase. This holds for every odd residue class of \(r\).

\noindent\minipage{\linewidth}

For clarity, write \(T_n=\lambda_r(n)I_{r,n}\), where \(\lambda_r(n)=e^{-\pi i(r-2)\sum_jn_j/2}\), and let \(n^{(i)}\) reverse the \(i\)th coordinate of \(n\). For \(n_i\ne0\), the entire cancellation is summarized by

\[\begin{gathered}
I_{r,n^{(i)}}=-I_{r,n},\qquad
\lambda_r(n^{(i)})=(-1)^{n_i}\lambda_r(n),\\[3pt]
\{n,n^{(i)}\}\ \longrightarrow\
T_n+T_{n^{(i)}}=\bigl(1-(-1)^{n_i}\bigr)T_n
\ \longrightarrow\
\begin{cases}
0,&n_i\text{ even},\\
2T_n,&n_i\text{ odd}.
\end{cases}
\end{gathered}\]

If \(n_i=0\), oddness gives \(I_{r,n}=0\) directly. Repeating the calculation in each coordinate leaves only odd frequencies. The cancellation is exact in the continuous even model; \eqref{eq-2-4} retains the exponentially small difference from the finite invariant.

\endminipage

\section{3. Asymptotics of the reduced sum}\label{sec-3}

The algebraic reduction is now complete: the invariant is an odd-frequency Fourier sum, up to the error controlled in \hyperref[sec-2-3]{Section 2.3}. We state the saddle-point result needed to evaluate it, explain why its first orbit dominates, and give the expansion coefficients. \hyperref[sec-C]{Appendix C} contains the contour and frequency estimates.

\subsection{3.1 Asymptotics of the odd Fourier terms}\label{sec-3-1}

The parity calculation leaves infinitely many odd Fourier frequencies. We now show that, when all \(|P_i|\) are sufficiently large, the leading frequencies are exactly \((\pm1,\ldots,\pm1)\). The theorem is stated for an analytic phase and amplitude so that the estimates can be read independently of the surgery formula. In \hyperref[sec-3-2]{Section 3.2} we use the functions constructed above.

\leavevmode\phantomsection\label{stmt-theorem-3-1}\textbf{Theorem 3.1 (Expansion for a signed quadratic lattice).}\par\nopagebreak[4]

Let \(\Phi\) be holomorphic near \(0\in\mathbb C^k\), real and independently even on the real locus, with

\[\Phi(0)=A,\qquad
\Phi''(0)=-\pi\operatorname{diag}(\nu_1,\ldots,\nu_k),\qquad\nu_i>0.\]

Here the double prime denotes the Hessian matrix, and \(\operatorname{diag}\) forms a diagonal matrix from its entries. Let \(a_r\) be independently even and holomorphic on a fixed complex neighborhood, with a uniform expansion to every fixed order in powers of \(r^{-1}\) there, and \(a_0(0)\ne0\). We write \(a_j\) for its coefficient of order \(r^{-j}\). Fix odd holomorphic functions \(\omega_i\) with \(\omega_i'(0)\ne0\). Choose a sufficiently small product of even smooth cutoffs \(\chi\), equal to one on a smaller box. Let \(\mathbf1=(1,\ldots,1)\in\mathbb R^k\) and \(\mathbf P=(P_1,\ldots,P_k)\). For odd \(r\), define

\[\mathcal L_r=\frac2r\left(\mathbb Z^k-\frac{r-2}{4}\mathbf1\right),\]

\[\mathcal Z_r(\mathbf P)=
\left(\frac2{\sqrt r}\right)^k
\sum_{x\in\mathcal L_r}\chi(x)a_r(x)\prod_i\omega_i(x_i)
\exp\left(r\left[\Phi(x)+i\pi\sum_iP_ix_i^2\right]\right).
\tag{3.1}\label{eq-3-1}\]

There is \(P_0>0\) such that for each fixed real vector \(\mathbf P\) with \(\min_i|P_i|>P_0\), the phase

\[S_{\mathbf P}(x)=\Phi(x)+i\pi\sum_iP_ix_i^2-i\pi\sum_i x_i\]

has a unique small critical point \(x_{\mathbf P}\), and for every fixed \(J\ge1\),

\[\mathcal Z_r(\mathbf P)=e^{-\pi i(r-2)k/2}e^{rS_{\mathbf P}(x_{\mathbf P})}
\left(\sum_{j=0}^{J-1}C_{\mathbf P,j}r^{-j}
+O_{\mathbf P,J}(r^{-J})\right).
\tag{3.2}\label{eq-3-2}\]

The leading coefficient is

\[C_{\mathbf P,0}=2^k(2\pi)^{k/2}
\frac{a_0(x_{\mathbf P})\prod_i\omega_i(x_{\mathbf P,i})}
{\sqrt{\det(-S_{\mathbf P}''(x_{\mathbf P}))}}\ne0.
\tag{3.3}\label{eq-3-3}\]

The determinant branch is selected by the Gaussian deformation. The parameters may have mixed signs and unrelated magnitudes. The remainder is for a fixed \(\mathbf P\); no uniform joint limit in \(r,\mathbf P\) is asserted.

\textbf{Why the first odd orbit dominates.} Near the complete structure, \(\Phi\) has a strict real maximum and diagonal negative Hessian. The quadratic filling term places the critical point of frequency \(n\) near \((n_i/(2P_i))_i\). Among the odd frequencies whose saddles remain in this neighborhood, the orbit \((\pm1,\ldots,\pm1)\) has the largest real critical value; \hyperref[stmt-lemma-c-1]{Lemma C.1} gives a strict quantitative gap. The remaining frequencies admit a summable bound with a fixed exponential loss. These estimates are proved in \hyperref[sec-C]{Appendix C}.

At the leading saddle, the loss from the complete value is \(O(\sum_iP_i^{-2})\). Taking the filling coefficients sufficiently large therefore makes this loss smaller than the exponential gaps in all previous replacement estimates. The \(2^k\) sign-related terms add, and Gaussian integration gives \eqref{eq-3-3}. In particular, the factor \(r^{-k/2}\) from the Gaussian cancels the \(r^{k/2}\) from Poisson summation. The leading coefficient is nonzero because each odd factor has a simple zero at zero and the leading saddle has nonzero coordinates. This proves the assertion once the estimates in \hyperref[sec-C]{Appendix C} are supplied.

\subsection{3.2 Expansion of the invariant and its coefficients}\label{sec-3-2}

For the shadow invariant, take \(\omega_i(x)=\sin(2\pi x)\), \(\Phi,a_r\) as in \hyperref[sec-2-3]{Section 2.3}, and \(\mathbf P=\mathbf p+\boldsymbol\iota/2\). All localization and reflection losses were fixed before the filling vector. Increasing \(P_0(N)\) once more makes \eqref{eq-A-33} hold for each of them. Write \(x_{\mathbf p}=x_{\mathbf P}\) and \(S_{\mathbf p}=S_{\mathbf P}(x_{\mathbf p})\). For every fixed \(J\ge1\),

\[\begin{aligned}
Z_r(N(\mathbf p))={}&\Theta_{r,\mathbf p}e^{-\pi i(r-2)k/2}
\mathcal A_{\mathbf p}e^{rS_{\mathbf p}}\\
&\times\left(\sum_{j=0}^{J-1}c_{\mathbf p,j}r^{-j}
+O_{N,\mathbf p,J}(r^{-J})\right),\qquad c_{\mathbf p,0}=1.
\end{aligned}\tag{3.4}\label{eq-3-4}\]

where

\[\mathcal A_{\mathbf p}=2^k(2\pi)^{k/2}
\frac{B_0(x_{\mathbf p})\prod_i\sin(2\pi x_{\mathbf p,i})}
{\sqrt{\det(-S_{\mathbf P}''(x_{\mathbf p}))}}\ne0.
\tag{3.5}\label{eq-3-5}\]

Formula \eqref{eq-3-4} separates the two geometric contributions. The exponential \(e^{rS_{\mathbf p}}\) will give the volume and Chern--Simons term. In \eqref{eq-3-5}, the Gaussian determinant combines with the block coefficient \(B_0\) and the sine filling factors to give the torsion coefficient; the determinant alone is only the filling-Jacobian part of that formula. The square root is fixed by continuous Gaussian deformation. The leading coefficient is nonzero by \hyperref[stmt-theorem-3-1]{Theorem 3.1}. The threshold was chosen so that every replacement error in \hyperref[sec-2-3]{Section 2.3} is exponentially smaller than the displayed leading term, which is why \eqref{eq-3-4} describes the original finite invariant.

\textbf{A formula for every coefficient.} There are two successive saddle-point calculations. Integrating the internal Racah variable of one block gives the coefficients \(b_j\). Multiplying the block expansions gives \(a_j\), and integrating the component variables gives the coefficients \(c_{\mathbf p,j}\) of the filled manifold. Thus the first correction combines the internal block corrections with derivatives of the outer phase and amplitude. The formulas below implement these three operations; \eqref{eq-3-5f} gives their one-variable form.

Put \(h=1/r\) and \(\Omega(x)=\prod_i\sin(2\pi x_i)\). Define the amplitude coefficients \(f_\ell(x)=\Omega(x)a_\ell(x)\). Since the block coefficients are independently even, the symmetrization in \hyperref[sec-2-3]{Section 2.3} does not change their products, and

\[a_\ell(x)=\sum_{j_1+\cdots+j_c=\ell}\prod_{b=1}^c b_{j_b}(x_b),
\qquad
\frac{f_1(x)}{f_0(x)}=\sum_{b=1}^c\frac{b_1(x_b)}{b_0(x_b)}.
\tag{3.5a}\label{eq-3-5a}\]

Here all \(j_b\) are nonnegative integers, and repeated component labels are included in \(x_b\). The ratios below are taken near the leading saddle, where \(f_0\) does not vanish. Set \(x_*=x_{\mathbf p}\), and write \(Q=-S_{\mathbf P}''(x_*)\) and \(V=Q^{-1}\). Thus \(V\) is the complex symmetric covariance matrix of the local Gaussian, with the determinant branch already fixed in \eqref{eq-3-5}. For a polynomial \(p\) in \(z=(z_1,\ldots,z_k)\), define its Gaussian contraction by

\[\langle p\rangle_V=
\left[\exp\!\left(\frac12\sum_{a,b=1}^kV_{ab}\partial_{z_a}\partial_{z_b}\right)p(z)\right]_{z=0}.
\tag{3.5b}\label{eq-3-5b}\]

For example, \(\langle z_i z_j\rangle_V=V_{ij}\) and \(\langle z_i z_j z_\ell z_m\rangle_V=V_{ij}V_{\ell m}+V_{i\ell}V_{jm}+V_{im}V_{j\ell}\). The differential series terminates on each polynomial, and we apply it coefficientwise to formal series with polynomial coefficients. Equivalently, odd monomials have zero contraction and an even monomial is evaluated by summing the products of covariances over all pairings of its factors. For \(m\ge3\), let \(S_m(z)=D^mS_{\mathbf P}(x_*)[z,\ldots,z]/m!\) be the homogeneous term of degree \(m\) in the Taylor expansion of the phase; \(D^m\) denotes the \(m\)th total differential. Let \(t\) be a formal variable, and let \([t^{2j}]\) denote extraction of its coefficient of degree \(2j\). Then

\[c_{\mathbf p,j}=\frac1{f_0(x_*)}[t^{2j}]
\left\langle
\left(\sum_{\ell\ge0}t^{2\ell}f_\ell(x_*+tz)\right)
\exp\!\left(\sum_{m\ge3}t^{m-2}S_m(z)\right)
\right\rangle_V.
\tag{3.5c}\label{eq-3-5c}\]

For fixed \(j\), only \(\ell\le j\), phase derivatives through order \(2j+2\), and finitely many amplitude derivatives occur. The same rule in one internal variable computes \(b_j/b_0\). Here \(S(x,u)\) is the single-block phase \eqref{eq-A-14}, and \(u_*(x)\) solves its internal critical equation \(S_u=0\). Keep the six external coordinates \(x\) fixed, expand in \(u\) about \(u_*(x)\), and replace \(f_\ell\) by the coefficient of \(h^\ell\) in the formal amplitude \(\exp(\log K_h^+-S/h)\). Equations \eqref{eq-A-11} and \eqref{eq-A-13} determine those coefficients by Bernoulli polynomials and finite products. Thus \eqref{eq-3-5c} gives a procedure involving only finite differentiation and Gaussian contractions at each requested order.

To justify the rule, expand about \(x_*\) and formally substitute \(x=x_*+tz\) with \(t^2=h\). The quadratic part of \(rS_{\mathbf P}\) becomes \(-z^TQz/2\), and the higher terms are exactly the exponential in \eqref{eq-3-5c}. The factor \(2^k(2\pi)^{k/2}/\sqrt{\det Q}\) is common to all coefficients and cancels after division by \eqref{eq-3-5}. The cutoff is one near the saddle, so it contributes no derivatives. The fixed-order remainder estimates of \hyperref[stmt-theorem-3-1]{Theorem 3.1} justify each finite truncation; the complementary contours and other Fourier orbits are exponentially smaller. Odd powers of \(t\) vanish under the Gaussian contraction. The contraction rule is the Gaussian pairing formula of Isserlis \cite{ISS}, extended algebraically to the complex symmetric covariance \(V\).

\textbf{The first correction.} To calculate \(c_{\mathbf p,1}\), one needs the block ratios \(b_1/b_0\), the first two derivatives of the outer amplitude, and the third and fourth derivatives of the outer phase, all at the leading saddle. The inverse Hessian \(V\) pairs these derivatives according to \eqref{eq-3-5b}.

Choose a local logarithm \(L(x)=\log f_0(x)\) near the nonzero value \(f_0(x_*)\). In the next formula, \(x_{*,b}\) is the six-tuple induced by \(x_*\) on block \(b\). All derivatives are evaluated at \(x_*\), with \(L_i=\partial_iL\), \(L_{ij}=\partial_i\partial_jL\), and \(S_{i_1\cdots i_m}=\partial_{i_1}\cdots\partial_{i_m}S_{\mathbf P}\). Repeated tensor indices \(i,j,\ell,m,n,s\) are summed from \(1\) to \(k\). Expanding \eqref{eq-3-5c} through degree two gives

\[\begin{aligned}
c_{\mathbf p,1}={}&\sum_{b=1}^c\frac{b_1(x_{*,b})}{b_0(x_{*,b})}
+\frac12V_{ij}(L_{ij}+L_iL_j)\\
&+\frac12 L_iS_{j\ell m}V_{ij}V_{\ell m}
+\frac18S_{ij\ell m}V_{ij}V_{\ell m}\\
&+\frac1{12}S_{ij\ell}S_{mns}V_{im}V_{jn}V_{\ell s}
+\frac18S_{ij\ell}S_{mns}V_{ij}V_{mn}V_{\ell s}.
\end{aligned}\tag{3.5d}\label{eq-3-5d}\]

The first term uses the explicit block correction \eqref{eq-A-30a}. The next two terms come from the quadratic amplitude term and the product of the linear amplitude term with the cubic phase. The quartic phase has three Gaussian pairings. The square of the cubic phase has fifteen pairings: six pair all three factors across the two cubic terms, and nine pair two factors within each cubic term. These counts give the coefficients \(1/8\), \(1/12\), and \(1/8\) in the last three terms.

All amplitude derivatives in \eqref{eq-3-5d} can be computed from the leading block formula \eqref{eq-A-28}:

\[\begin{aligned}
L_i&=\partial_i\log B_0+2\pi\cot(2\pi x_i),\\
L_{ij}&=\partial_i\partial_j\log B_0
-4\pi^2\delta_{ij}\csc^2(2\pi x_i).
\end{aligned}\tag{3.5e}\label{eq-3-5e}\]

Here \(\delta_{ij}\) is one for \(i=j\) and zero otherwise. The derivatives of order at least three of \(S_{\mathbf P}\) are precisely those of \(\Phi=\sum_b\phi(x_b)\), because the surgery term is quadratic. Derivatives of the block saddle are obtained by differentiating \(S_u(x,u_*(x))=0\); for a block coordinate \(x_a\), this gives \(\partial_a u_*=-S_{ux_a}/S_{uu}\). The logarithm branches do not affect any derivative in these formulas.

\textbf{One-variable example.} For a single external variable, put \(q=-S_{\mathbf P}''(x_*)\). The two cubic-square contractions combine, and \eqref{eq-3-5d} reduces to the worked expression

\[\begin{aligned}c_{\mathbf p,1}={}&\sum_{b=1}^c\frac{b_1(x_{*,b})}{b_0(x_{*,b})}
+\frac{L''+(L')^2}{2q}\\
&+\frac{S_{\mathbf P}'''L'}{2q^2}
+\frac{S_{\mathbf P}''''}{8q^2}
+\frac{5(S_{\mathbf P}''')^2}{24q^3},
\end{aligned}\tag{3.5f}\label{eq-3-5f}\]

with all quantities evaluated at \(x_*\). As an exactly solvable check, take the Gaussian model of \hyperref[stmt-example-a-5]{Example A.5}: \(a_r=1\), \(\Phi(x)=A-\pi\nu x^2/2\), and \(D=\nu-2iP\). Then \(q=\pi D\), \(f_0(x)=\sin(2\pi x)\), \(f_1=0\), and all higher phase derivatives vanish. The first relative correction is therefore \(f_0''/(2qf_0)=-2\pi/D\), in agreement with the exact factor \(\exp[-2\pi/(rD)]\) in \eqref{eq-A-34}.

\subsection{3.3 Insertion of an even function}\label{sec-3-3}

\leavevmode\phantomsection\label{stmt-corollary-3-2}\textbf{Corollary 3.2 (Fixed even insertions).}\par\nopagebreak[4]

Let \(F\) be a fixed independently even holomorphic function on the neighborhood in \hyperref[stmt-theorem-3-1]{Theorem 3.1}. Denote the lattice sum with amplitude \(a_rF\) by \(\mathcal Z_r[F]\). Then the quotient has an expansion to every fixed order, beginning with

\[\frac{\mathcal Z_r[F]}{\mathcal Z_r[1]}
=F(x_{\mathbf P})+O_{\mathbf P,F}(r^{-1}).
\tag{3.6}\label{eq-3-6}\]

\textbf{Proof.} The proof of \hyperref[stmt-theorem-3-1]{Theorem 3.1} applies to the numerator even when its leading coefficient vanishes: all parity and remainder estimates still hold. Its first coefficient is \(F(x_{\mathbf P})C_{\mathbf P,0}\). Divide its expansion by that of the denominator, whose first coefficient is nonzero. \(\square\)

For an insertion in the original finite invariant, we also need a bound on the full real color range to preserve the localization estimate of \hyperref[sec-2-3]{Section 2.3}. Fixed trigonometric polynomials satisfy this requirement. This is the case for the colored core curves treated in \hyperref[sec-4-3]{Section 4.3}.

\section{4. Hyperbolic geometry and the main theorem}\label{sec-4}

\hyperref[sec-3]{Section 3} has proved an expansion for the finite invariant in terms of a critical point, a critical value, and a Hessian. We now identify these quantities geometrically. The critical equation is the smooth Dehn-filling equation; its value gives complex volume, and its Hessian gives the filling Jacobian in the torsion formula. These identifications prove the \hyperref[stmt-main-theorem]{main theorem}. We then treat fixed colors on the core curves of the filling solid tori.

\subsection{4.1 The filling equations and complex volume}\label{sec-4-1}

Recall that \(u_i\) and \(v_i\) are the logarithmic meridian and longitude holonomies in the preferred marking. We first translate the block phase into the classical potential, including its additive constant. We then add the surgery term and impose its critical equation.

\textbf{Coordinates and branches.} Write \(\alpha_i\) for the classical angular variable and \(\xi_b\) for the internal variable in the potential of \cite[Section 4]{WY}. To distinguish the internal variable from a peripheral holonomy, denote the second argument of our block phase \(S\) by \(w\) in this paragraph. The coordinate dictionary is

\[
\begin{aligned}
\alpha_i&=\pi+2\pi x_i,&
u_i&=2i(\alpha_i-\pi)=4\pi i x_i,\\
\xi_b&=\frac{7\pi}{4}+2\pi(s_b+w),&
s_b&=\frac14\sum_{j=1}^6x_{b,j}.
\end{aligned}
\]

At the complete structure, \(x=0\), \(u_i=v_i=0\), and the internal critical variable is \(w=0\), so \(\xi_b=7\pi/4\). The positive real chamber \(x_i>0\) has tetrahedral dihedral angles \(\alpha_i-\pi=2\pi x_i\) and meridian cone angles \(4\pi x_i\). We continue the logarithmic holonomies from this chamber to a simply connected neighborhood of the complete structure. This fixes their signs; the other sign choices are the reflections already treated in \hyperref[sec-2-4]{Section 2.4}.

The dictionary concerns the classical phase. At finite level our centered colors give \(\alpha_i=2\pi(a_i+1)/r\), whereas the unshifted angular variable in the cited quantum formula is \(2\pi a_i/r\). Likewise, at an internal lattice point \(w=u_{z_b}\), the displayed \(\xi_b\) equals \(2\pi z_b/r+7\pi/(2r)\). These order-\(r^{-1}\) shifts are already included in the midpoint amplitude and its corrections in \hyperref[sec-A]{Appendix A}.

\noindent\minipage{\linewidth} \textbf{Lemma (Comparison of the potentials).} Let \(\Psi_{\rm NZ}\) be the Neumann--Zagier potential in the peripheral marking of \hyperref[sec-1-3-2]{Section 1.3.2}, with \(\partial_{u_i}\Psi_{\rm NZ}=v_i/2\). Its lift can be chosen so that, throughout the chosen neighborhood,

\[
\Psi_{\rm NZ}(u)=4\pi i\Phi(x)
-2\pi^2\sum_i\iota_i x_i^2
+\frac{\pi^2}{2}\sum_i\iota_i,
\qquad u_i=4\pi i x_i.
\]

In particular, the complete-point value of this lift is \(\Psi_{\rm NZ}(0)=i\operatorname{Vol}(N)+(\pi^2/2)\sum_i\iota_i\). \endminipage

\textbf{Proof.} We first check the normalization of one block. Let \(\Lambda(\theta)=-\int_0^\theta\log|2\sin t|\,dt\) be the Lobachevsky function on the real line. The normalization \(F(1/4)=0\) in \hyperref[sec-A]{Appendix A} gives \(F(t)=-\Lambda(2\pi t)/(2\pi)\) for \(0<t<1/2\). Denote the face and quadrilateral half-sums of the angular variables by \(\tau_f^\alpha=\frac12\sum_{i\in f}\alpha_i\) and \(\eta_g^\alpha=\frac12\sum_{i\in g}\alpha_i\). The coordinate change gives

\[
\begin{aligned}
\tau_f^\alpha-\alpha_i&=2\pi t_{e_{f,i}},&
\tau_f^\alpha-\pi&=\pi-2\pi t_{e_{f,0}},\\
\xi-\tau_f^\alpha&=2\pi t_{d_f},&
\eta_g^\alpha-\xi&=2\pi t_{d_g},&
2\pi-\xi&=2\pi t_{d_0}.
\end{aligned}
\]

Here the \(t_e,t_d\) are exactly the arguments in \eqref{eq-A-14}. Let \(U(\alpha,\xi)\) be the complex classical tetrahedral potential of the cited formulation, and let \(V(\alpha,\xi)\) be its real volume function on the real parameter domain. Using \(\Lambda(\theta+\pi)=\Lambda(\theta)\) and \(\Lambda(\pi-\theta)=-\Lambda(\theta)\) gives

\[
V(\alpha,\xi)=-\frac12\sum_e\Lambda(2\pi t_e)
+\sum_d\Lambda(2\pi t_d)=2\pi S(x,w),
\qquad U(\alpha,\xi)=2\pi^2+4\pi iS(x,w).
\]

This equality holds first for real parameters near the central point. The branch of \(U\) has \(U(\pi,\ldots,\pi,7\pi/4)=2\pi^2+2iv_8\). On our side, each logarithm in \(F\) is real on \(0<t<1/2\) and is continued in the strip \(0<\Re t<1/2\). Holomorphy therefore extends the equality to a common complex neighborhood, with this fixed additive constant. Since \(\partial\xi/\partial w=2\pi\), the internal critical equations agree as well.

Write \(\alpha_b\) for the six angular variables of block \(b\) and \(\xi_b(\alpha_b)\) for its internal critical coordinate \(\xi\). The assembled potential of \cite[Proposition 4.1]{WY}, equivalently \cite[Proposition 5.6]{PW}, is

\[
\mathcal U(\alpha)=\sum_bU(\alpha_b,\xi_b(\alpha_b))
-\frac12\sum_i\iota_i(\alpha_i-\pi)^2
+\frac{\pi^2}{2}\sum_i\iota_i
=2c\pi^2+\Psi_{\rm NZ}(u).
\]

Substitution of the block identity proves the stated formula. At zero, \(4\pi\Phi(0)=\operatorname{Vol}(N)=2cv_8\). The normalization of the cited potential is \(\Psi_{\rm NZ}(0)=i(\operatorname{Vol}(N)+iC_N)\) modulo \(\pi^2\mathbb Z\). The complete Chern--Simons class is \((\pi^2/2)\sum_i\iota_i\) modulo \(\pi^2\mathbb Z\); choosing its representative \(C_N=-(\pi^2/2)\sum_i\iota_i\) gives precisely our lift. These two representatives differ by an integer multiple of \(\pi^2\). Thus both the derivatives and the additive constant have been fixed. \(\square\)

\textbf{The filling equation.} Differentiating the lemma and using \(\partial_{u_i}\Psi_{\rm NZ}=v_i/2\) gives

\[\partial_i\Phi=v_i/2-\pi i\iota_i x_i.\]

The outer phase is obtained by adding the quadratic surgery term and the linear term of the leading Fourier frequency. We write \(S_{\mathbf P,\mathbf1}=S_{\mathbf P}\) to indicate that all frequencies are one. In its derivative, the mutation terms cancel: \(-\pi i\iota_i x_i+2\pi i(p_i+\iota_i/2)x_i=2\pi ip_ix_i\). Using \(u_i=4\pi ix_i\) gives

\[\partial_iS_{\mathbf P,\mathbf1}=\tfrac12(v_i+p_i u_i-2\pi i).\]

Thus the stationary equations are exactly \(v_i+p_i u_i=2\pi i\), the smooth filling equations from \hyperref[sec-1-3-3]{Section 1.3.3}. Hyperbolic Dehn filling identifies their small solution with the geometric structure on \(N(\mathbf p)\) when the fixed coefficients are sufficiently long. The critical-value identity of Neumann--Zagier and Yoshida \cite{NZ,Y}, in the convention of \cite[Section 2.2]{WY}, then gives

\[\Re S_{\mathbf p}=\operatorname{Vol}(N(\mathbf p))/(4\pi).\]

Together with \eqref{eq-3-4}--\eqref{eq-3-5} and the coefficient calculation in \hyperref[stmt-proposition-4-2]{Proposition 4.2} below, this gives the following theorem.

\leavevmode\phantomsection\label{stmt-theorem-4-1}\textbf{Theorem 4.1 (Fixed long integral fillings).}\par\nopagebreak[4]

For a fixed marked fundamental shadow link complement \(N\), there is \(P_0(N)>0\) such that every fixed integral filling vector satisfying

\[\min_i|p_i+\iota_i/2|>P_0(N)\]

gives a closed oriented hyperbolic manifold \(M=N(\mathbf p)\) with

\[
|Z_r(M)|=\frac{e^{r\operatorname{Vol}(M)/(4\pi)}}{\sqrt{|\mathcal T_M|}}
(1+O_{N,\mathbf p}(r^{-1}))
\]

along all sufficiently large odd levels. Thus

\[
TV_r(M)=\frac{1}{|\mathcal T_M|}
e^{r\operatorname{Vol}(M)/(2\pi)}(1+O_{N,\mathbf p}(r^{-1})).
\]

Here \(\mathcal T_M=\operatorname{Tor}(M;\operatorname{Ad}\rho_{\rm geom})\) in the convention of \cite{WY-T}. The phase-resolved expansion \eqref{eq-3-4} holds to every fixed order. The theorem implies eventual nonvanishing; the filling vector is held fixed in every limit.

\textbf{The complex critical value.} Keep the lift of \(\Psi_{\rm NZ}\) fixed in the lemma. For an integral filling, the oriented core longitude can be taken as \(-\mu_i\), since \((p_i\mu_i+\lambda_i)\cdot(-\mu_i)=1\). Its logarithmic holonomy is therefore \(-u_i\). The complex-volume identity in \cite[Section 2.2]{WY} gives the Chern--Simons value \(C\in\mathbb R/\pi^2\mathbb Z\) by

\[
\operatorname{Vol}(M)+iC=
\Psi_{\rm NZ}/i-\sum_i\frac{u_iv_i}{4i}-\frac\pi2\sum_i u_i
\pmod{i\pi^2\mathbb Z},
\]

where all quantities on the right are evaluated at the filling solution. To compare this with our critical value, substitute the lemma into \(S_{\mathbf P}\) and use \(P_i=p_i+\iota_i/2\). The mutation quadratics cancel, leaving

\[
S_{\mathbf P}(x)=\frac{\Psi_{\rm NZ}(u)}{4\pi i}
+\sum_i\frac{p_i u_i^2}{16\pi i}
-\frac14\sum_i u_i+\frac{\pi i}{8}\sum_i\iota_i.
\]

At the critical point, \(v_i=2\pi i-p_i u_i\). Dividing the complex-volume identity by \(4\pi\) and making this substitution gives

\[
\frac{\operatorname{Vol}(M)+iC}{4\pi}
=\frac{\Psi_{\rm NZ}(u)}{4\pi i}
+\sum_i\frac{p_i u_i^2}{16\pi i}
-\frac14\sum_i u_i\pmod{\pi i/4}.
\]

Comparison proves, with every normalization term retained,

\[S_{\mathbf p}=\frac{\operatorname{Vol}(M)+iC}{4\pi}+\frac{\pi i}{8}\sum_i\iota_i\pmod{\pi i/4}.\]

\textbf{The surgery phase.} It remains to combine the geometric phase with the constants extracted from the finite sum. Use the signature difference \(\sigma\) defined in \hyperref[sec-2-2]{Section 2.2} and put \(W=\sum_iP_i\). The anomaly is

\[\kappa_r=e^{-\pi i((r+1)/4+3/r)},\qquad \varepsilon_{r,\mathbf p}=\kappa_r^{-\sigma}.\]

Expand the signature factor already contained in \(\Theta_{r,\mathbf p}\), then combine \(\Theta_{r,\mathbf p}\) with the leading Fourier phase and the displayed critical value. Collecting the exponent terms proportional to \(r\), constant in \(r\), and proportional to \(r^{-1}\) gives a fixed \(B_{\mathbf p}\ne0\) and an integer

\[m=2(c-k)+\sigma-\sum_i p_i\]

such that, for a chosen representative of \(C\),

\[
Z_r(M)=B_{\mathbf p}
 e^{r(\operatorname{Vol}(M)+iC)/(4\pi)}
 e^{\pi i r m/4}
 e^{\pi i(3\sigma-W)/r}
 \left(\sum_{j=0}^{J-1}c_{\mathbf p,j}r^{-j}
 +O_{N,\mathbf p,J}(r^{-J})\right),
\tag{4.1}\label{eq-4-1}
\]

Here \(J\ge1\) is fixed, \(c_{\mathbf p,0}=1\), and \(|B_{\mathbf p}|=|\mathcal T_M|^{-1/2}\). The exact factor \(e^{\pi i(3\sigma-W)/r}\) is retained rather than included in the coefficients. Thus \(c_{\mathbf p,1}\) is exactly \eqref{eq-3-5d}; expanding this phase factor instead would replace the first relative coefficient by \(c_{\mathbf p,1}+\pi i(3\sigma-W)\). Changing the representative of \(C\) changes the integer \(m\) accordingly. In particular the original, unrephased invariant satisfies the branch-free ratio

\[
\lim_{r\to\infty,\ r\text{ odd}}\frac{Z_{r+8}(M)}{Z_r(M)}
=\exp\left[\frac2\pi(\operatorname{Vol}(M)+iC)\right].
\tag{4.2}\label{eq-4-2}
\]

The step of eight in \eqref{eq-4-2} removes the factor \(e^{\pi i r m/4}\) and also makes the result independent of the representative of \(C\). For a logarithmic formulation, fix a logarithm of the nonzero leading prefactor in each odd residue class modulo eight, and use the analytic logarithm of the relative factor near one. Its slope is \((\operatorname{Vol}+iC)/(4\pi)\) modulo \(\pi i/4\). In the opposite Chern--Simons sign convention, the corresponding invariant is \(CS_{\rm target}=-C\); the subscript records this comparison convention.

\subsection{4.2 The torsion coefficient}\label{sec-4-2}

\leavevmode\phantomsection\label{stmt-proposition-4-2}\textbf{Proposition 4.2 (The leading coefficient and adjoint twisted Reidemeister torsion).}\par\nopagebreak[4]

Let \(G_b\) be the block Gram matrix defined in \hyperref[sec-D-2]{Appendix D.2}, evaluated at \(x_{\mathbf p}\), and continue \(d_b=\sqrt{-\det G_b}\) from its positive value at real central parameters. Then

\[\mathcal A_{\mathbf p}^{\,2}\mathcal T_M=\pm i^{\,c-k},
\qquad
|\mathcal A_{\mathbf p}|=|\mathcal T_M|^{-1/2}.
\tag{4.3}\label{eq-4-3}\]

\textbf{Proof.} We compare the three factors in \eqref{eq-3-5} with the three factors in the torsion formula: the block amplitudes, the filling Hessian, and the sine terms.

\hyperref[sec-D-2]{Appendix D.2} gives \(b_0(x_b)^2=1/(8d_b)\), so \(B_0^2=8^{-c}/\prod_b d_b\). For the Hessian, define the filling Jacobian matrix by

\[\mathsf J=\left(\frac{\partial(v_i+p_i u_i)}{\partial u_j}\right).
\qquad S_{\mathbf P}''=2\pi i\mathsf J.
\tag{4.4}\label{eq-4-4}\]

The second identity follows by differentiating the geometric critical equation at the beginning of this section. For integral filling the core longitude may be taken as \(-\mu_i\), so its logarithmic holonomy is \(-u_i\). The established Gram-matrix and filling formulas \cite[Theorems 1.1 and 1.4]{WY-T} give

\[\mathcal T_M=\pm2^{3c-2k}\det\mathsf J
\prod_{b=1}^c\sqrt{\det G_b}
\prod_{i=1}^k\sinh(u_i/2)^{-2}.
\tag{4.5}\label{eq-4-5}\]

Their regularity hypotheses hold in the chosen neighborhood of the complete character \cite[Remark 1.3]{WY-T}. The restricted geometric character lies in its distinguished component, and sufficiently long smooth hyperbolic fillings have nondegenerate filling Jacobian and loxodromic core curves. Thus the factors used here are nonzero; the closed adjoint twisted Reidemeister torsion is the acyclic torsion, defined up to sign.

Square \eqref{eq-3-5}, use \eqref{eq-4-4} and \(\sin(2\pi x_i)=-i\sinh(u_i/2)\), and obtain

\[\mathcal A_{\mathbf p}^{\,2}
=(-i)^k2^{2k-3c}
\frac{\prod_i\sinh(u_i/2)^2}{\det\mathsf J\prod_b d_b}.
\tag{4.6}\label{eq-4-6}\]

Choose \(\sqrt{\det G_b}=id_b\) in \eqref{eq-4-5}. Changing those choices only changes its permitted sign. Multiplication of \eqref{eq-4-5} and \eqref{eq-4-6} proves \eqref{eq-4-3}. \(\square\)

This calculation fixes the absolute coefficient in the normalization of \hyperref[sec-1-3-4]{Section 1.3.4}. The remaining unit-modulus factors are the signature and Chern--Simons phases already displayed in \eqref{eq-4-1}.

\subsection{4.3 Fixed colors on the core curves}\label{sec-4-3}

Let \(\gamma_1,\ldots,\gamma_k\) be the core curves of the filling solid tori. Give them fixed even colors \(m_i\ge0\), using the framings represented by zero-framed meridians in the surgery diagram. Write \(Z_r(M,\gamma;\mathbf m)\) for the relative Reshetikhin--Turaev invariant of this colored pair in the normalization of \hyperref[sec-1-3-4]{Section 1.3.4}.

\leavevmode\phantomsection\label{stmt-theorem-4-3}\textbf{Theorem 4.3 (Geometric core characters).}\par\nopagebreak[4]

For a fixed filling in \hyperref[stmt-theorem-4-1]{Theorem 4.1} and fixed even \(\mathbf m\), the quotient \(Z_r(M,\gamma;\mathbf m)/Z_r(M)\) has an expansion at every fixed order. Its leading term is

\[\frac{Z_r(M,\gamma;\mathbf m)}{Z_r(M)}
=\prod_i\frac{\sinh((m_i+1)u_i/2)}{\sinh(u_i/2)}
+O_{N,\mathbf p,\mathbf m}(r^{-1}).
\tag{4.7}\label{eq-4-7}\]

Here \(u_i\) is evaluated at the geometric filling point. Each factor is the \(\operatorname{Sym}^{m_i}\) character of the holonomy of the corresponding core curve and is nonzero. Thus the relative invariant has the same volume exponent.

\textbf{Proof.} We first express a colored core curve as an insertion in the surgery sum, then apply \hyperref[stmt-corollary-3-2]{Corollary 3.2}. A core curve is represented in the surgery diagram by a meridian of the corresponding surgery component. The exact Hopf kernel in \cite[Section 1]{WY} is

\[H(a,m)=(-1)^{a+m}[(a+1)(m+1)].\]

A colored meridian replaces the original surgery dimension \(d_a\) by this kernel. For \(a,m\) even, and \(x=(a-(r-2)/2)/r\), their ratio is

\[F_m(x)=\frac{H(a,m)}{d_a}
=\frac{\sin(2\pi(m+1)x)}{\sin(2\pi x)}
=\sum_{j=0}^m e^{2\pi i(m-2j)x}.
\tag{4.8}\label{eq-4-8}\]

The centering uses \((a+1)/r=x+1/2\); the two sine signs cancel because \(m\) is even. The finite exponential sum shows that \(F_m\) is holomorphic and even. It also gives \(|F_m|\le m+1\) on the real color range and a uniform bound on a fixed complex neighborhood when \(m\) is fixed. These are precisely the two requirements from \hyperref[sec-3-3]{Section 3.3}: the global bound preserves localization, and evenness preserves the Fourier cancellation. \hyperref[stmt-corollary-3-2]{Corollary 3.2} therefore gives the quotient expansion for the finite invariant.

At the geometric saddle \(u_i=4\pi ix_i\), so \(F_{m_i}(x_i)=\sinh((m_i+1)u_i/2)/\sinh(u_i/2)\). This is the \(\operatorname{Sym}^{m_i}\) character evaluated at the core curve, whose logarithmic holonomy is \(-u_i\). It is unchanged by a choice of lift because \(m_i\) is even. A loxodromic core curve has complex length with nonzero real part, so neither the numerator nor denominator vanishes. This proves the theorem. \(\square\)

For color \(m=2\) on a single core curve, \eqref{eq-4-7} becomes

\[\lim_{r\to\infty,\ r\ {\rm odd}}
\frac{Z_r(M,\gamma;2)}{Z_r(M)}
=1+2\cosh u
=\operatorname{tr}(\operatorname{Ad}\rho_{\rm geom}(\gamma)).
\tag{4.9}\label{eq-4-9}\]

A different fixed integer framing multiplies the relative invariant by a fixed power of \(\theta_m=e^{\pi im(m+2)/r}\). It leaves the leading ratio unchanged and modifies the later coefficients. Colors proportional to \(r\) change the exponential phase and are not covered by this argument.

\section{Appendix A. Construction and reflection of the block integral}\label{sec-A}

This appendix proves \hyperref[stmt-theorem-2-2]{Theorem 2.2}. It contains the special-function construction and the estimates used as inputs in \hyperref[sec-2-3]{Section 2.3}. We interpolate the finite factorials, compare the resulting block integral with an exactly symmetric double sine integral, and bound the remainder strongly enough to obtain exponential reflection. The variables \(h,x\) are those of \hyperref[sec-2-1]{Section 2.1}; within this appendix, \(x\) always denotes the six external coordinates of a single block.

The three claims needed from the construction are separated below: \eqref{eq-A-13} specifies the block integral, \eqref{eq-A-30} identifies its even coefficients, and \eqref{eq-A-31} bounds the remainder. \hyperref[sec-A-1]{Section A.1} supplies the factorial estimate used in that final bound. The auxiliary double sine integral in \hyperref[sec-A-3]{Section A.3} is used to prove coefficient symmetry.

For one block the real reflection domain is

\[\mathscr D=\left\{x\in\mathbb R^6:
\sum_{i\in f}|x_i|<\frac12\ \text{for every }f\in\mathcal F\right\}.\]

On the positive orthant, the angles \(2\pi x_i\) are the dihedral angles of a strictly hyperideal tetrahedron; the zero coordinates are included as analytic limiting configurations. The signed domain is invariant under independent coordinate reflections, as explained in \hyperref[sec-B]{Appendix B}. The notation \(K\Subset\mathscr D\) means a compact set contained in this open domain. Estimates on a complex neighborhood \(U_K\) are uniform over that fixed neighborhood, or over a specified smaller one.

\textbf{Comparison of two integrals.} The block integral \(\mathcal B_h\) comes directly from the interpolated quantum factorials, but its reflection symmetry is not apparent from its integrand. We compare it with a double sine integral \(\mathcal C_h\) whose reflection symmetry is exact. Their asymptotic coefficients agree up to alternating signs, which forces all coefficients of \(\mathcal B_h\) to be even. A bound on the growth of the remainders then gives the exponentially small difference needed in the surgery sum.

\subsection{A.1 Midpoint factorial interpolation}\label{sec-A-1}

We begin with the positive sine factorial \(\prod_{j=1}^n2\sin(2\pi j/r)\) for \(n<r/2\). We need an interpolation that is holomorphic in its argument and whose asymptotic remainder is uniform on vertical strips. The gamma product below provides both properties. Here \(\Gamma\) is Euler's gamma function, with \(\Gamma(n+1)=n!\). For a positive real \(r\), put \(\omega=r/2\) and define

\[
\mathcal H_r(z)=
\left(\frac{2\pi}{\omega}\right)^z\Gamma(z+1)
\prod_{m=1}^{\infty}
\frac{\Gamma(\omega m+z+1)}
{\Gamma(\omega m-z)(\omega m)^{2z+1}}.
\tag{A.1}\label{eq-A-1}
\]

The product converges locally uniformly in \(-1<\Re z<\omega\); the logarithm of its \(m\)th factor is \(O(m^{-2})\) on compact subsets. It is positive for real \(z\) in this interval. Its two shift equations are

\[
\frac{\mathcal H_r(z+1)}{\mathcal H_r(z)}
=2\sin\frac{2\pi(z+1)}r,
\qquad
\frac{\mathcal H_r(z+\omega)}{\mathcal H_r(z)}
=2\sin\pi(z+1).
\tag{A.2}\label{eq-A-2}
\]

The first follows from the gamma functional equation and Euler's sine product. For the second, apply the shift to the finite product through \(m=M\). The ratio telescopes to

\[
\left(\frac{2\pi}{\omega}\right)^\omega
\frac{\Gamma(\omega(M+1)+z+1)\Gamma(\omega M-z)}
{\Gamma(z+1)\Gamma(-z)\,\omega^{2\omega M}(M!)^{2\omega}}.
\]

Stirling's formula gives the limit \(2\pi/[\Gamma(z+1)\Gamma(-z)]=2\sin\pi(z+1)\). Meromorphic continuation gives the equation beyond the initial strip where necessary. A similar telescoping product gives

\[
\mathcal H_r(z)\mathcal H_r(\omega-1-z)=\omega.
\tag{A.3}\label{eq-A-3}
\]

At nonnegative integers \(n<r/2\),

\[
\mathcal H_r(n)=\prod_{j=1}^{n}2\sin\frac{2\pi j}r.
\tag{A.4}\label{eq-A-4}
\]

The midpoint choice centers the factorial argument at a half integer. It pairs the complementary arguments in \eqref{eq-A-3} exactly and will make the formal logarithm odd in the small parameter. With \(h=1/r\), define

\[
M_h(t)=\sqrt{2h}\,\mathcal H_{1/h}(t/h-1/2).
\tag{A.5}\label{eq-A-5}
\]

In these coordinates, \eqref{eq-A-3} reads \(M_h(t)M_h(1/2-t)=1\), with \(M_h(1/4)=1\). Thus the reflection center is built into the normalized factorial before we form the block integral. Use the holomorphic functions

\[g(t)=\log(2\sin2\pi t),\qquad F'(t)=g(t),\qquad F(1/4)=0\]

on \(0<\Re t<1/2\), with their real branches on the real interval.

\textbf{Asymptotic precision.} The interpolation is defined for positive real \(h\); quantum invariants are recovered at \(h=1/r\) for odd \(r\). An expansion to every fixed order means that, after extracting the leading exponential and power, truncation after \(J\) terms has error \(O(h^J)\) for each fixed \(J\). For a normalized amplitude \(a_h(x)\) with formal coefficient functions \(a_j(x)\), a uniform Gevrey-one remainder bound has the stronger form

\[\left|a_h(x)-\sum_{j=0}^{J-1}a_j(x)h^j\right|
\le C A^J J!h^J,\]

with constants independent of the truncation order and of \(h\) in the stated range. Choosing \(J\) proportional to \(1/h\), with a sufficiently small proportionality constant, then gives an exponentially small remainder. Equality at every fixed algebraic order alone does not do this: \(e^{-1/\sqrt h}\) is smaller than every power of \(h\), but larger than \(e^{-\eta/h}\) for any fixed \(\eta>0\) as \(h\) tends to zero. The reflection error must also remain below the scale of the final signed sum; \hyperref[stmt-proposition-a-4]{Proposition A.4} makes this comparison precise. Subscripts on \(O\)-terms indicate the fixed data on which their constants may depend.

In the next estimate, \(B_j(z)\) is the Bernoulli polynomial of degree \(j\), in the convention \(te^{zt}/(e^t-1)=\sum_{j\ge0}B_j(z)t^j/j!\). The function \(\zeta\) is the Riemann zeta function, with \(\zeta(s)=\sum_{n\ge1}n^{-s}\) for real \(s>1\). Superscripts in parentheses denote derivatives.

\leavevmode\phantomsection\label{stmt-lemma-a-1}\textbf{Lemma A.1 (Explicit midpoint remainder).}\par\nopagebreak[4]

If \(0<\delta\le1/4\), \(\delta\le\Re t\le1/2-\delta\), and the truncation order is an integer \(N\ge2\), then

\[
\log M_h(t)=\frac{F(t)}h+
\sum_{m=1}^{N-1}\frac{B_{2m}(1/2)}{(2m)!}
h^{2m-1}g^{(2m-1)}(t)+R_N(t,h),
\tag{A.6}\label{eq-A-6}
\]

where

\[
|R_N(t,h)|\le
\frac{10\zeta(2N)(2N-2)!}{(2\pi)^{2N}}
\left(\frac h\delta\right)^{2N-1}.
\tag{A.7}\label{eq-A-7}
\]

The constants are independent of \(\Im t\). This uniformity will allow us to estimate the whole vertical contour in \hyperref[sec-A-3]{Section A.3}.

\textbf{Proof.}

\textbf{Step 1: a remainder for one gamma factor.} The midpoint form of Binet's formula \cite[Section 5.9]{DLMF} is

\[
\log\Gamma(w+1/2)=w(\log w-1)+\tfrac12\log(2\pi)
+\int_0^\infty e^{-wv}K_{1/2}(v)\,dv,
\]

where

\[
K_{1/2}(v)=\frac1{2v\sinh(v/2)}-\frac1{v^2}
=2\sum_{j\ge1}\frac{(-1)^j}{v^2+(2\pi j)^2}.
\tag{A.8}\label{eq-A-8}
\]

The subtraction removes the singularity at zero. Expand each denominator by a finite geometric series. The integral of the remainder is bounded by

\[
\frac{2\zeta(2N)(2N-2)!}
{(2\pi)^{2N}(\Re w)^{2N-1}}.
\tag{A.9}\label{eq-A-9}
\]

\textbf{Step 2: sum over the gamma product.} In \eqref{eq-A-1} with \(z=t/h-1/2\), all gamma arguments are midpoint arguments with

\[w=t/h,\qquad w=(m/2+t)/h,\qquad w=(m/2-t)/h.\]

For \(p=2N-1\ge3\),

\[
(\Re t)^{-p}+\sum_{m\ge1}
[(m/2+\Re t)^{-p}+(m/2-\Re t)^{-p}]
\le5\delta^{-p}.
\]

Multiplying this bound by the constant in \eqref{eq-A-9} proves \eqref{eq-A-7}. In particular, taking absolute values in Binet's integral has removed the imaginary part of each gamma argument, so the estimate remains valid at arbitrary imaginary height.

\textbf{Step 3: identify the coefficients.} Differentiate the leading terms of the gamma expansion and apply Euler's sine product. Their sum has derivative \(g(t)/h\). At \(t=1/4\), \eqref{eq-A-3} and \eqref{eq-A-5} give \(M_h(1/4)=1\), which fixes the integration constant and gives \(F(t)/h\). Applying the cotangent partial fraction expansion to the terms of order \(h^{2m-1}\) gives \(B_{2m}(1/2)g^{(2m-1)}(t)/(2m)!\). This proves \eqref{eq-A-6} on the real interval. Holomorphy gives the identity throughout the strip. \(\square\)

Summing the Binet remainders also gives an exact integral formula for the difference between the logarithm and its leading term:

\[
\log M_h(t)-\frac{F(t)}h
=h\int_0^\infty K_{1/2}(hv)
\frac{e^{-tv}-e^{-(1/2-t)v}}{1-e^{-v/2}}\,dv.
\tag{A.10}\label{eq-A-10}
\]

This follows by summing the Binet remainders in \eqref{eq-A-1} as a geometric series and setting the old integration variable equal to \(hv\). Near zero the quotient has a finite limit and \(K_{1/2}\) is bounded; at infinity the quotient decays exponentially uniformly on closed sub-strips. Thus the identity holds holomorphically throughout \(0<\Re t<1/2\). Its right-hand side vanishes at \(t=1/4\), in agreement with the exact midpoint normalization. Expanding the kernel reproduces \eqref{eq-A-6} with the same uniform remainder bounds.

\textbf{Consequences of the midpoint normalization.} The expansion in \hyperref[stmt-lemma-a-1]{Lemma A.1} has two properties that we will use. First, its formal logarithm

\[
\mathfrak L_h(t)=F(t)/h+
\sum_{m\ge1}\frac{B_{2m}(1/2)}{(2m)!}
h^{2m-1}g^{(2m-1)}(t)
\tag{A.11}\label{eq-A-11}
\]

is odd in \(h\):

\[\mathfrak L_{-h}(t)=-\mathfrak L_h(t).\tag{A.12}\label{eq-A-12}\]

Second, its coefficients and remainders have factorial bounds. On compact sub-strips, these bounds persist under the finite products, reciprocals, small \(h\)-translations, and square roots used below. Thus the resulting amplitudes have uniform Gevrey-one expansions in the sense defined above. The notation \(-h\) in \eqref{eq-A-12} refers only to the formal series; the functions themselves continue to be evaluated at \(h>0\).

To see why exponentiation preserves the required bound, if the logarithmic coefficient of degree \(j\) is bounded by \(A^{j+1}j!\), the recurrence for coefficients of its exponential gives a bound \(C^{j+1}j!\) after increasing \(C\). Truncate at degree \(N\), with \(Nh\) smaller than a fixed constant. The logarithmic remainder is \(O(C^{N+1}N!h^N)\), and the truncated lower-order logarithm stays bounded on a disk of radius proportional to \(1/N\). Cauchy's estimate on that disk bounds the degree-\(N\) exponential remainder by the same type of bound. Products and analytic changes of variables are handled in the same way.

\subsection{A.2 The midpoint block integral}\label{sec-A-2}

We now apply the interpolation to the factorials in a single tetrahedral weight. Its six external variables are \(x_1,\ldots,x_6\); the internal variable \(u\) corresponds to the Racah index. For fixed edge colors, summation in this internal variable evaluates the single block. Gluing identifies edge colors across blocks; the subsequent component-color sum performs the filling pairings. We use the face and quadrilateral sets from \hyperref[sec-1-3-5]{Section 1.3.5}:

\[\mathcal F=\{123,156,246,345\},\qquad
\mathcal Q=\{1245,1346,2356\}.\]

For \(x\in\mathbb C^6\) and \(u\in\mathbb C\), define

\[
s=\frac14\sum_i x_i,\quad
 t_f=\frac12\sum_{i\in f}x_i,\quad
 q_g=\frac12\sum_{i\in g}x_i.
\]

The four triangle factors give sixteen numerator forms

\[e_{f,0}=-t_f,\qquad e_{f,i}=t_f-x_i\quad(i\in f),\]

and factorial reflection of the Racah numerator leaves eight denominator forms

\[d_0=-s-u,\qquad d_f=s+u-t_f,\qquad d_g=q_g-s-u.\]

Write \(t_e=1/4+e\) and \(t_d=1/8+d\). Products indexed by \(e\) and \(d\) run over these sixteen and eight forms, respectively, with multiplicity. At \(x=u=0\), all numerator arguments are \(1/4\) and all unshifted denominator arguments are \(1/8\). \hyperref[sec-B-6]{Appendix B.6} checks the exact correspondence with the factorials in \(\tau_r\). The following table separates the geometric indices from the arguments of the interpolating function.

\begin{center}
\begin{tabular}{@{}>{\raggedright\arraybackslash}p{0.12\linewidth}>{\raggedright\arraybackslash}p{0.40\linewidth}>{\raggedright\arraybackslash}p{0.40\linewidth}@{}}
\toprule
Quantity & Range or definition & Role in the block weight \\
\midrule
$f\in\mathcal F$ & Four face triples & Index for $t_f$ and the four triangle factors \\[3pt]
$g\in\mathcal Q$ & Three quadrilateral index sets & Index for $q_g$ and the upper Racah bounds \\[3pt]
$e$ & $e_{f,0}$ and $e_{f,i}$, $i\in f$; sixteen forms & Numerator forms from the triangle factors \\[3pt]
$d$ & $d_0$, $d_f$, and $d_g$; eight forms & Denominator forms after factorial reflection \\[3pt]
$t_e$ & $1/4+e$ & Numerator argument of $M_h$ \\[3pt]
$t_d$ & $1/8+d$ & Denominator argument before the shift $h/4$ \\
\bottomrule
\end{tabular}
\end{center}

The variable \(u\) is the centered, rescaled Racah index: at finite level it takes the values \(u_z=(z-(r-2)/8-\frac14\sum_i a_i)/r\). Increasing \(z\) by one changes \(u_z\) by \(h\), with all six edge colors held fixed. The kernel interpolates the normalized internal summands, and the integral replaces their sum near the central critical point:

\[
K_h^+(x,u)=
\frac{\prod_e M_h(t_e)^{1/2}}
{\prod_d M_h(t_d+h/4)},
\qquad
\mathcal B_h(x)=h^{-1/2}\int_\mathbb R\chi(u)K_h^+(x,u)\,du.
\tag{A.13}\label{eq-A-13}
\]

Here \(\chi\) is fixed, smooth, one near zero, and supported in a sufficiently small real interval. Square roots are chosen positive and logarithms real on the central real domain; these branches are continued in a small simply connected complex neighborhood. The factors in \eqref{eq-A-13} are evaluated in a compact sub-strip, so they are holomorphic there.

The leading logarithm of the kernel is governed by the phase below. Here \(x\) fixes the six external colors, while \(u\) still varies the internal Racah index of this single block:

\[
S(x,u)=\frac12\sum_eF(t_e)-\sum_dF(t_d).
\tag{A.14}\label{eq-A-14}
\]

Direct differentiation, or the polynomial identity

\[\sum_dd^2=8u^2+\tfrac12\sum_i x_i^2,\]

gives

\[
S(x,u)=c_8-\frac\pi2\sum_i x_i^2-8\pi u^2+O(|(x,u)|^3),
\quad c_8=\frac{v_8}{2\pi}.
\tag{A.15}\label{eq-A-15}
\]

Subscripts on a phase denote partial derivatives in the indicated variables. The implicit-function theorem gives a unique analytic \(u_*(x)\) near zero with \(S_u(x,u_*)=0\). Put

\[\phi(x)=S(x,u_*(x)).\tag{A.16}\label{eq-A-16}\]

The equation \(S_u=0\) therefore evaluates the internal block sum. The Dehn-filling equations arise from differentiating with respect to the component colors, as in \hyperref[sec-4-1]{Section 4.1}. For real small \(x\), \(u_*(x)\) and \(\phi(x)\) are real, and \(\lambda(x)=-S_{uu}(x,u_*)>0\).

The local reflection statement is \hyperref[stmt-theorem-2-2]{Theorem 2.2}. Its proof occupies \hyperref[sec-A-3]{Sections A.3}--\hyperref[sec-A-5]{A.5} below: we construct an exactly symmetric integral, compare its coefficients with those of \(\mathcal B_h\), and estimate the remainder. The coordinates here are the six external variables of one block; the outer component variables were introduced in \hyperref[sec-2-1]{Section 2.1} of the main text.

\textbf{Extension to the hyperideal domain.} \hyperref[sec-B]{Appendix B} extends \hyperref[stmt-theorem-2-2]{Theorem 2.2} to every compact subset \(K\) of the reflection domain \(\mathscr D\) defined at the start of this appendix. Let \(\sigma x\) be obtained by changing any chosen coordinate signs. Here \(B_{h,K}\) is the internal integral localized near its saddle, as in \eqref{eq-B-4-1}, and \(\phi(x)\) is its critical phase. There are a complex neighborhood \(U_K\) of the reflection orbit of \(K\) and positive constants such that

\[|B_{h,K}(x)-B_{h,K}(\sigma x)|\le C_Kh^{-m_K}\exp\!\left(\frac{\Re\phi(x)-\eta_K}{h}\right),\qquad x\in U_K.\]

The factor \(e^{-\eta_K/h}\) is the essential part of this estimate: it controls errors even after larger terms have canceled. \hyperref[sec-B-5]{Appendix B.5} proves the corresponding comparison for finite internal sums, including their endpoints and lattice offsets. A reflected color need not belong to the set of even colors, so this comparison is made through the analytic interpolation.

\subsection{A.3 The double sine integral and its vertical tails}\label{sec-A-3}

The double sine provides the exactly symmetric integral used in the proof. Its difference equations belong to the analytic framework of Ruijsenaars \cite{RU}, and the associated continuous Racah--Wigner coefficients were developed by Ponsot--Teschner \cite{PT}. We use the integral and normalization of \cite{LMSWY} and \cite{MY}. The definitions below express that integral in the same variables as \(\mathcal B_h\).

\textbf{The double sine.} For \(b>0\), put \(Q=b+b^{-1}\) and let \(S_b\) be the normalized double sine function. Its shift equations are

\[S_b(z+b^{\pm1})=2\sin(\pi b^{\pm1}z)S_b(z),\]

and its reflection and normalization are \(S_b(z)S_b(Q-z)=1\) and \(S_b(Q/2)=1\). The \(b\)-\(6j\) integral is a continuous analogue of the tetrahedral weight with six parameters \(A_i\). It is invariant under independent substitutions \(A_i\mapsto Q-A_i\) when the parameters are \(b\)-admissible: all differences between the quadrilateral and face sums have real parts in \((0,Q)\). See \cite[Proposition 3.4]{LMSWY} for positive real \(b\) and \cite[Definition 1.1 and Proposition 2.3]{MY} for the continued formulation. We now write these sums explicitly and verify the inequalities for the integral used here.

Let \(b=\sqrt{2h}>0\), \(Q=b+b^{-1}\). The shift equations and the midpoint normalizations imply

\[
M_h(t)=S_b(2t/b+b/2).
\tag{A.17}\label{eq-A-17}
\]

For a direct uniqueness verification, compare \(S_b(b(z+1))/b\) with \eqref{eq-A-1}. They have the same shifts \(1\) and \(r/2\). If \(r/2\) is irrational their quotient has two incommensurable real periods and hence is constant as a meromorphic function; the reflection midpoint fixes it as one. Continuity in the positive parameter proves the identity at rational \(r/2\) as well. Meromorphic extension used in this argument follows from the two shifts.

Put \(A_i=Q/2+x_i/b\) in the exact \(b\)-\(6j\) integral. In this paragraph the continuous face sums are \(T_f=\sum_{i\in f}A_i\), and the four quadrilateral sums are \(Q_g=\sum_{i\in g}A_i\) for \(g\in\mathcal Q\), together with \(Q_0=2Q\). Thus \(j\) runs through \(\{0\}\cup\mathcal Q\); these are sums of the continuous parameters, while \hyperref[sec-1-3-5]{Section 1.3.5} used the same letters for half-sums of finite colors. The admissibility condition is \(0<\Re(Q_j-T_f)<Q\) for every pair of indices. For small real \(x\), these parameters are \(b\)-admissible: every difference \(Q_j-T_f\) is \(Q/2+2e/b\) for one of the sixteen forms \(e\), and lies strictly between zero and \(Q\). This remains true after arbitrary sign changes of the \(x_i\).

Writing \(U\) for the integration variable of the double sine formula, make the change of variable

\[U=7Q/4+(2s+2u)/b.\]

The eight numerator double sine arguments are \(Q/4+2d/b\). Using \eqref{eq-A-17}, the resulting integrand is

\[K_h^-(x,u)=\frac{\prod_d M_h(t_d-h/4)}{\prod_e M_h(t_e)^{1/2}}.\tag{A.18}\label{eq-A-18}\]

The multiset of sixteen denominator arguments agrees with the sixteen forms \(e\) above. \hyperref[sec-D]{Appendix D} displays the pair-sum identification.

Define

\[
\mathcal C_h(x)=\frac1{i\sqrt h}
\int_{u_*(x)+i\mathbb R}K_h^-(x,u)\,du.
\tag{A.19}\label{eq-A-19}
\]

For real small \(x\), this is a nonzero constant multiple, independent of \(x\), of the exact \(b\)-\(6j\) integral; the constant accounts for the normalization of the contour measure. Consequently its exact reflection identity gives

\[\mathcal C_h(\sigma x)=\mathcal C_h(x).\tag{A.20}\label{eq-A-20}\]

To use \eqref{eq-A-20} for the saddle expansion, we need to know that the part of the vertical contour away from the saddle is exponentially smaller. The next lemma gives this estimate. Its decay also justifies moving the vertical line within the admissibility strip, since the horizontal connecting segments then tend to zero.

\leavevmode\phantomsection\label{stmt-lemma-a-2}\textbf{Lemma A.2 (Uniform bound along vertical lines).}\par\nopagebreak[4]

There is \(c>0\) such that, for real \(x\) in a small fixed cube and all \(v\in\mathbb R\),

\[
\Re S(x,u_*(x)+iv)-\phi(x)
\ge c\min(v^2,|v|).
\tag{A.21}\label{eq-A-21}
\]

\textbf{Proof.}

Write \(a_d=t_d(x,u_*(x))\). After shrinking the cube,

\[1/16\le a_d\le3/16.\]

The coefficient of \(u\) in \(d\) is \(\varepsilon_d\in\{+1,-1\}\). There are four of each. For \(v>0\),

\[
\frac{d}{dv}\Re S(x,u_*+iv)
=\sum_d\varepsilon_d\Im g(a_d+i\varepsilon_dv)
=\sum_d\arctan\big(\cot(2\pi a_d)\tanh(2\pi v)\big).
\tag{A.22}\label{eq-A-22}
\]

Every summand is positive; on the displayed compact interval it is bounded below by \(c_1\min(v,1)\). Integrating proves \eqref{eq-A-21} for \(v>0\), and conjugation proves it for \(v<0\). \(\square\)

The midpoint bound \eqref{eq-A-7}, used with \(N=2\), is uniform on the entire vertical strip. Since \(g'\) is uniformly bounded there, it gives

\[
\log K_h^-(x,u)
=-S(x,u)/h-\tfrac14\sum_dg(t_d)+O(h)
\tag{A.23}\label{eq-A-23}
\]

uniformly along this vertical contour. Also

\[\exp[-\tfrac14\sum_d\Re g(a_d+i\varepsilon_dv)]\le C e^{-4\pi|v|}.\]

Together with \eqref{eq-A-21}, this proves absolute convergence and, for every fixed small \(v_0>0\), the tail bound

\[
\frac1{\sqrt h}\int_{|v|\ge v_0}|K_h^-(x,u_*+iv)|\,dv
\le C h^{-K}e^{-(\phi(x)+\kappa)/h}
\tag{A.24}\label{eq-A-24}
\]

for some \(\kappa>0\), uniformly in the real parameter cube. Thus the local saddle expansion describes the entire integral \(\mathcal C_h\), with its complementary contour contributing an exponentially smaller term.

\subsection{A.4 Reflection of the saddle coefficients}\label{sec-A-4}

We now prove that the phase and every coefficient of the block integral are invariant under reflection. The argument has two parts: compare the two saddle expansions, then use exact symmetry of the double sine integral.

\textbf{Comparison of the expansions.} The oddness of the midpoint logarithm in \eqref{eq-A-11}--\eqref{eq-A-12}, together with the opposite shifts in the two kernels, gives the formal identity

\[K_h^-(x,u)=K_{-h}^+(x,u).\tag{A.25}\label{eq-A-25}\]

Here \eqref{eq-A-25} means equality of the formal expansions after replacing \(h\) by \(-h\). Both contour integrals used in the proof have positive real \(h\) and \(b=\sqrt{2h}\).

Write the logarithmic expansion

\[\log K_h^+=S/h+A_0+hA_1+\cdots.\]

The finite shifts in \eqref{eq-A-13} yield

\[
A_0=-\tfrac14\sum_dg(t_d),\qquad
A_1=-\tfrac1{48}\sum_eg'(t_e)+\tfrac1{96}\sum_dg'(t_d).
\tag{A.26}\label{eq-A-26}
\]

There are holomorphic coefficient functions \(b_j(x)\) such that the real-contour saddle expansion is

\[\mathcal B_h(x)\sim e^{\phi(x)/h}\sum_{j\ge0}b_j(x)h^j.\tag{A.27}\label{eq-A-27}\]

The leading coefficient is

\[
b_0(x)=\sqrt{\frac{2\pi}{\lambda(x)}}
\exp[-\tfrac14\sum_dg(t_d(x,u_*(x)))],
\quad\lambda=-S_{uu}(x,u_*).
\tag{A.28}\label{eq-A-28}
\]

It is positive for real small \(x\) and nowhere zero in a smaller complex neighborhood.

Apply ordinary finite-order analytic stationary phase \cite{WONG}, \cite[Section 2.4]{DLMF} to the vertical integral \eqref{eq-A-19}, with its exponentially small tails controlled by \eqref{eq-A-24}. Equation \eqref{eq-A-25} gives

\[\mathcal C_h(x)\sim e^{-\phi(x)/h}\sum_{j\ge0}(-1)^jb_j(x)h^j.\tag{A.29}\label{eq-A-29}\]

To verify the sign without a continuation argument, use the analytic Morse coordinate \(v\) with \(S=\phi-v^2/2\). The real-contour formal Gaussian contractions have variance \(h\); on the vertical contour, put \(v=it\). Their variance becomes \(-h\), while every amplitude coefficient is also evaluated formally at \(-h\). Every complete contribution of total order \(j\) therefore changes by \((-1)^j\). The contour Jacobian \(i\) is canceled by the normalization in \eqref{eq-A-19}. The same reasoning is valid for every finite order.

\textbf{Use of exact reflection.} For real \(x\), the leading coefficient in \eqref{eq-A-29} is positive. Consequently \(h\log|\mathcal C_h(x)|\) tends to \(-\phi(x)\). Exact reflection \eqref{eq-A-20} first gives \(\phi(\sigma x)=\phi(x)\). After removing this common exponential, taking the limit identifies \(b_0\) at the two points. Subtracting that term and dividing by \(h\) identifies \(b_1\). Repeating the argument identifies every coefficient. Hence

\[
\phi(\sigma x)=\phi(x),\qquad b_j(\sigma x)=b_j(x)\quad(j\ge0)
\tag{A.30}\label{eq-A-30}
\]

on the real cube. The functions involved are holomorphic in \(x\), so the identities extend to a fixed complex polydisc. This proves the coefficient symmetry needed for the next section.

\textbf{The first two coefficients.} At the central point, the formulas give

\[b_0(0)=\frac1{4\sqrt2},\qquad \frac{b_1(0)}{b_0(0)}=\frac{13\pi}{24}.\]

More generally the first relative correction \(b_1(x)/b_0(x)\) equals

\[
A_1+\frac{A_0''+(A_0')^2}{2\lambda}
+\frac{S'''A_0'}{2\lambda^2}
+\frac{S''''}{8\lambda^2}
+\frac{5(S''')^2}{24\lambda^3},
\tag{A.30a}\label{eq-A-30a}
\]

where primes denote derivatives in the internal variable \(u\), evaluated at \(u_*(x)\). The five terms arise from the first amplitude correction, the quadratic amplitude term, the product of the linear amplitude and cubic phase terms, the quartic phase term, and the square of the cubic phase term, respectively. \hyperref[sec-D-1]{Appendix D.1} gives the derivatives that yield the displayed value at \(x=0\).

\subsection{A.5 Uniform remainders and exponential reflection}\label{sec-A-5}

We complete the proof of \hyperref[stmt-theorem-2-2]{Theorem 2.2} by estimating the remainder in \eqref{eq-A-27}. The point is to allow the truncation order to increase with \(1/h\).

\textbf{A factorial remainder for the saddle expansion.} \hyperref[stmt-lemma-a-1]{Lemma A.1} and the closure argument in \hyperref[sec-A-1]{Section A.1} give a holomorphic Gevrey-one expansion for \(e^{-S/h}K_h^+\). The analytic Morse change of variable makes the phase quadratic near the saddle. On a fixed disk in the Morse variable, Cauchy's estimate bounds the Taylor coefficient of degree \(2m\) in the amplitude coefficient of order \(h^j\) by

\[C^{j+1}j!\,\rho^{-2m}.\]

Here \(\rho>0\) is a fixed radius of the analytic Morse-coordinate disk. Gaussian integration contributes \((2m-1)!!h^m\), where \((2m-1)!!=1\cdot3\cdots(2m-1)\) and \((-1)!!=1\). Since

\[j!m!\le(j+m)!,\]

summing contributions with \(j+m<N\) and estimating the Taylor remainders gives constants \(A,C,\epsilon,\kappa>0\) such that

\[
\left|\mathcal B_h(x)-e^{\phi(x)/h}\sum_{j=0}^{N-1}b_j(x)h^j\right|
\le A C^NN!h^N e^{\Re\phi(x)/h}
+A h^{-K}e^{(\Re\phi(x)-\kappa)/h},
\tag{A.31}\label{eq-A-31}
\]

uniformly for \(x\) in a fixed smaller complex polydisc and \(1\le N\le\epsilon/h\).

\textbf{Treatment of the cutoff.} On a fixed inner interval where \(\chi=1\), the contour may be deformed holomorphically to the local saddle and Morse coordinate. On the connecting pieces, the real part of the phase is uniformly below its value at the saddle. The remaining support of \(\chi\) already has that separation on the original contour, and is bounded absolutely. Differentiating a nonanalytic cutoff or formally moving it is not required for \eqref{eq-A-31}.

\textbf{Derivation of the remainder bound.} We use finite Taylor expansions at each amplitude order. For each amplitude order \(j<N\), expand its analytic Morse amplitude only through degree \(2(N-j)-1\). Its Taylor remainder is at most a geometric constant times \(|v|^{2(N-j)}\) on a smaller disk. Integrating that power against the Gaussian gives \(C^{N-j}(N-j)!h^{N-j}\). Multiply by the \(j!h^j\) bound and sum over \(j\). The remaining amplitude remainder contributes \(C^NN!h^N\) directly. Outside the smaller disk the gap in the real part of the phase controls every term as long as \(Nh\) is sufficiently small. This proves the stated uniform range in \(N\).

\textbf{Optimal truncation.} Choose \(N=\lfloor\epsilon_0/h\rfloor\) with \(\epsilon_0<\min(\epsilon,1/(2C))\). Stirling's formula makes the first term in \eqref{eq-A-31} exponentially smaller than \(e^{\Re\phi/h}\). By \eqref{eq-A-30}, the two truncated series for \(x\) and \(\sigma x\) agree exactly. Subtract their expansions and bound both remainders. This proves \eqref{eq-2-3}, after decreasing \(\eta>0\) to absorb a polynomial prefactor.

\subsection{A.6 Error estimates for finite sums}\label{sec-A-6}

We record two consequences for later use. \hyperref[stmt-proposition-a-3]{Proposition A.3} abstracts the argument that transfers exact symmetry to an exponential estimate. \hyperref[stmt-proposition-a-4]{Proposition A.4} determines when that estimate is strong enough to survive a further finite summation. \hyperref[stmt-example-a-5]{Example A.5} shows why an expansion to every fixed order would not suffice.

\leavevmode\phantomsection\label{stmt-proposition-a-3}\textbf{Proposition A.3 (Transfer from an exactly symmetric auxiliary integral).}\par\nopagebreak[4]

Let a finite group \(G\) act holomorphically near a real domain, preserving that domain. Suppose two families have expansions

\[B_h(x)\sim e^{\phi(x)/h}\sum_{j\ge0}b_j(x)h^j,
\qquad C_h(x)\sim e^{-\phi(x)/h}\sum_{j\ge0}(-1)^jb_j(x)h^j.\]

Assume that \(\phi\) is real on the real domain, that \(\phi,b_j\) are holomorphic, and that \(b_0\) never vanishes. Suppose \(C_h(gx)=C_h(x)\) exactly for every \(g\in G\), and that its displayed expansion describes the whole auxiliary integral, including control of its complementary contour. Finally assume that \(B_h\) satisfies a bound of the form \eqref{eq-A-31} on a fixed \(G\)-invariant complex neighborhood, for \(1\le N\le\epsilon/h\). Then \(\phi\) and every \(b_j\) are \(G\)-invariant, and for some \(\eta>0\),

\[|B_h(gx)-B_h(x)|\le C_1h^{-m_1}e^{(\Re\phi(x)-\eta)/h}.
\tag{A.32}\label{eq-A-32}\]

\textbf{Proof.} The logarithmic rate \(h\log|C_h(x)|\to-\phi(x)\) on the real domain shows that \(\phi(gx)=\phi(x)\), since \(C_h(gx)=C_h(x)\). Successive subtraction identifies their coefficients. Holomorphic uniqueness extends these identities to a smaller common complex neighborhood. Take \(N=\lfloor a/h\rfloor\) with \(0<a<\min(\epsilon,1/(2C))\). The inequality \(N!\le N^N\) bounds the factorial term by \(A2^{-N}e^{\Re\phi/h}\). The truncated series at \(x\) and \(gx\) agree; their remainders give \eqref{eq-A-32}, for any sufficiently small \(\eta<\min(\kappa,a\log2)\). \(\square\)

For the block integral, the hypotheses of \hyperref[stmt-proposition-a-3]{Proposition A.3} are supplied by the kernel relation \eqref{eq-A-25}, the vertical tail estimate of \hyperref[stmt-lemma-a-2]{Lemma A.2}, and the factorial remainder \eqref{eq-A-31}. Thus the proposition summarizes the final step of the reflection proof.

\leavevmode\phantomsection\label{stmt-proposition-a-4}\textbf{Proposition A.4 (Error under finite products and sums).}\par\nopagebreak[4]

Suppose that, for real parameters, finitely many block families satisfy bounds given by \(e^{\phi_b/h}\) times a power of \(h^{-1}\), and their replacements have errors bounded by

\[|B_{b,h}-\widetilde B_{b,h}|\le Ch^{-M}e^{(\phi_b-\eta_b)/h},\qquad\eta_b>0.\]

Assume the remaining finite lattice sum has total absolute weight bounded by a power of \(h^{-1}\). Set \(A=\max\sum_b\phi_b\) on its real parameter range and \(\eta=\min_b\eta_b\). The error in the resulting sum is \(O(h^{-M'}e^{(A-\eta)/h})\). Here \(M'\) is a fixed nonnegative exponent accounting for the polynomial weights. In particular, it is exponentially smaller than a nonzero leading term \(h^\nu e^{s/h}\), with fixed real power \(\nu\) and real exponential rate \(s\), whenever

\[\eta>A-s.\tag{A.33}\label{eq-A-33}\]

\textbf{Proof.} Telescope the finite product, using one difference in each summand and the individual exponential bounds on the remaining factors. The exponential rate of each term is at most \(\sum_b\phi_b-\eta\). Apply the absolute weight bound and divide by the proposed leading scale. \(\square\)

Finite sign averaging preserves the estimate. A repeated external color is reflected at each occurrence in a block. For a different subsequent operation with exponential norm cost \(\Lambda\), the same proof requires \(\eta>A+\Lambda-s\); that norm bound must first be established. The number of blocks is fixed in this statement.

\leavevmode\phantomsection\label{stmt-example-a-5}\textbf{Example A.5 (Agreement at every algebraic order is insufficient).}\par\nopagebreak[4]

Consider the one-dimensional shifted lattice of \hyperref[sec-3-1]{Section 3.1}, with the whole Gaussian line in place of a cutoff, phase \(\Phi(x)=A-\pi\nu x^2/2\), odd factor \(\sin(2\pi x)\), and \(D=\nu-2iP\), where \(\nu>0\) and \(P\ne0\) are real and \(A\) is a real constant. Let \(n\in\mathbb Z\) be the Fourier index, with Fourier kernel \(e^{-i\pi rnx}\). After removing \(e^{rA}\), the Fourier integral is exactly

\[I_{r,n}=-i\sqrt{\frac2{rD}}
\exp\!\left(-\frac{\pi rn^2}{2D}-\frac{2\pi}{rD}\right)
\sinh\frac{2\pi n}{D}.\tag{A.34}\label{eq-A-34}\]

The unperturbed amplitude \(1\) cancels every even Fourier mode. For the leading \(n=\pm1\) pair, the real part of the critical value is \(A-d\), with

\[d=\frac{\pi\nu}{2(\nu^2+4P^2)}>0.\]

Replace the amplitude by \(1+\varepsilon_r x\), where \(\varepsilon_r=e^{-\sqrt r}\). On every fixed complex neighborhood this has the same power expansion as \(1\) to every fixed order. Its odd part nevertheless admits a zero Fourier mode. Including the Poisson factor \(\sqrt r\), that contribution is

\[\varepsilon_r\frac{2\sqrt{2/D}}{rD}
e^{rA-2\pi/(rD)}.\tag{A.35}\label{eq-A-35}\]

Its magnitude relative to the unperturbed sum is asymptotic to a positive constant times \(r^{-1}e^{dr-\sqrt r}\) and tends to infinity. For \(\varepsilon_r=e^{-\eta r}\) the relative change instead tends to zero if \(\eta>d\), and grows exponentially if \(\eta<d\). These formulas follow by integrating the two exponentials in the sine; differentiating the Fourier Gaussian gives \eqref{eq-A-35}. Thus the relevant comparison is between the reflection error exponent and the decrease \(d\) caused by cancellation. The perturbation \(e^{-\sqrt r}x\) has no contribution at any fixed algebraic order, but it is too large for this comparison. It also fails the uniform Gevrey-one remainder hypothesis of \hyperref[stmt-proposition-a-3]{Proposition A.3}.

\section{Appendix B. Reflection on the hyperideal domain and finite sums}\label{sec-B}

This appendix proves two extensions needed for the applications. The first is uniform reflection on a compact subset of the reflection domain \(\mathscr D\), rather than only near the complete point \(x=0\). The second is an exponentially accurate comparison between the block integral and the complete finite Racah sum, including its endpoints.

\hyperref[sec-B-1]{Sections B.1}--\hyperref[sec-B-3]{B.3} give the real and vertical phase estimates. \hyperref[sec-B-4]{Section B.4} uses them to extend the reflection argument of \hyperref[sec-A]{Appendix A}. \hyperref[sec-B-5]{Section B.5} proves the finite-sum comparison, and \hyperref[sec-B-6]{Section B.6} identifies its exact quantum normalization.

\subsection{B.1 The domain and the factorial arguments}\label{sec-B-1}

Keep \(F\), \(g\), \(M_h\), \(S\), \(K_h^+\) and the face and quadrilateral lists from \hyperref[sec-A-1]{Sections A.1}--\hyperref[sec-A-2]{A.2}. Thus \(F'=g=\log(2\sin 2\pi t)\) on \(0<\Re t<1/2\), with the real branches and \(F(1/4)=0\).

For real \(x=(x_1,\ldots,x_6)\) set

\[
\mathscr D=\left\{x:\sum_{i\in f}|x_i|<\frac12
\quad\text{for each }f\in\{123,156,246,345\}\right\}.
\tag{B.1.1}\label{eq-B-1-1}
\]

This is a convex open domain invariant under all independent sign changes. Its positive orthant is the strict hyperideal dihedral-angle domain after \(\theta_i=2\pi x_i\); zero coordinates are included in \(\mathscr D\) as analytic limiting configurations.

Use the linear forms \(s,t_f,q_g\) of \hyperref[sec-A-2]{Section A.2} and write

\[
\alpha_f=\frac18+s-t_f,\qquad
\beta_0=\frac18-s,\qquad
\beta_g=\frac18+q_g-s.
\]

The purpose of this notation is to separate the denominator arguments that increase with \(u\) from those that decrease with \(u\). There are four of each. In \hyperref[sec-B-1]{Sections B.1}--\hyperref[sec-B-3]{B.3}, the symbols \(a_f,b_j\) denote these factorial arguments; the coefficient functions of the saddle expansion will be written \(b_j(x)\) in \hyperref[sec-B-4]{Section B.4}. The index \(j\) for a factorial argument runs through \(\{0\}\cup\mathcal Q\), and the eight arguments before their \(h\) shift are

\[a_f=\alpha_f+u,\qquad b_j=\beta_j-u.\]

Pairing any increasing argument with any decreasing argument eliminates \(u\). The sixteen sums \(\alpha_f+\beta_j\) are exactly the numerator arguments \(t_e=1/4+e\) of \hyperref[sec-A-2]{Section A.2}, with multiplicity; \hyperref[sec-D-1]{Appendix D.1} displays the correspondence. Also

\[\sum_f\alpha_f+\sum_j\beta_j=1.\tag{B.1.2}\label{eq-B-1-2}\]

For each face, the bounds \(0<t_e<1/2\) say that the four expressions \(-x_i-x_j-x_k\), \(-x_i+x_j+x_k\), \(x_i-x_j+x_k\), and \(x_i+x_j-x_k\) have absolute value less than \(1/2\). Including their negatives gives all eight sign choices. Consequently these bounds are equivalent to \eqref{eq-B-1-1}.

Define

\[I_x=(-\min_f\alpha_f,\ \min_j\beta_j).\tag{B.1.3}\label{eq-B-1-3}\]

The lower pair bounds make this interval nonempty. For \(u\in I_x\),

\[a_f>0,\quad b_j>0,\quad a_f+b_j<1/2,\quad
\sum_fa_f+\sum_jb_j=1.\tag{B.1.4}\label{eq-B-1-4}\]

Every argument therefore lies in \((0,1/2)\). Individual arguments may exceed \(1/4\), so their cotangents may be negative. In the next two sections we control these signs by pairing an \(a\)-argument with a \(b\)-argument and using their sum.

\subsection{B.2 The unique real critical point}\label{sec-B-2}

We show that the internal phase is strictly concave throughout \(I_x\). The estimate also keeps the saddle a positive distance from both endpoints when \(x\) varies over a compact subset of \(\mathscr D\).

For \(a,b>0\) with \(a+b<1/2\),

\[
\cot(2\pi a)+\cot(2\pi b)
=\frac{\sin(2\pi(a+b))}{\sin(2\pi a)\sin(2\pi b)}
\ge2\cot(\pi(a+b)).
\tag{B.2.1}\label{eq-B-2-1}
\]

The purpose of the pairing is to use positivity of \(a+b\), even when one of the individual cotangents is negative. The last inequality follows from \(\sin(2\pi a)\sin(2\pi b)\le\sin^2(\pi(a+b))\). Pair the four \(a\)'s bijectively with the four \(b\)'s, and call their sums \(c_1,\ldots,c_4\). They lie in \((0,1/2)\) and sum to one. Since \(c\mapsto\cot(\pi c)\) is strictly convex on this interval, Jensen's inequality gives

\[
\sum_f\cot(2\pi a_f)+\sum_j\cot(2\pi b_j)
\ge2\sum_{i=1}^4\cot(\pi c_i)\ge8.
\]

Differentiating the phase proves the sharp real-concavity inequality:

\[S_{uu}=-2\pi\left(\sum_f\cot(2\pi a_f)+
\sum_j\cot(2\pi b_j)\right)\le-16\pi.\]

Moreover

\[S_u=-\sum_fg(a_f)+\sum_jg(b_j).\]

As \(u\) approaches the left endpoint of \(I_x\), at least one \(a_f\) tends to zero and all \(b_j\) stay in \((0,1/2)\), so \(S_u\to+\infty\). At the right endpoint \(S_u\to-\infty\) instead. Thus there is exactly one \(u_*(x)\in I_x\) with \(S_u=0\). It depends real analytically on \(x\) throughout \(\mathscr D\). Put

\[\phi(x)=S(x,u_*(x)),\qquad\lambda(x)=-S_{uu}(x,u_*(x))\ge16\pi.\]

Integrating the sharp real-concavity inequality twice yields the useful global internal gap

\[S(x,u)\le\phi(x)-8\pi(u-u_*(x))^2\quad(u\in I_x).\tag{B.2.2}\label{eq-B-2-2}\]

For any compact \(K\Subset\mathscr D\), the distance of \(u_*(x)\) from both endpoints has a positive minimum. All numerator and saddle denominator arguments therefore stay in a fixed compact subinterval of \((0,1/2)\).

\subsection{B.3 Growth of the phase on vertical lines}\label{sec-B-3}

The double sine integral has exponential factor \(e^{-S/h}\). We therefore need a lower bound on the increase of \(\Re S\) along a vertical line. The pairing from \hyperref[sec-B-2]{Section B.2} gives a bound valid on the whole domain.

Fix real \(x\in\mathscr D\), real \(u\in I_x\), and \(v>0\). In the strip used here \(2\sin(2\pi(a+iv))\) has positive real part, so its logarithm has imaginary part in \((-\pi/2,\pi/2)\). Direct differentiation gives

\[
\frac{d}{dv}\Re S(x,u+iv)=
\sum_f f_t(a_f)+\sum_j f_t(b_j),\qquad
f_t(a)=\arctan(t\cot(2\pi a)),\quad t=\tanh(2\pi v).
\tag{B.3.1}\label{eq-B-3-1}
\]

Some summands may be negative. The paired sums are controlled as follows. For fixed \(0<c<1/2\), the function \(f_t(a)+f_t(c-a)\) on \(0<a<c\) has its minimum at \(a=c/2\). Indeed,

\[
f_t'(a)=\frac{-2\pi t}{\sin^2(2\pi a)+t^2\cos^2(2\pi a)}.
\]

For \(a>c/2\) the difference of the two denominators, at \(a\) and \(c-a\), is

\[
(1-t^2)\sin(2\pi c)\sin(2\pi(2a-c))>0,
\]

so \(f_t'(a)-f_t'(c-a)>0\). Symmetry handles \(a<c/2\). Hence

\[f_t(a)+f_t(b)\ge2\arctan(t\cot(\pi(a+b))).\]

The function \(q_t(c)=\arctan(t\cot(\pi c))\) is convex for \(0<c<1/2\):

\[
q_t''(c)=
\frac{2\pi^2t(1-t^2)\sin(\pi c)\cos(\pi c)}
{(\sin^2(\pi c)+t^2\cos^2(\pi c))^2}>0.
\]

Pairing as in \hyperref[sec-B-2]{Section B.2} and using \(\sum c_i=1\) now proves

\[\frac{d}{dv}\Re S(x,u+iv)\ge8\arctan(\tanh(2\pi v)).\]

Integration and complex conjugation give \(\Re S(x,u+iv)-S(x,u)\ge\kappa(v)\), where \(\kappa(v)=8\int_0^{|v|}\arctan(\tanh(2\pi t))\,dt\). Since \(\arctan t\ge\pi t/4\) for \(0\le t\le1\), integration also gives the simpler explicit bound \(\log\cosh(2\pi v)\).

The sharp real-concavity bound \(S_{uu}\le-16\pi\) in \hyperref[sec-B-2]{Section B.2} and the sharp bound along vertical lines \(\Re S(x,u+iv)-S(x,u)\ge\kappa(v)\) both become equalities at \(x=0,u=0\). In particular, \(\kappa(v)=8\pi v^2+O(v^4)\) near zero and \(\kappa(v)=2\pi|v|+O(1)\) at infinity.

\subsection{B.4 Uniform reflection on compact subsets}\label{sec-B-4}

We now repeat the coefficient comparison and optimal truncation from \hyperref[sec-A]{Appendix A}. The preceding estimates provide the uniform saddle and tail bounds that were previously proved only near zero.

Let \(K\Subset\mathscr D\). Enlarge it to its finite sign orbit if necessary. Choose \(\epsilon>0\) so that \([u_*(x)-2\epsilon,u_*(x)+2\epsilon]\subset I_x\) for every \(x\in K\). For a fixed smooth cutoff \(\chi(w)\) supported in \(|w|<2\epsilon\) and equal to one for \(|w|\le\epsilon\), define

\[
B_{h,K}(x)=h^{-1/2}\int_{\mathbb R}\chi(w)
K_h^+(x,u_*(x)+w)\,dw.
\tag{B.4.1}\label{eq-B-4-1}
\]

The analytic saddle and its arguments extend holomorphically to a fixed sign-invariant neighborhood \(U_K\) of \(K\). After shrinking \(U_K\), \eqref{eq-B-4-1} is holomorphic there and has a uniform nondegenerate saddle at \(w=0\).

\textbf{The exactly symmetric integral.} Use the double sine integral of \hyperref[sec-A-3]{Section A.3}, with \(b=\sqrt{2h}>0\):

\[
\mathcal C_h(x)=\frac1{i\sqrt h}\int_{u_*(x)+i\mathbb R}
\frac{\prod_dM_h(t_d-h/4)}{\prod_eM_h(t_e)^{1/2}}\,du.
\tag{B.4.2}\label{eq-B-4-2}
\]

All \(b\)-admissibility inequalities follow from \(0<t_e<1/2\). Thus the exact reflection identity of Meng--Yang \cite[Proposition 2.3]{MY} applies for every \(x\in K\) and positive real \(b\).

Parametrize the vertical contour by \(u=u_*(x)+iv\), with \(v\in\mathbb R\). The midpoint remainder of \hyperref[stmt-lemma-a-1]{Lemma A.1} is uniform on vertical strips. On the contour in \eqref{eq-B-4-2} it gives, uniformly for \(x\in K\),

\[
\log K_h^-=-S/h-\tfrac14\sum_dg(t_d)+O_K(h),
\qquad
\exp\left(-\tfrac14\sum_d\Re g(t_d)\right)
\le C_K e^{-4\pi|v|}.
\]

Consequently the bound along vertical lines proves absolute convergence and, for fixed \(v_0>0\),

\[
\frac1{\sqrt h}\int_{|v|\ge v_0}|K_h^-(x,u_*+iv)|\,dv
\le C_Kh^{-1/2}\exp\left(-\frac{\phi(x)+\kappa(v_0)}h\right).
\tag{B.4.3}\label{eq-B-4-3}
\]

The same decay justifies any required horizontal change of the vertical line inside its admissibility strip. Thus \(\mathcal C_h(\sigma x)=\mathcal C_h(x)\) exactly.

\textbf{Equality of the coefficients.} The midpoint logarithm is odd in \(h\), and the two kernels obey \(K_h^-=K_{-h}^+\) as formal series. Finite-order saddle expansion therefore gives

\[
B_{h,K}(x)\sim e^{\phi(x)/h}\sum_{j\ge0}b_j(x)h^j,
\qquad
\mathcal C_h(x)\sim e^{-\phi(x)/h}\sum_{j\ge0}(-1)^jb_j(x)h^j,
\]

with

\[b_0(x)=\sqrt{2\pi/\lambda(x)}
\exp\left(-\tfrac14\sum_dg(t_d(x,u_*))\right)>0.\]

The vertical tail bound makes the second expansion an asymptotic of the exactly symmetric integral. Uniqueness of its real exponential asymptotic, followed by coefficient comparison, proves

\[\phi(\sigma x)=\phi(x),\qquad b_j(\sigma x)=b_j(x)\quad(j\ge0).\tag{B.4.4}\label{eq-B-4-4}\]

\textbf{The exponential remainder.} Holomorphy extends \eqref{eq-B-4-4} to \(U_K\). The midpoint bounds and the Morse-coordinate calculation of \hyperref[sec-A-5]{Section A.5} apply uniformly on \(K\), because the factorial arguments stay in a compact sub-strip and \(\lambda\ge16\pi\). In particular there are \(A,C,\epsilon_1,\eta>0\) and \(m\ge0\) such that for \(1\le N\le\epsilon_1/h\),

\[
\left|B_{h,K}-e^{\phi/h}\sum_{j<N}b_jh^j\right|
\le AC^NN!h^N e^{\Re\phi/h}+Ah^{-m}e^{(\Re\phi-\eta)/h}.
\tag{B.4.5}\label{eq-B-4-5}
\]

As in \hyperref[sec-A-5]{Section A.5}, analytic stationary phase is applied only where the cutoff equals one. Equation \eqref{eq-B-2-2} bounds the transition region by an exponentially smaller term. Optimizing \(N\) and using \eqref{eq-B-4-4} proves

\[
|B_{h,K}(x)-B_{h,K}(\sigma x)|
\le C_Kh^{-m_K}e^{(\Re\phi(x)-\eta_K)/h}\qquad(x\in U_K).
\tag{B.4.6}\label{eq-B-4-6}
\]

This proves \hyperref[eq-2-3]{(ER)} uniformly on compact subsets of \(\mathscr D\).

\subsection{B.5 Comparison with the finite Racah sum}\label{sec-B-5}

Exponential convergence for analytic trapezoidal sums has a general contour-integral theory \cite{TW}. Here the finite endpoints and arbitrary lattice offsets require the following explicit rectangle estimate. The conclusion \eqref{eq-B-5-2} concerns real external colors, which is sufficient for the finite surgery sum. The holomorphic object used later in the outer contour deformation is the block integral and its even symmetrization.

For real \(x\in K\) put

\[I_{x,h}=[-\min\alpha+h/4,\ \min\beta-h/4],\]

and for any real lattice offset \(\vartheta\) define

\[
\mathcal Q_{h,\vartheta}(x)=\sqrt h
\sum_{u\in h(\mathbb Z+\vartheta)\cap I_{x,h}}K_h^+(x,u).
\tag{B.5.1}\label{eq-B-5-1}
\]

The offset may depend on \(h\) and \(x\). There are constants independent of \(\vartheta\) such that

\[
|\mathcal Q_{h,\vartheta}(x)-B_{h,K}(x)|
\le C_Kh^{-m_K}e^{(\phi(x)-\eta_K)/h}.
\tag{B.5.2}\label{eq-B-5-2}
\]

We prove \eqref{eq-B-5-2} in two steps. First, an absolute bound removes the terms near the finite endpoints. Second, a cotangent contour compares the remaining inner lattice with its real integral. Both estimates are uniform in the offset \(\vartheta\).

\textbf{Step 1: the endpoint terms.} At zero, \(F\) is understood by its continuous limit, which exists because its derivative \(g(t)\) has an integrable logarithmic singularity there. For a fixed upper bound \(t_{\max}<1/2\), uniformly for \(0\le t\le t_{\max}\),

\[\log M_h(t+h/4)=F(t)/h+O_{t_{\max}}(1+|\log h|).\tag{B.5.3}\label{eq-B-5-3}\]

To see this, apply \hyperref[stmt-lemma-a-1]{Lemma A.1} with \(N=2\) at \(t+h/4\) and a strip margin bounded below by a constant times \(h\). Its remainder and \(h g'(t+h/4)\) are bounded. Also \([F(t+h/4)-F(t)]/h=O(1+|\log h|)\), including \(t=0\), because \(g(t)=\log t+O(1)\) there. The upper endpoints of the arguments stay below \(1/2\) by a fixed amount on \(K\).

Applying \eqref{eq-B-5-3} to the denominators, with the ordinary compact-substrip bound for the numerators, yields

\[|K_h^+(x,u)|\le C_Kh^{-m_K}e^{S(x,u)/h}
\quad(u\in\overline I_x).\tag{B.5.4}\label{eq-B-5-4}\]

There are \(O_K(h^{-1})\) lattice points. Equation \eqref{eq-B-2-2} therefore controls the entire lattice outside any fixed inner saddle interval exponentially, including the moving endpoints in \eqref{eq-B-5-1}.

\textbf{Step 2: the inner sum and its integral.} Choose endpoints \(L_h,R_h\) within \(O(h)\) of \(u_*\mp\epsilon\) such that \(L_h/h-\vartheta\) and \(R_h/h-\vartheta\) are half integers. On the rectangle with these endpoints and heights \(\pm s\), use the meromorphic kernel

\[W_h(z)=\frac1{2i}\cot\big(\pi(z/h-\vartheta)\big).\]

Its residue at a lattice point is \(h/(2\pi i)\), so the residue theorem gives \(h\sum K_h^+\) as the rectangle integral of \(W_hK_h^+\). On the lower and upper horizontal sides, respectively,

\[W_h=\tfrac12+O(e^{-2\pi s/h}),\qquad
W_h=-\tfrac12+O(e^{-2\pi s/h}).\]

On the vertical sides \(|W_h|\le1/2\) by the half-integer choice. Cauchy's theorem compares the average of the two horizontal integrals with the real integral; the remaining terms lie on those same vertical sides. All constants are independent of \(\vartheta\).

On this fixed analytic rectangle, \(\Re S(x,u+iv)\le S(x,u)+C_Kv^2\). Therefore horizontal errors are bounded by \(C_K\exp((\phi+C_Ks^2-2\pi s)/h)\), and vertical errors by \(C_K\exp((\phi-c_K\epsilon^2+C_Ks^2)/h)\). Choose \(s>0\) sufficiently small for both to have a fixed positive gap. The analytic amplitude is bounded uniformly there by the midpoint expansion. Multiplication by \(h^{-1/2}\), and the endpoint bound above, prove \eqref{eq-B-5-2}.

Combining \eqref{eq-B-4-6} and \eqref{eq-B-5-2}, for arbitrary offsets \(\vartheta,\vartheta'\),

\[
|\mathcal Q_{h,\vartheta}(x)-\mathcal Q_{h,\vartheta'}(\sigma x)|
\le C_Kh^{-m_K}e^{(\phi(x)-\eta_K)/h}.
\tag{B.5.5}\label{eq-B-5-5}
\]

Allowing two different offsets is useful because independent reflection can take an even color outside the set of even colors. Formula \eqref{eq-B-5-5} compares the interpolated finite sums on their respective lattices. \hyperref[sec-B-6]{Section B.6} identifies the actual Racah sums among these interpolated sums.

\subsection{B.6 Quantum normalization and a finite-difference identity}\label{sec-B-6}

For an actual even-color tuple in \(K\), let \(R=r-2\), \(h=1/r\), and \(x_i=(a_i-R/2)/r\). For the integer Racah index \(z\) set

\[u_z=\frac1r\left(z-\frac R8-\frac14\sum_i a_i\right).\]

The sixteen numerator factorial arguments are \(rt_e-1/2\); the eight denominator arguments are

\[R-z,\quad z-T_f,\quad Q_g-z,\]

which equal \(rt_d-1/4\). Their nonnegativity gives precisely \(I_{x,h}\). All numerator indices are below \(r/2\) for sufficiently large \(r\), uniformly on \(K\). The nonzero summands have a constant sign, since \(z+1>r/2\); the triangle prefactors contribute only a fixed phase for the fixed tuple.

In the convention \([n]=\sin(2\pi n/r)/\sin(2\pi/r)\), direct cancellation of factorial scales and \(H_r(r-1-n)=r/H_r(n)\) gives

\[
|\tau_r(\mathbf a)|=\frac{2\sin(2\pi/r)}r\sum_zK_h^+(x,u_z),
\qquad D_r|\tau_r(\mathbf a)|=\mathcal Q_{h,\vartheta}(x).
\tag{B.6.1}\label{eq-B-6-1}
\]

For clarity, the four triangle factors together give \((2\sin(2\pi/r))^2r^{-2}\prod_eH_r(rt_e-1/2)^{1/2}\); the reflected Racah numerator supplies \(r/(2\sin(2\pi/r))\). The factors \(\sqrt{2h}\) introduced in \(M_h\) cancel between its numerator and denominator. These cancellations give \eqref{eq-B-6-1}, including its power of \(r\).

This uses the symbol of \cite[Definition 2.3 and Proposition 2.6]{WY}. It is the normalization used when \(D_r\) is absorbed into each block in \hyperref[sec-2-2]{Section 2.2}.

The actual lattice offset is \(\vartheta=-R/8-\frac14\sum_i a_i\) modulo one, as follows directly from the formula for \(u_z\). Its dependence on the external colors is the reason \eqref{eq-B-5-2} was proved uniformly for every offset.

\textbf{A finite-difference identity.} The same normalization gives an exact identity for the finite internal sum. It is useful for checking the internal stationary equation. Define

\[
\mathsf P_h(u)=\prod_f2\sin(2\pi(\alpha_f+u+3h/4)),\qquad
\mathsf Q_h(u)=\prod_j2\sin(2\pi(\beta_j-u-h/4)).
\]

The exact shift of \(M_h\) proves

\[\mathsf P_h(u)K_h^+(x,u+h)=\mathsf Q_h(u)K_h^+(x,u).\tag{B.6.2}\label{eq-B-6-2}\]

For every function \(f\) on the finite lattice and its next point,

\[
\sum_{u=u_z}K_h^+(x,u)
\left[\mathsf Q_h(u)f(u+h)-\mathsf P_h(u-h)f(u)\right]=0.
\tag{B.6.3}\label{eq-B-6-3}
\]

Indeed, reindex the first term using \eqref{eq-B-6-2}. At the upper endpoint a minus factor in \(\mathsf Q_h\) is \(\sin0\), and at the lower endpoint a plus factor in \(\mathsf P_h(u-h)\) is \(\sin0\). Thus both boundary terms vanish exactly. This endpoint assertion is for the actual factorial lattice; arbitrary offsets in \eqref{eq-B-5-1} need not have the same exact vanishing.

Its classical equation is \(\mathsf Q_0=\mathsf P_0\), equivalently \(S_u=0\) on \(I_x\). If \(z=e^{4\pi iu}\), \(A_f=e^{4\pi i\alpha_f}\), \(B_j=e^{4\pi i\beta_j}\), and \(c=e^{2\pi i\sum\alpha_f}\), it becomes

\[c\prod_j(B_j-z)-c^{-1}\prod_f(A_fz-1)=0.\]

Because \(\prod A_f=c^2\) and \(\prod B_j=c^{-2}\), the constant and quartic coefficients cancel. After division by \(z\ne0\) this is a quadratic, possibly of lower degree at special parameters. This quadratic describes the internal stationary equation. The Dehn-filling equations arise separately from the derivatives with respect to the component colors, as in \hyperref[sec-4-1]{Section 4.1}.

\section{Appendix C. Estimates for the odd Fourier modes}\label{sec-C}

This appendix proves \hyperref[stmt-theorem-3-1]{Theorem 3.1}. The phase, amplitude, lattice, and odd factors are those in its statement. The parity reduction is the same as in \hyperref[sec-2-4]{Section 2.4}. We include it to keep the analytic proof self-contained, then compare the moderate frequencies and bound the sum of the high frequencies.

The order of the choices matters. First fix the local neighborhood and cutoff. Next choose the frequency cutoff \(\delta\) and the small contour displacement parameter \(s_0\). Their exponential losses are then fixed. Finally choose \(P_0\) large enough that the leading critical value lies above all these error rates, fix \(\mathbf P\), and let \(r\) tend to infinity. A moderate frequency satisfies \(|n_i/(2P_i)|<\delta\) for every \(i\); every other frequency is called high. Thus the moderate set is finite for each fixed filling vector.

\textbf{Proof of \hyperref[stmt-theorem-3-1]{Theorem 3.1}.} We first use parity to remove the even modes. We then divide the odd modes into two ranges: those whose critical points lie in the fixed local neighborhood and the remaining high frequencies. \hyperref[stmt-lemma-c-1]{Lemma C.1} compares the critical values in the first range; a contour estimate controls the second. Finally we evaluate the orbit of \((1,\ldots,1)\).

\subsection{C.1 Parity and comparison of critical values}\label{sec-C-1}

Poisson summation in \eqref{eq-3-1} gives \(r^{k/2}\sum_n \lambda_r(n)I_{r,n}\), where

\[\lambda_r(n)=e^{-\pi i(r-2)\sum_i n_i/2},\qquad
I_{r,n}=\int\chi a_r\prod_i\omega_i\,
e^{r(\Phi+i\pi\sum_iP_ix_i^2-i\pi n\cdot x)}\,dx.\]

For each fixed \(r\) the smooth compact support makes the Fourier series absolutely summable. Reflection of \(x_i\) negates \(I_{r,n}\) and reverses \(n_i\), while the quotient of lattice phases is \((-1)^{n_i}\). An orbit with an even coordinate therefore cancels; if \(n_i=0\), the corresponding integral is already zero. For odd coordinates the sign orbit adds. It suffices to compare \(n\) with positive odd coordinates, retaining all \(2^k\) orbit members for the leading coefficient.

\leavevmode\phantomsection\label{stmt-lemma-c-1}\textbf{Lemma C.1 (Comparison of critical values).}\par\nopagebreak[4]

For sufficiently small \(\delta>0\) and sufficiently large \(P_0\), whenever \(|P_i|>P_0\) and \(|n_i/(2P_i)|<\delta\) for all \(i\), the phase

\[S_{\mathbf P,n}(x)=\Phi(x)+i\pi\sum_iP_ix_i^2-i\pi n\cdot x\]

has one small critical point \(x^{(n)}\). For positive odd \(n\) in this range,

\[\Re S_{\mathbf P,\mathbf1}(x^{(\mathbf1)})
-\Re S_{\mathbf P,n}(x^{(n)})
\ge\frac\pi{16}\sum_i\nu_i\frac{n_i^2-1}{P_i^2}.
\tag{C.1}\label{eq-C-1}\]

\textbf{Proof.} Set \(t_i=P_i^{-1}\) and \(a_i=n_i/(2P_i)\). The critical equation is

\[x=a+\frac{i}{2\pi}\operatorname{diag}(t)\nabla\Phi(x).\]

The analytic implicit function theorem gives a unique small solution \(x(a,t)\), equivariant under independent signs of \(a\). Because \(x_i-a_i\) is divisible by \(t_i\), the function

\[H(a,t)=\Phi(x(a,t))+i\pi\sum_i
\frac{(x_i(a,t)-a_i)^2}{t_i}\]

extends holomorphically to \(t_i=0\). It is independently even in \(a\), and thus analytic in \(z_i=a_i^2\). Its real part for real \(a,t\) is the real part of the critical value of \(S_{\mathbf P,n}\), since the omitted term \(-i\pi\sum_i a_i^2/t_i\) is purely imaginary. At \(z=t=0\),

\[\partial_{z_i}\Re H=-\pi\nu_i/2.\]

In a smaller real parameter box this derivative is at most \(-\pi\nu_i/4\). Integrate successively from \(z_i=1/(4P_i^2)\) to \(z_i=n_i^2/(4P_i^2)\). This proves \eqref{eq-C-1} without a relation between the sizes or signs of the \(P_i\). \(\square\)

\subsection{C.2 Estimates for the moderate frequencies}\label{sec-C-2}

For the frequencies in \hyperref[stmt-lemma-c-1]{Lemma C.1}, we move the contour through the corresponding critical point. Since the cutoff is smooth, we keep track of its derivative terms by the following identity. Fix a real displacement vector \(v\in\mathbb R^k\). Let \(f\) be holomorphic on the tube swept out by \(x+itv\), \(0\le t\le1\), and let \(\psi\) be smooth with compact real support. For \(\xi\in\mathbb R^k\), put \(F(z)=f(z)e^{-i\xi\cdot z}\). All unspecified real integrals in this lattice argument are over \(\mathbb R^k\). Differentiation and real integration by parts give

\[\begin{aligned}
\int\psi(x)F(x)\,dx
={}&\int\psi(x)F(x+iv)\,dx\\
&+i\int_0^1\int(v\cdot\nabla\psi(x))F(x+itv)\,dx\,dt.
\end{aligned}\tag{C.2}\label{eq-C-2}\]

Indeed the derivative of \(\int\psi(x)F(x+itv)\,dx\) is \(-i\int(v\cdot\nabla\psi)F(x+itv)\,dx\). No continuation of \(\psi\) is used.

Choose the support and the region where \(\chi=1\) inside a fixed strict-concavity neighborhood of \(\Phi\). Next choose \(\delta\) much smaller than the distance to the cutoff-transition region. For the moderate frequencies of \hyperref[stmt-lemma-c-1]{Lemma C.1}, use \eqref{eq-C-2} with \(v=\Im x^{(n)}\). On the final slice, \(\Re S_{\mathbf P,n}(y+iv)\) has negative definite Hessian \(\Re\Phi''(y+iv)\) and its stationary real point is \(\Re x^{(n)}\). It is therefore a strict maximum. On the intermediate slices the surgery term contributes a real linear function with coefficients \(O(\delta)\): the critical equation gives \(P_iv_i=O(\delta)\). The Fourier term is constant in \(y\) and has size \(O(\delta^2)\). Consequently the fixed decrease in the real part of the phase on the cutoff-transition region persists after shrinking \(\delta\). Its size is fixed before choosing the large filling parameters.

There are only finitely many moderate modes for a fixed \(\mathbf P\). Holomorphic stationary phase at their nondegenerate saddles therefore gives expansions at every fixed order; \eqref{eq-C-2} makes the omitted cutoff terms exponentially smaller. A Taylor expansion in analytic Morse coordinates with Gaussian remainder bounds proves this assertion using exactly the fixed-order amplitude hypotheses in \hyperref[stmt-theorem-3-1]{Theorem 3.1}. No uniformity over an increasing set of filling vectors is required.

\subsection{C.3 A summable bound for the high frequencies}\label{sec-C-3}

A remaining frequency has \(|n_j/P_j|\ge2\delta\) for some \(j\). Split the integral with an even cutoff \(\rho(x_j)\) equal to one for \(|x_j|\le\delta/2\) and supported in \(|x_j|\le3\delta/4\). The piece with \(1-\rho\) has a fixed decrease in the real part of the phase, since \(\Phi\) has a strict maximum at zero.

For the inner piece, apply \eqref{eq-C-2} with

\[v=-\frac{s_0\operatorname{sgn}(n_j)}{|P_j|}e_j,\]

where \(s_0>0\) is fixed and small, \(e_j\) is the \(j\)th standard basis vector, and \(\operatorname{sgn}\) is the sign of a nonzero real number. At the final translated point \(y+iv\), the exact change in the real surgery and Fourier phase is

\[2\pi s_0\operatorname{sgn}(P_jn_j)y_j
-\pi s_0|n_j/P_j|
\le-\pi\delta s_0/2.\]

The block phase costs \(O(s_0^2/P_j^2)\) because it is real on the real locus. The same inequality with \(s_0\) replaced by \(ts_0\) holds along the transport. The derivative of \(\chi\rho\) in direction \(v\) is supported in the \(\rho\) transition, where \(|x_j|\ge\delta/2\): the original cutoff is identically one in coordinate \(j\) on the support of \(\rho\). Thus the transport terms also retain a fixed positive loss. Choose \(P_0\) large enough to absorb the quadratic block cost.

The exponential bound just obtained must also be summable over \(n\). We obtain the required decay in \(n\) by integration by parts, including in the transport terms of \eqref{eq-C-2}. In \eqref{eq-C-2} take

\[f_r(z)=a_r(z)\prod_i\omega_i(z_i)
e^{r(\Phi(z)+i\pi\sum_iP_iz_i^2)},\qquad \xi=\pi r n.\]

For a fixed \(\mathbf P\) and fixed derivative order \(L\), derivatives of \(f_r\) have bounds \(C_{\mathbf P,L}r^L\) times the absolute value of its exponential factor on smaller tubes. This follows by the product rule and Cauchy bounds for the amplitude and phase. Integrate by parts \(L\) times in a coordinate where \(|n_\ell|\) is largest, leaving \(e^{-i\pi rn\cdot x}\) as the Fourier factor. Its denominators \((\pi r|n_\ell|)^L\) cancel the derivative powers of \(r\); the translation factor \(e^{\pi rn\cdot tv}\) is independent of \(x\). Derivatives of the cutoffs remain in the same regions where the real part of the phase is lower. Apply this argument also to the unshifted outer piece and then integrate in \(t\). For \(L>k\) it gives

\[|I_{r,n}|\le C_{\mathbf P,L}r^M(1+|n|)^{-L}
e^{r(A-\eta_0)}
\tag{C.3}\label{eq-C-3}\]

for all high frequencies, with \(\eta_0>0\) fixed before the sufficiently large \(\mathbf P\) is fixed. The sum over these frequencies converges and has the same exponential loss.

\subsection{C.4 The leading orbit and its coefficient}\label{sec-C-4}

The leading saddle satisfies

\[A-\Re S_{\mathbf P}(x_{\mathbf P})
=O\left(\sum_iP_i^{-2}\right),\qquad
x_{\mathbf P,i}=\frac1{2P_i}\left(1+O(P_{\min}^{-1})\right).\]

The second estimate follows from the critical equation and the divisibility of \(\partial_i\Phi\) by \(x_i\). Choose \(P_0\) so that the first loss is smaller than half the fixed cutoff and high-frequency losses. \hyperref[stmt-lemma-c-1]{Lemma C.1} separates each of the other finitely many moderate odd modes. The remaining orbit consists of the \(2^k\) reflections of \(n=\mathbf1\), all with the same contribution after their lattice phases are included. Gaussian integration supplies \((2\pi/r)^{k/2}\); the Poisson prefactor cancels \(r^{-k/2}\). Every \(\omega_i(x_{\mathbf P,i})\) is nonzero by its simple zero at zero, while \(a_0\) and the Hessian remain nonzero. This proves \eqref{eq-3-2}--\eqref{eq-3-3} and \hyperref[stmt-theorem-3-1]{Theorem 3.1}. \(\square\)

\section{Appendix D. Algebraic identities and the Gram determinant}\label{sec-D}

\subsection{D.1 Pairing identities and central coefficients}\label{sec-D-1}

We record the finite calculations used in the reflection and normalization arguments. For the pairing in \hyperref[sec-B]{Appendix B}, the entry in row \(f\) and column \(j\) of the following table is the numerator form \(e\) for which \(\alpha_f+\beta_j=1/4+e\):

\begin{center}
\begin{tabular}{@{}lllll@{}}
\toprule
$f$ & $j=0$ & $j=1245$ & $j=1346$ & $j=2356$ \\
\midrule
$123$ & $-t_{123}$ & $t_{345}-x_3$ & $t_{246}-x_2$ & $t_{156}-x_1$ \\
$156$ & $-t_{156}$ & $t_{246}-x_6$ & $t_{345}-x_5$ & $t_{123}-x_1$ \\
$246$ & $-t_{246}$ & $t_{156}-x_6$ & $t_{123}-x_2$ & $t_{345}-x_4$ \\
$345$ & $-t_{345}$ & $t_{123}-x_3$ & $t_{156}-x_5$ & $t_{246}-x_4$ \\
\bottomrule
\end{tabular}
\end{center}

For example, the entry in row \(123\) and column \(1245\) reads \(\alpha_{123}+\beta_{1245}=1/4+t_{345}-x_3\), obtained by substituting the definitions of \(\alpha,\beta,s,t_f,q_g\). The other entries use the same linear substitution. Thus the sixteen entries, with multiplicity, are exactly the sixteen forms of \hyperref[sec-A-2]{Section A.2}. Expanding the eight denominator squares gives \(\sum_dd^2=8u^2+\frac12\sum_ix_i^2\). Both calculations are identities of polynomials with rational coefficients.

At the central point, \(\lambda=16\pi\), \(A_0'=S'''=0\), \(A_0''=16\pi^2\), \(S''''=-256\pi^3\), and \(A_1=\pi/6\). Substitution in the first-correction formula gives \(b_0(0)=1/(4\sqrt2)\) and \(b_1(0)/b_0(0)=13\pi/24\).

More explicitly, the odd-derivative terms vanish at the central point, and the remaining terms in the relative first-correction formula are

\[
\frac{b_1(0)}{b_0(0)}
=\frac\pi6+\frac{16\pi^2}{2(16\pi)}
-\frac{256\pi^3}{8(16\pi)^2}
=\frac{13\pi}{24}.
\]

The next subsection computes the same leading coefficient away from the central point and identifies the square-root branch needed in the torsion formula.

\subsection{D.2 The local Gram normalization}\label{sec-D-2}

Chen--Murakami \cite[Theorem 2]{CM} identified the Gram determinant in the leading hyperbolic quantum \(6j\) coefficient. For \hyperref[stmt-proposition-4-2]{Proposition 4.2}, we need its constant and square-root branch in our normalization. We derive them here from the internal critical equation.

Use the parameters of \hyperref[sec-B]{Appendix B} \(\alpha_f,\beta_j,u_*(x)\), and set

\[
a_f=\alpha_f+u_*,\quad b_j=\beta_j-u_*,
\quad
p_*=\prod_f\sin(2\pi a_f)=\prod_j\sin(2\pi b_j).
\]

The last equality is the internal stationary equation. For real \(x\in\mathscr D\), all eight factorial arguments \(a_f,b_j\) lie in \((0,1/2)\), so \(p_*>0\). Let

\[
T_*=\sum_f\cot(2\pi a_f)+\sum_j\cot(2\pi b_j),\qquad
\lambda=2\pi T_*.
\]

The leading saddle coefficient in \eqref{eq-A-28} gives

\[
b_0(x)^2
=\frac{2\pi}{\lambda}
\exp\left(-\frac12\sum_d g(t_d)\right)
=\frac1{16T_*p_*}.
\tag{D.1}\label{eq-D-1}
\]

Define

\[
G(x)=
\begin{pmatrix}
1&-c_1&-c_2&-c_6\\
-c_1&1&-c_3&-c_5\\
-c_2&-c_3&1&-c_4\\
-c_6&-c_5&-c_4&1
\end{pmatrix},\qquad c_i=\cos(2\pi x_i).
\]

Then

\[
\boxed{b_0(x)^2=\frac1{8\sqrt{-\det G(x)}}.}
\tag{D.2}\label{eq-D-2}
\]

The square root is positive on the real domain and equals \(4\) at \(x=0\). The identity extends holomorphically to a sufficiently small complex neighborhood, and locally to the compact domains of \hyperref[sec-B]{Appendix B} with the continued branch.

\textbf{Algebraic proof.} We first express the critical equation as a quadratic, then compare its discriminant with the Gram determinant.

\textbf{Step 1: the quadratic equation.} Put \(A_f=e^{4\pi i\alpha_f}\), \(B_j=e^{4\pi i\beta_j}\), \(z=e^{4\pi iu}\). Here \(\sum\alpha_f=\sum\beta_j=1/2\), so both products of exponentials equal one. Set

\[
R(z)=\prod_f(A_fz-1)-\prod_j(B_j-z)
=z(q_2z^2+q_1z+q_0).
\]

Write \(A=(A_f)_{f\in\mathcal F}\) and \(B=(B_j)_{j\in\{0\}\cup\mathcal Q}\) for the two quadruples of auxiliary exponentials. In terms of the degree-\(j\) elementary symmetric polynomials \(e_j\) in four variables,

\[
q_2=e_1(B)-e_3(A),\quad
q_1=e_2(A)-e_2(B),\quad
q_0=e_3(B)-e_1(A).
\]

\textbf{Step 2: its discriminant.} Expanding the symmetric polynomials and the determinant gives the Laurent-polynomial identity

\[
q_1^2-4q_2q_0=16\det G(x).
\tag{D.3}\label{eq-D-3}
\]

To expand the identity explicitly, write \(t_i=e^{\pi ix_i}\). The indicator \(1_{\ell\in f}\) is one when \(\ell\in f\) and zero otherwise, and likewise for a quadrilateral. The symbols \(A_f,B_j\) in this calculation are auxiliary exponentials; in particular their \(B_0\) is the \(j=0\) member of that quadruple. Then

\[
A_f=i\prod_{\ell=1}^6t_\ell^{1-2\,1_{\ell\in f}},\qquad
B_0=i\prod_\ell t_\ell^{-1},\qquad
B_g=i\prod_\ell t_\ell^{2\,1_{\ell\in g}-1},
\]

and \(c_i=(t_i^2+t_i^{-2})/2\). The three \(q\)'s have respectively eight, twelve, and eight Laurent monomials; substituting into \eqref{eq-D-3} and collecting coefficients gives zero residual.

\textbf{Step 3: the branch and the normalization.} Return to the real trigonometric products

\[
P(u)=\prod_f2\sin(2\pi(\alpha_f+u)),\qquad
Q(u)=\prod_j2\sin(2\pi(\beta_j-u)).
\]

Then \(Q-P=R(z)/z^2=q_2z+q_1+q_0/z\). At a zero,

\[
(Q-P)'(u)^2=-16\pi^2(q_1^2-4q_2q_0).
\]

At the real stationary point \((P-Q)'=32\pi p_*T_*>0\). Combining this with \eqref{eq-D-3} gives \(\sqrt{-\det G}=2p_*T_*\), proving \eqref{eq-D-2}. It also shows directly that the Gram determinant stays nonzero on this real domain. \(\square\)

This fixes the branch and factor used with the torsion formula of \cite{WY-T} in \hyperref[stmt-proposition-4-2]{Proposition 4.2}. At \(x=0\), it gives \(b_0(0)^2=1/32\), as in \hyperref[sec-D-1]{Section D.1}.

\section{Appendix E. One-edge central state sums}\label{sec-E}

This appendix gives an independent application of the reflection estimate; its results are not used in the proof of the \hyperref[stmt-main-theorem]{main theorem}. Here all six labels of a block are equal, and a fixed number of these block weights are multiplied together. The outer quantum dimension is a sine, hence a difference of two exponentials. We prove that the two corresponding sums have the same complete expansion at the central maximum. Their difference is therefore exponentially smaller than either summand.

\hyperref[stmt-theorem-e-2]{Theorem E.2} proves this cancellation for every fixed block count \(t\ge2\). Determining the smaller term requires more precision. \hyperref[sec-F]{Appendix F} supplies explicit constants, which allow us to evaluate it for every fixed \(t\ge8192\) in \hyperref[stmt-theorem-e-4]{Theorem E.4}. This number comes from comparing the remaining critical-value loss with the explicit reflection and cutoff gaps in \eqref{eq-E-19}; it is not a geometric transition or a claimed sharp threshold.

Throughout this appendix, the outer colors are restricted to the central set \(\mathcal C_r\) defined below, while each block contains its complete internal Racah sum. The results concern this central state sum. The estimates for outer colors omitted from \(\mathcal C_r\) would be a further requirement for a manifold application.

\subsection{E.1 Definition and normalization of the central sum}\label{sec-E-1}

In this section angular coordinates make the one-variable phase simpler. Let \(m\) be a nonnegative integer, so that the repeated even color is \(2m\), and put

\[\varepsilon=1/r,\qquad H=2\pi\varepsilon,\qquad
\sigma_r=(-1)^{(r+1)/2},\qquad X_m=(2m+1)H-\pi.\]

Thus the earlier small parameter \(h\) is \(\varepsilon\) in this section, and \(H=2\pi h\). For the repeated color \(a=2m\), the earlier centered variable is \(x=X_m/(2\pi)\). These are changes of units in the same block calculation. For odd \(r\ge9\) define

\[\mathcal C_r=\{m\in\mathbb Z_{\ge0}:3m\le r-2,\quad 3r<16m<5r\}.\]

The first condition in \(\mathcal C_r\) is the admissibility condition for three equal labels; the second selects the central range. Define the equal-label symbol and the outer sum by

\[
\begin{aligned}
J_r(2m)&=\frac{[m]!^6}{[3m+1]!^2}
\sum_{k=3m}^{\min(4m,r-2)}
\frac{(-1)^k[k+1]!}{[k-3m]!^4[4m-k]!^3},\\
S_{r,t}&=\sum_{m\in\mathcal C_r}[2m+1]J_r(2m)^t,
\qquad t\in\mathbb Z_{\ge2}.
\end{aligned}\tag{E.1}\label{eq-E-1}
\]

Here \([n]=\sin(nH)/\sin H\), as in \hyperref[sec-1-3-4]{Section 1.3.4}. All notation introduced in \hyperref[sec-E]{Appendices E}--\hyperref[sec-F]{F} is local to this one-edge application. The summation variable \(k\) in the first line is a local Racah index, and \(t\) is the fixed number of repeated block factors. With \(d=r-2-3m\), factorial reflection gives

\[
B_r(m):=\sigma_rJ_r(2m)
=\frac{[m]!^6[d]!^2}{|[r-1]!|}
\sum_{j=0}^{\min(m,d)}
\frac1{[j]!^4[m-j]!^3[d-j]!}>0.
\tag{E.2}\label{eq-E-2}
\]

Indeed, \([r-j]=-[j]\) and \([r-1-n]!=(-1)^n[r-1]!/[n]!\) give \eqref{eq-E-2}; all remaining factorial arguments lie below \(r/2\). In the convention of \hyperref[sec-1-3-5]{Section 1.3.5}, \(J_r(2m)=\tau_r(2m,\ldots,2m)\): the four triangle factors give \([m]!^6/[3m+1]!^2\), and \(i^{-12m}=1\). Write \(\mathcal Q_{h,\vartheta}\) for the normalized finite internal sum on the shifted lattice \(h(\mathbb Z+\vartheta)\), defined in \eqref{eq-B-5-1}; \(\vartheta\) is its real offset. Equation \eqref{eq-B-6-1} then identifies the normalization used here with that of the block integral:

\[D_r B_r(m)=\mathcal Q_{\varepsilon,\vartheta}(y\mathbf1),
\qquad y=X_m/(2\pi),\qquad \mathbf1=(1,\ldots,1),
\tag{E.3}\label{eq-E-3}\]

with the actual internal offset of that formula. The angular central interval tends to \([-\pi/4,\pi/4]\), a compact subset of the diagonal reflection domain \(|X|<\pi/3\).

For a direct comparison of kernels, let \(X\) be the continuous angular color coordinate and let \(b^2=H/\pi=2\varepsilon\) and

\[\mathsf S_H(z)=S_b(z/(\pi b)),\qquad
\mu=(\pi+X)/2,\qquad\nu=(\pi-3X)/2.\]

The shifts are \(\mathsf S_H(z+H)=2\sin z\,\mathsf S_H(z)\) and \(\mathsf S_H(z)\mathsf S_H(\pi+H-z)=1\). With angular internal variable \(v\), the one-variable kernel is

\[K_H(X,v)=
\frac{\mathsf S_H(\mu+H/2)^6\mathsf S_H(\nu+H/2)^2}
{\mathsf S_H(v+H/2)^4\mathsf S_H(\mu-v+H)^3\mathsf S_H(\nu-v+H)}.
\tag{E.4}\label{eq-E-4}\]

Since \(\mathsf S_H(H)=\sqrt{H/\pi}\), its finite sum is exactly

\[B_r(m)=\frac{H\sin H}{\pi}
\sum_{j=0}^{\min(m,d)}K_H(X_m,(j+1/2)H).\tag{E.5}\label{eq-E-5}\]

In the notation of \hyperref[sec-A-2]{Section A.2}, the factor-by-factor identity is

\[K_\varepsilon^+(y\mathbf1,u)
=K_H(X,\pi/4+2\pi u+H/4),\qquad y=X/(2\pi).
\tag{E.6}\label{eq-E-6}\]

For the actual index \(k=3m+j\), \eqref{eq-B-6-1} has \(u=j\varepsilon-1/8+\varepsilon/4\), giving \(v=(j+1/2)H\) in \eqref{eq-E-6}. Equations \eqref{eq-E-3} and \eqref{eq-E-6} fix the normalization, the integration variable, and the internal lattice offset. They allow us to apply the reflection and finite-sum estimates of \hyperref[sec-A]{Appendix A} and \hyperref[sec-B]{Appendix B} in the angular coordinates.

\subsection{E.2 The critical value and an even approximation}\label{sec-E-2}

Let \(F_{\mathrm{ang}}'(z)=\log(2\sin z)\), \(F_{\mathrm{ang}}(\pi/2)=0\), with the real branch on \((0,\pi)\). The phase of \eqref{eq-E-4} is

\[P(X,v)=6F_{\mathrm{ang}}(\mu)+2F_{\mathrm{ang}}(\nu)-4F_{\mathrm{ang}}(v)
-3F_{\mathrm{ang}}(\mu-v)-F_{\mathrm{ang}}(\nu-v).\tag{E.7}\label{eq-E-7}\]

It is strictly concave for real \(0<v<\min(\mu,\nu)\). At \(X=0\), its maximum is \(v_* =\pi/4\), with value \(v_8\) and \(P_{vv}=-8\). Its critical value is

\[f(X)=v_8-3\int_0^X L(s)\,ds,\qquad
L(z)=\operatorname{arcosh}\frac{\cos z}{2\cos z-1},\qquad L(z)=z+O(z^3).
\tag{E.8}\label{eq-E-8}\]

The branch specified here is odd, so \(f\) is even. It is strictly decreasing for \(0<X<\pi/3\). For example, putting \(k=\tan(X/2)\), the root of \(P_v=0\) near \(v=\pi/4\) is

\[\tan v_*(X)=
\frac{(1+k^2)\sqrt{1-2k^2}-4k^3}{1+3k^2+6k^4}.\]

Substitution and differentiation give \(dP(X,v_*)/dX=-3L(X)\), proving \eqref{eq-E-8}. The classical and leading-amplitude normalizations are

\[
\begin{aligned}
\phi(y\mathbf1)&=f(X)/(2\pi),\\
a_0(X)&=\sqrt{2/\pi}\,b_0(y\mathbf1)
=\frac{[(3\cos X-1)(1+\cos X)^3]^{-1/4}}{2\sqrt\pi}.
\end{aligned}\tag{E.9}\label{eq-E-9}
\]

The fourth root is continued from its positive value at zero. In particular, \(a_0(0)=1/(4\sqrt\pi)\).

\leavevmode\phantomsection\label{stmt-lemma-e-1}\textbf{Lemma E.1 (An even model for the repeated-color block).}\par\nopagebreak[4]

There are a fixed \(\rho>0\), a complex neighborhood of \([-\rho,\rho]\), and positive constants such that, for all sufficiently small allowed \(H\), an exactly even holomorphic function \(a_H\) satisfies

\[
\left|H^{-3/2}B_r(m)-a_H(X_m)e^{f(X_m)/H}\right|
\le C H^{-M}e^{(f(X_m)-\eta)/H},\qquad |X_m|\le\rho.
\tag{E.10}\label{eq-E-10}
\]

It has a uniform Gevrey-one expansion with leading term \(a_0\), and all its coefficients are even. We denote its coefficients in powers of \(H\) by \(a_j(X)\); their relation to the block coefficients \(b_j\) is \eqref{eq-E-12}. Thus \(a_H(X)\sim\sum_{j\ge0}a_j(X)H^j\). On the whole central range there is the absolute bound

\[B_r(m)\le C H^{-M}e^{f(X_m)/H}.\tag{E.11}\label{eq-E-11}\]

For \hyperref[stmt-theorem-e-4]{Theorem E.4} we use the explicit choices \(\rho=\delta=1/256\) and \(\eta=1/1024\). With these choices, \(a_H\) is holomorphic on \(|X|<1/64\) and \(a_H=a_0+O(H)\) uniformly on \(|X|\le1/128\). \hyperref[sec-F]{Appendix F} proves these constants. The level is still required to be sufficiently large; we do not give a numerical level threshold.

\textbf{Proof.} Apply the finite-lattice estimate of \hyperref[sec-B]{Appendix B} to \eqref{eq-E-3}, and optimally truncate the even coefficients from \hyperref[sec-A-4]{Sections A.4}--\hyperref[sec-A-5]{A.5} on the diagonal. The normalization \(D_r^{-1}=2\sin H\sqrt\varepsilon\) gives the formal relation

\[a_H(X)\sim\sqrt{2/\pi}\,\frac{\sin H}{H}
\sum_{n\ge0}b_n(y\mathbf1)\left(\frac{H}{2\pi}\right)^n.
\tag{E.12}\label{eq-E-12}\]

Thus \(a_1(0)/a_0(0)=13/48\). Independent reflection implies evenness on the diagonal. The uniform absolute bounds on compact real sets and the endpoint estimates in \hyperref[sec-B]{Appendix B} give \eqref{eq-E-11}, also at central colors approaching the ends of the interval. \hyperref[sec-F]{Appendix F} makes the remainder constants effective and supplies the stated \(\delta,\eta\). Thus the even approximation follows from the six-variable reflection theorem, with the change of scale recorded in \eqref{eq-E-12}. \(\square\)

\subsection{E.3 Cancellation of the two exponential terms}\label{sec-E-3}

Writing the sine quantum dimension as a difference of exponentials splits the outer sum into two terms. We call them the Weyl terms and denote them by \(I_{r,t}^{+}\) and \(I_{r,t}^{-}\). Their common expansion is given by the next theorem.

\leavevmode\phantomsection\label{stmt-theorem-e-2}\textbf{Theorem E.2 (Cancellation of the two Weyl terms).}\par\nopagebreak[4]

Define

\[I_{r,t}^{s}=\frac1{2i\sin H}\sum_{m\in\mathcal C_r}
e^{-isX_m}B_r(m)^t,\qquad s\in\{+1,-1\}.\]

Then

\[\sigma_r^tS_{r,t}=I_{r,t}^{+}-I_{r,t}^{-}.\tag{E.13}\label{eq-E-13}\]

For each fixed integer \(T\ge2\), both terms have the same complete Gevrey-one expansion, uniformly for \(2\le t\le T\) and both signs. More precisely, writing \(p_t=(3t-3)/2\), there are \(C,K,c,\kappa>0\) and coefficients \(c_n(t)\) such that, for integer truncation orders \(N\ge1\) in the indicated range,

\[
\left|H^{-p_t}e^{-tv_8/H}I_{r,t}^{s}
-\sum_{n=0}^{N-1}c_n(t)H^n\right|
\le C K^N N!H^N+C e^{-\kappa/H},\qquad NH\le c.
\tag{E.14}\label{eq-E-14}
\]

Their first two coefficients are

\[
\begin{aligned}
c_0(t)=C_t&=\frac{(4\sqrt\pi)^{-t}}{4i}\sqrt{\frac{2\pi}{3t}},\\
\frac{c_1(t)}{c_0(t)}=D_t&=\frac{13t^2+6t-13}{48t}.
\end{aligned}\tag{E.15}\label{eq-E-15}
\]

In particular, there is an \(\epsilon_T>0\) such that

\[|S_{r,t}|\le C_T\exp\!\left(\frac{tv_8-\epsilon_T}{H}\right),
\qquad 2\le t\le T,\quad r\text{ sufficiently large and odd}.
\tag{E.16}\label{eq-E-16}\]

\textbf{Proof.} Since \([2m+1]=-\sin X_m/\sin H\), the difference of the two exponentials gives \eqref{eq-E-13}. We will compare each finite sum with an integral. Reflection then identifies the two integrals, and the remainder estimate gives cancellation in the sums.

\textbf{Step 1: replace the finite sums by integrals.} The color lattice has spacing \(2H\) and an offset depending on \(r\) modulo four. Reflection \(X\mapsto-X\) need not preserve this lattice, which is why we first use the finite rectangle formula.

Choose a fixed \(0<d<\rho\) and use symmetric grid boundaries \(d_H=H\lfloor d/H\rfloor\). Equation \eqref{eq-E-11} and the strict decrease of \(f\) bound all terms outside \(|X_m|<d_H\) by an exponential smaller than \(e^{tv_8/H}\). In the interior, replace each block by \eqref{eq-E-10}. The product error has the same type of gap by \hyperref[stmt-proposition-a-4]{Proposition A.4}; \(t\) ranges over a fixed finite set and the dimensions and color count cost only powers of \(H\).

Apply the finite rectangle formula \eqref{eq-F-1} to the remaining lattice, of spacing \(2H\). On horizontal contours at a small fixed height \(s_0\), the nonconstant part of the cotangent kernel has decay \(e^{-\pi s_0/H}\), whereas the real phase can increase by at most \(C_Ts_0^2\). The vertical sides remain strictly below \(tv_8\) if \(s_0\) is small enough. Both offsets of the even-color lattice have a bounded cotangent kernel on these grid boundaries. Consequently, up to an error with a positive exponential gap, the corresponding term equals

\[\frac{H^{3t/2}}{4iH\sin H}
\int_{-d_H}^{d_H}e^{-isX}a_H(X)^t e^{tf(X)/H}\,dX.
\tag{E.17}\label{eq-E-17}\]

\textbf{Step 2: compute their common expansion.} The two integrals are exactly equal: the amplitude and phase are even, so their sine parts vanish and their cosine parts agree. Their common phase has \(f''(0)=-3\) and a unique real maximum. The uniform factorial remainders of \hyperref[sec-A-5]{Section A.5}, applied to this one-dimensional integral, give \eqref{eq-E-14}. Fixed endpoint changes have an exponential gap. Its Gaussian term gives \(C_t\).

For the first correction, \eqref{eq-E-12} gives \(a_1(0)/a_0(0)=13/48\), while

\[ (\log a_0)''(0)=3/4,\qquad f^{(4)}(0)=-15/2.\]

The contraction is therefore

\[\frac{13t}{48}+\frac{3t/4-1}{6t}
+\frac{t(-15/2)}{8(3t)^2}=D_t.\]

\textbf{Step 3: subtract the two terms.} Choose \(N\) proportional to \(H^{-1}\) in \eqref{eq-E-14} and subtract the two sums. The common truncated expansion cancels. Each remainder is exponentially small, which proves \eqref{eq-E-16} after decreasing the gap to absorb polynomial factors. \(\square\)

\leavevmode\phantomsection\label{stmt-corollary-e-3}\textbf{Corollary E.3 (Fixed even insertions).}\par\nopagebreak[4]

The cancellation bound \eqref{eq-E-16} remains valid after inserting an exactly even holomorphic factor \(\psi_H(X_m)\) into \eqref{eq-E-1}, provided these factors are uniformly bounded on a fixed complex neighborhood of the central interval. If \(\psi_H\) has a uniform Gevrey-one expansion, the two terms again have a common Gevrey-one expansion. For a fixed \(\psi\) its leading coefficient is \(\psi(0)C_t\).

\textbf{Proof.} The insertion preserves the equality of the two integrals in \eqref{eq-E-17}; its uniform bound preserves all error estimates. \(\square\)

\subsection{E.4 Evaluation of the smaller surviving term}\label{sec-E-4}

\leavevmode\phantomsection\label{stmt-theorem-e-4}\textbf{Theorem E.4 (Effective asymptotic for large fixed block count).}\par\nopagebreak[4]

For each fixed integer \(t\ge8192\), put

\[\theta_t=\frac{\pi}{3t},\qquad \ell_t=L(\theta_t),\qquad
V_t=t f(\theta_t),\qquad\Delta_t=tv_8-V_t,\]

and let

\[\vartheta(\ell)=\arccos\frac{\cosh\ell}{2\cosh\ell-1},
\qquad \vartheta(\ell)=\ell+O(\ell^3).\]

The displayed expansion selects the odd analytic branch at zero. This angle function is distinct from the lattice offset denoted by \(\vartheta\) in \eqref{eq-E-3}.

Then, along all odd levels,

\[
\begin{aligned}
S_{r,t}&=\sigma_r^{t+1}\mathcal A_t H^{(3t-3)/2}
e^{V_t/H}(1+O_t(H)),\\
\mathcal A_t&=\sinh\ell_t\,a_0(i\ell_t)^t
\sqrt{\frac{2\pi}{3t\vartheta'(\ell_t)}}>0.
\end{aligned}\tag{E.18}\label{eq-E-18}
\]

In particular \(H\log|S_{r,t}|\to V_t\). The integer \(t\) is fixed before \(r\to\infty\); no joint limit or uniform level threshold is asserted.

\textbf{Proof.} The proof separates the error comparison from the evaluation of the surviving term.

\textbf{Step 1: compare the error with the proposed leading term.} For \(0\le s\le1/4\), differentiation of \eqref{eq-E-8} gives \(s\le L(s)\le2s\). Thus

\[\Delta_t\le\frac{\pi^2}{3t},\qquad\ell_t\le\frac{2\pi}{3t}.\]

With the constants of \hyperref[stmt-lemma-e-1]{Lemma E.1}, for \(t\ge8192\) these imply

\[
\begin{aligned}
\Delta_t&<0.000401596<\eta/2,\qquad
\ell_t<0.000255664<\delta,\\
t(v_8-f(\delta))&\ge\tfrac32t\delta^2\ge3/16>\Delta_t.
\end{aligned}\tag{E.19}\label{eq-E-19}
\]

Therefore both the local replacement error and the remaining colors within \(\mathcal C_r\) are exponentially smaller than \(e^{V_t/H}\). Moving the cutoff to \(d_H=H\lfloor\delta/H\rfloor\) changes the phase at the cutoff by only \(O_t(H)\) and preserves these strict gaps. This is precisely the margin condition in \hyperref[stmt-proposition-a-4]{Proposition A.4}.

\textbf{Step 2: express the sum by residues.} The local model is even, and the sine multiplying it is odd. We exploit this parity in the residue formula for the finite color lattice. Its meromorphic kernel is

\[Q_r(z)=\frac\pi{2H}\cot\frac{\pi(z+\pi-H)}{2H}
=\frac\pi{2H}\left(\sigma_r\sec\frac{\pi z}{H}-\tan\frac{\pi z}{H}\right).
\tag{E.20}\label{eq-E-20}\]

Use residues on the rectangle with real endpoints \(\pm d_H\) and height \(\ell_t\). The factor \(-\sin z\,a_H(z)^te^{tf(z)/H}\) is odd. Its product with the tangent term is even and meromorphic, so that term has zero integral on this symmetric rectangle. Indeed the involution \(z\mapsto-z\) preserves the rectangle and changes the sign of the integral of an even meromorphic function. The secant term is therefore the only contribution. On the upper side,

\[\sec(\pi z/H)=2e^{i\pi z/H}
\bigl(1+O(e^{-2\pi\ell_t/H})\bigr).\]

\textbf{Step 3: locate the saddle and bound the rest of the contour.} The phase of the first secant term is \(\Phi_t(z)=tf(z)+i\pi z\). It satisfies

\[\Phi_t'(i\ell_t)=0,\qquad\Phi_t(i\ell_t)=V_t,
\qquad\Phi_t''(i\ell_t)=-3t\vartheta'(\ell_t)<0.
\tag{E.21}\label{eq-E-21}\]

Here \(L(i\ell)=i\vartheta(\ell)\) and \(L(\vartheta(\ell))=\ell\). Differentiating both sides proves the value identity

\[f(i\ell)=f(\vartheta(\ell))+3\ell\vartheta(\ell).\]

The whole rectangle lies in \(|z|<1/128\). On \(|z|\le1/64\) the formula

\[L'(z)=\frac{\sqrt{1+\cos z}}{(2\cos z-1)\sqrt{3\cos z-1}}\]

gives \(|L'(z)-1|<1/1000\), and hence \(\Re f''(z)<0\). For an explicit bound, regard the right side as a function of \(c=\cos z\). Use \(|\cos z-1|\le\cosh(1/64)-1<1/8191\) and \(|dL'/dc|<3\) on that \(c\) disk. Thus the upper horizontal phase has its unique real maximum at \(z=i\ell_t\) and a strict loss at both ends. On a vertical side, \(y\mapsto t\Re f(d_H+iy)-\pi y\) is convex for \(0\le y\le\ell_t\). At \(y=0\) it is below \(V_t\) by \eqref{eq-E-19}, and at \(y=\ell_t\) by horizontal strict concavity. Its maximum is therefore below \(V_t\). Moreover, \(d_H/H\) is an integer, so the secant on that side has absolute value \(\operatorname{sech}(\pi y/H)\le2e^{-\pi y/H}\). The lower half follows by reflection. All higher secant modes have a further positive exponential gap.

\textbf{Step 4: compute the coefficient.} The upper side is traversed from right to left. Combining its residue factor \(1/(2\pi i)\), the coefficient \(\sigma_r\pi/H\), and \(-\sin(i\ell_t)=-i\sinh\ell_t\) gives \(\sigma_r\sinh\ell_t/(2H)\). The lower side gives the same contribution. Multiplication by \(\sigma_r^tH^{3t/2}/\sin H\) and the Gaussian factor \(\sqrt{2\pi H/(3t\vartheta'(\ell_t))}\) yields exactly \eqref{eq-E-18}. The amplitude estimate \(a_H=a_0+O(H)\) is uniform on this contour. \(\square\)

The two theorems answer different questions about the central sum. \hyperref[stmt-theorem-e-2]{Theorem E.2} removes the expansion at the largest real critical value for every fixed \(t\ge2\). \hyperref[stmt-theorem-e-4]{Theorem E.4} evaluates the next contribution when the explicit error comparison \eqref{eq-E-19} holds. For smaller \(t\), determining that contribution requires sharper error bounds or a different contour analysis.

\section{Appendix F. Effective constants for the one-edge application}\label{sec-F}

This appendix proves the numerical choices in \hyperref[stmt-lemma-e-1]{Lemma E.1}. The reflection theorem already gives even coefficients; the task here is to keep an explicit exponential error after replacing the finite internal sum by its saddle expansion. We use the angular variables \(H,X,v\) and the phase \(P\) from \hyperref[sec-E]{Appendix E}.

The estimates proceed in the order in which they are used. \hyperref[sec-F-1]{Section F.1} gives the finite rectangle formula. \hyperref[sec-F-2]{Section F.2} constructs a Morse coordinate on specified disks. \hyperref[sec-F-3]{Sections F.3} and \hyperref[sec-F-4]{F.4} bound the factorial and Gaussian remainders. \hyperref[sec-F-5]{Section F.5} combines these bounds with the finite-sum comparison and proves the gap \(\eta=1/1024\).

The principal constants are collected here for reference:

{
\begin{longtable}[]{@{}
  >{\raggedright\arraybackslash}p{(\linewidth - 4\tabcolsep) * \real{0.3333}}
  >{\raggedright\arraybackslash}p{(\linewidth - 4\tabcolsep) * \real{0.3333}}
  >{\raggedright\arraybackslash}p{(\linewidth - 4\tabcolsep) * \real{0.3333}}@{}}
\toprule\noalign{}
\begin{minipage}[b]{\linewidth}\raggedright
Quantity
\end{minipage} & \begin{minipage}[b]{\linewidth}\raggedright
Choice
\end{minipage} & \begin{minipage}[b]{\linewidth}\raggedright
Used for
\end{minipage} \\
\midrule\noalign{}
\endhead
\bottomrule\noalign{}
\endlastfoot
External complex radius \(U_0\) & \(1/64\) & Holomorphy of the phase and amplitude \\
Morse radius \(R_w\) & \(1/8\) & Cauchy estimates for the transformed amplitude \\
Integration endpoint \(w_0\) & \(3/32\) & Gaussian truncation \\
Quadrature height \(s_q\) & \(1/256\) & Comparison with the internal finite sum \\
Truncation order \(N(H)\) & \(\lfloor1/(1024H)\rfloor\) & Exponential remainder \\
\end{longtable}
}

These choices are convenient sufficient bounds. The proof requires small enough \(H\), but does not specify a numerical lower bound for the level.

\subsection{F.1 A finite residue rectangle}\label{sec-F-1}

We first record the version of the finite-lattice identity used in \hyperref[sec-E]{Appendix E}. Let \(\Lambda=\xi+aH\mathbb Z\), where \(a>0\) is the spacing factor and \(\xi\) is a real offset, and choose real endpoints \(l<u\) off \(\Lambda\). Suppose \(g\) is holomorphic on a neighborhood of the rectangle with these endpoints and height \(s>0\). Set

\[q(z)=\frac\pi{aH}\cot\frac{\pi(z-\xi)}{aH}.\]

Let \(U\) be the upper path from \(u\) to \(l\), including the two upper vertical segments, and let \(L\) be the lower path from \(l\) to \(u\). Then residues and Cauchy's theorem give the exact identity

\[
\begin{aligned}
\sum_{x\in\Lambda\cap(l,u)}g(x)-\frac1{aH}\int_l^u g(x)\,dx
={}&\frac1{2\pi i}\int_U g(z)\left(q(z)+\frac{i\pi}{aH}\right)\,dz\\
&+\frac1{2\pi i}\int_L g(z)\left(q(z)-\frac{i\pi}{aH}\right)\,dz.
\end{aligned}\tag{F.1}\label{eq-F-1}
\]

On each horizontal side, the cotangent kernel differs from its limiting constant by at most

\[\frac{2\pi}{aH}\frac{e^{-2\pi s/(aH)}}{1-e^{-2\pi s/(aH)}}.\]

For the auxiliary lattice, \(a=1\), \(\xi=H/2\), and endpoints in \(H\mathbb Z\) give \(|q|\le\pi/H\) on the vertical sides. For either outer offset, \(a=2\), \(\xi=\pm H/2\) modulo \(2H\), the same choice gives \(|q|=\pi/(2H)\) on those sides. Thus if \(|g|\le C H^{-M}e^{B/H}\) on the horizontal sides and \(|g|\le C H^{-M}e^{E/H}\) on the vertical sides, the right side of \eqref{eq-F-1} is bounded by

\[C'H^{-M-1}\left(e^{(B-2\pi s/a)/H}+e^{E/H}\right).
\tag{F.2}\label{eq-F-2}\]

The contour lengths are bounded and \(s\) is fixed. The first term in \eqref{eq-F-2} is the horizontal error; the second accounts for the finite endpoints.

\subsection{F.2 A uniform one-variable Morse coordinate}\label{sec-F-2}

We need an analytic change of variable that makes \(P(X,v)\) quadratic in the internal variable, uniformly for complex \(X\) in a specified disk. We first locate the critical point and then bound the inverse coordinate. Use the polydisc

\[|X|\le U_0=1/64,\qquad |v-\pi/4|\le R=1/4.\]

Each of \(v,\mu-v,\nu-v\) is within \(R+3U_0/2\) of \(\pi/4\). Elementary trigonometric bounds on these disks give

\[|\csc^2 z|<9/2,\qquad |(\csc^2 z)'|<20.\]

At \(X=0\), \(P_{vv}(0,v)=-8/\cos(2(v-\pi/4))\), and the first \(X\) derivative of \(P_{vv}\) vanishes. Taylor's theorem therefore gives

\[
\begin{aligned}
|P_{vv}(X,v)-P_{vv}(0,v)|&\le30|X|^2,\\
|P_{vv}(X,v)+8|&\le
\frac{8(\cosh(2R)-1)}{2-\cosh(2R)}+30U_0^2<1.178<2.
\end{aligned}\tag{F.3}\label{eq-F-3}
\]

For the constant \(30\) in \eqref{eq-F-3}, differentiate twice in \(X\): the coefficients of the two derivatives of \(\csc^2\) have absolute values \(3/4\) and \(9/4\). Hence \(|\partial_X^2P_{vv}|\le(3/4+9/4)20=60\), and Taylor's theorem contributes the factor \(1/2\).

Also \(|P_v(X,\pi/4)|\le(15/2)|X|^2\). Here \(P_v(0,\pi/4)=\partial_XP_v(0,\pi/4)=0\), and \(|\partial_X^2P_v|\le(3/4+9/4)(9/2)<15\). On \(|v-\pi/4|=1/512\), the error relative to \(-8(v-\pi/4)\) is at most \(2/512+(15/2)U_0^2<8/512\). Rouché's theorem therefore gives a unique simple zero \(v_*(X)\) inside that circle. It is holomorphic in \(X\). Write \(\lambda=-P_{vv}(X,v_*)\), so \(|\lambda-8|<2\) and \(6<|\lambda|<10\).

\textbf{The change of variable.} For \(s=v-v_*(X)\), define the analytic Taylor factor

\[g(X,s)=\frac{2\bigl(f(X)-P(X,v_*(X)+s)\bigr)}{\lambda(X)s^2},
\qquad g(X,0)=1,\]

where the value at zero removes the apparent singularity. Since \(|P_{vvv}|\le36\), the local coordinate

\[w=(v-v_*)\sqrt\lambda\sqrt{g(X,v-v_*)},\qquad
P(X,v)=f(X)-w^2/2\]

satisfies \(|g(X,s)-1|\le2|s|\le1/4\) for \(|s|\le1/8\). Choose the roots continued from positive real values. Then \(|\sqrt g-1|<27/200\), and

\[\frac14+\frac{\sqrt{10}}8\frac{27}{200}<\frac{\sqrt6}{8}.
\tag{F.4}\label{eq-F-4}\]

Rouché's theorem on \(|s|=1/8\) gives a holomorphic inverse \(v(X,w)\) for \(|w|<1/4\), with \(|v-v_*|<1/8\). Applying Cauchy's estimate to \(v-v_*\) on disks of radius approaching \(1/8\) then gives

\[|v_w(X,w)|\le1\qquad (|w|\le1/8).\tag{F.5}\label{eq-F-5}\]

For real \(X,w\) the inverse is real and has positive derivative.

\subsection{F.3 Explicit bounds for the factorial expansion}\label{sec-F-3}

Let \(D_H(z)=\log\mathsf S_H(z+H/2)\). The midpoint and half-shift expansions are

\[D_H(z+cH)\sim\frac{F_{\mathrm{ang}}(z)}H+cF_{\mathrm{ang}}'(z)
+\sum_{j\ge1}\frac{B_{2j}(c+1/2)}{(2j)!}
F_{\mathrm{ang}}^{(2j)}(z)H^{2j-1},\qquad c\in\{0,1/2\}.
\tag{F.6}\label{eq-F-6}\]

The midpoint case follows from the integral and partial fractions used in \hyperref[stmt-lemma-a-1]{Lemma A.1}, after rescaling its variables. For \(c=1/2\), subtracting \(F_{\mathrm{ang}}'(z)/2\) leaves the corresponding kernel with \(\coth u\) in place of \(1/\sinh u\). Its partial-fraction remainder has the same absolute majorant. For a fixed strip margin \(0<a\le\pi/2\), after retaining the terms through \(j=J\ge1\), both cases have remainder at most

\[\frac{2\zeta(2J+2)}\pi
\left(1+\frac{a}{2\pi J}\right)(2J)!
\left(\frac{H}{2\pi a}\right)^{2J+1}
\quad\text{on }a\le\Re z\le\pi-a.
\tag{F.7}\label{eq-F-7}\]

For completeness, put \(p=2J+1\). The real parts of the scaled gamma arguments form two progressions, starting at \(\Re z\) and \(\pi-\Re z\), with spacing \(\pi\). Their reciprocal-power sum is bounded by

\[2\sum_{n\ge0}(a+n\pi)^{-p}
\le 2a^{-p}\left(1+\frac{a}{\pi(p-1)}\right).\]

The partial-fraction remainder for each gamma factor has coefficient \(2\zeta(p+1)(p-1)!H^p/(2\pi)^{p+1}\) times the corresponding reciprocal power; multiplication gives \eqref{eq-F-7}. The needed partial fractions and gamma expansions are recorded in \cite[Sections 4.36 and 5.11]{DLMF}. As in \hyperref[stmt-lemma-a-1]{Lemma A.1}, this is uniform in imaginary height: the oscillatory part of the integral is bounded before the geometric-series remainder is integrated. All unshifted arguments of \eqref{eq-E-4} in the preceding polydisc lie between \(1/2\) and \(\pi-1/2\) in real part, so we use \(a=1/2\).

On Cauchy disks of radius \(1/4\) about these arguments, \(|F_{\mathrm{ang}}'|<3\), and hence \(|F_{\mathrm{ang}}^{(n+1)}|\le3n!4^n\). Write

\[\log K_H=P/H+Q_0+q_H,\qquad
Q_0=-\tfrac32F_{\mathrm{ang}}'(\mu-v)-\tfrac12F_{\mathrm{ang}}'(\nu-v).\]

Denote the formal coefficient of \(H^n\) in \(q_H\) by \(q_n=q_n(X,v)\). There \(|e^{Q_0}|<2\). The even coefficients of \(q_H\) vanish; for odd \(n\) the Bernoulli bound gives

\[|q_n|\le\frac{192}{2\pi}\left(\frac2\pi\right)^n n!<32n!.\]

Applying \eqref{eq-F-7} to the sixteen factors, counted with absolute multiplicity, gives

\[\left|q_H-\sum_{n=1}^{N-1}q_nH^n\right|\le64N!H^N,
\qquad |q_H|\le64H,
\tag{F.8}\label{eq-F-8}\]

for sufficiently small positive \(H\). For \(N=1\), retain the first coefficient and use \eqref{eq-F-7} with \(J=1\); the general bound follows from the even and odd truncations of \eqref{eq-F-7}.

\textbf{Exponentiation.} We next pass from the logarithmic remainder to the amplitude itself, keeping a numerical factorial bound. Indeed, the finite polynomial \(q^{[N]}(z)=\sum_{n=1}^{N-1}q_nz^n\) has absolute value less than \(32\) for complex \(|z|\le1/(2N)\). Cauchy's estimate, \((2N)^N\le(2e)^NN!\), and \eqref{eq-F-8} imply that \(A_H=e^{-P/H}K_H\) has coefficients \(A_n\) satisfying

\[
\left|A_H-\sum_{n=0}^{N-1}A_nH^n\right|
\le C8^NN!H^N,\qquad NH\le1/4,
\qquad |A_n|\le C8^nn!.
\tag{F.9}\label{eq-F-9}
\]

Here Cauchy's estimate is applied to the finite polynomial \(q^{[N]}\). The prefactor may include \(e^{32}\); it is independent of \(H\) and \(N\) and does not affect the exponential rate. By \eqref{eq-F-5}, the same estimates hold for

\[\beta_H(X,w)=A_H(X,v(X,w))v_w(X,w)\]

on \(|w|\le R_w=1/8\).

\subsection{F.4 The Gaussian remainder}\label{sec-F-4}

Set \(w_0=3/32\) and

\[G_H(X)=\frac1{\sqrt{2\pi H}}\int_{-w_0}^{w_0}
e^{-w^2/(2H)}\beta_H(X,w)\,dw.\]

Let \(\beta_k(X,w)\) be the coefficient of \(H^k\) in the amplitude \(\beta_H(X,w)\) defined in \hyperref[sec-F-3]{Section F.3}. The saddle coefficients of \(G_H\) are

\[\gamma_n(X)=\sum_{k+j=n}
\frac{\partial_w^{2j}\beta_k(X,0)}{2^j j!}.\]

Cauchy's estimate and the Gaussian moments give

\[\left|\frac{\partial_w^{2j}\beta_k(X,0)}{2^j j!}\right|
\le C8^k k!128^j j!.\]

\textbf{Remainder estimate.} First truncate the \(H\) expansion with \eqref{eq-F-9}. In its coefficient of \(H^k\), Taylor-expand through degree \(2(N-k)-1\) in \(w\). Since \(w_0/R_w=3/4\), the remainder is at most four times the first omitted majorant term. Integrating with the full Gaussian moments and using \(\sum_{k=0}^n k!(n-k)!\le3n!\) gives \(C128^NN!H^N\).

Replacing the truncated moments by full moments costs \(e^{-w_0^2/(4H)}\). After extracting this factor, the remaining majorants are \(C(8H)^kk!(256H)^jj!\), whose sum over \(k+j<N\) is bounded if \(NH\le1/512\). Enlarging the factorial base gives

\[
\left|G_H-\sum_{n=0}^{N-1}\gamma_nH^n\right|
\le C256^NN!H^N+C e^{-9/(4096H)},\qquad NH\le1/512.
\tag{F.10}\label{eq-F-10}
\]

\textbf{Evenness of the coefficients.} The coefficients satisfy \(|\gamma_n|\le C256^nn!\) uniformly on the parameter disk. By the exact change of variables \eqref{eq-E-6}, their formal saddle series is the diagonal specialization of the series already proved even in \hyperref[sec-A-4]{Section A.4}, with an \(X\)-independent change of small parameter and measure. Thus every \(\gamma_n\) is even. The translation by \(H/4\) in \eqref{eq-E-6} does not change the formal saddle expansion: each positive-order translation term is an integral of a total derivative against the local Gaussian. Thus the parity of these coefficients follows from the already established block expansion.

\subsection{F.5 Comparison with the complete internal sum}\label{sec-F-5}

For real \(X\), let \(v_\pm(X)=v(X,\pm w_0)\). Then

\[P(X,v_\pm)=f(X)-d_0,\qquad d_0=w_0^2/2=9/2048.\]

We split the complete internal lattice at these two points. Strict real concavity bounds every term outside \([v_-,v_+]\) by an exponential with this decrease in the phase. The logarithmic bound for quantum factorials in \hyperref[sec-B]{Appendix B} includes endpoint factorials and costs only a power of \(H\). Inside this interval apply \eqref{eq-F-1} with spacing \(H\), vertical sides in \(H\mathbb Z\) within \(O(H)\) of \(v_\pm\), and height \(s_q=1/256\).

Since \(|P_{vv}|<10\) and the phase is real on the real axis,

\[\Re P(X,u+is_q)\le P(X,u)+5s_q^2.\]

The nonconstant part of the horizontal cotangent kernel has decay \(e^{-2\pi s_q/H}\); the real part of the phase on the vertical sides is at most \(f-d_0+5s_q^2+O(H)\). Both estimates are inside the pole-free polydisc of \hyperref[sec-F-2]{Section F.2}. The explicit inequalities

\[5s_q^2<d_0/2,\qquad2\pi s_q-5s_q^2>d_0/2\]

leave a fixed positive gap. The segments between the grid boundaries and \(v_\pm\) have the same decrease in the real part of the phase. Combining this with the exact prefactor \eqref{eq-E-5} gives

\[
\left|B_r(m)-\frac{\sin H}{\pi}\sqrt{2\pi H}
e^{f(X_m)/H}G_H(X_m)\right|
\le C H^{-M_0}e^{(f(X_m)-\eta_q)/H},
\quad \eta_q=9/8192.
\tag{F.11}\label{eq-F-11}
\]

The chosen \(\eta_q=d_0/4\) leaves room for the \(O(H)\) displacement and polynomial factors as \(H\to0\).

\textbf{Choice of the even approximation.} Set

\[N(H)=\left\lfloor\frac1{1024H}\right\rfloor,\qquad
a_H(X)=\frac{\sqrt{2\pi}}\pi\frac{\sin H}{H}
\sum_{n=0}^{N(H)-1}\gamma_n(X)H^n.
\tag{F.12}\label{eq-F-12}\]

This is exactly even and holomorphic on \(|X|<1/64\). The coefficient bounds give \(a_H=a_0+O(H)\) on \(|X|\le1/128\), with \(a_0\) as in \eqref{eq-E-9}. Moreover,

\[256^NN!H^N\le4^{-N}\le4e^{-(\log4)/(1024H)}.\]

Each available exponential gap exceeds \(\eta=1/1024\):

\[\frac{\log4}{1024}>\eta,\qquad \frac9{4096}>\eta,
\qquad\eta_q=\frac9{8192}>\eta.\tag{F.13}\label{eq-F-13}\]

There are three errors to combine: the factorial truncation, the omitted Gaussian tails, and the finite-lattice comparison. The three inequalities in \eqref{eq-F-13} put each of their exponential rates above \(\eta\). Equations \eqref{eq-F-10}--\eqref{eq-F-13} therefore prove \eqref{eq-E-10} with the constants stated in \hyperref[stmt-lemma-e-1]{Lemma E.1}.

\section{Acknowledgments}\label{acknowledgments}

The author is supported by BIMSA and the NSFC under Grant No. 12505090. The author thanks Shuang Ming, Tian Yang, and Yilong Wang for helpful discussions. The results were obtained with assistance from GPT-6 Astra. The author has checked the results and takes full responsibility for the content of this paper.

After the first version was posted, the author was informed by Yunpeng Meng, Baojun Wu, Tian Yang, and Yuxuan Yang that they had independently obtained the same results and were preparing an article on the subject.

\begingroup
\apptocmd{\thebibliography}{\setlength{\itemsep}{1pt}\setlength{\parsep}{0pt}\setlength{\parskip}{0pt}}{}{}
\interlinepenalty=10000
\bibliographystyle{amsplain-fullauthors}
\bibliography{references}
\endgroup

\end{document}